\documentclass[12pt]{article}

\usepackage{fullpage}
\usepackage{enumitem}
\usepackage[b]{esvect}
\usepackage{tikz}
\usepackage{caption}
\usepackage{subcaption}
\usepackage[font=footnotesize,labelfont=bf]{caption}
\usetikzlibrary{shapes.multipart,shapes.geometric,topaths,calc}
\usepackage{epsfig}
\usepackage{amsfonts}
\usepackage{amssymb}
\usepackage{amstext}
\usepackage{amsmath}
\usepackage{xspace}
\usepackage{theorem}
\usepackage{xcolor}
\usepackage{color}
\usepackage{graphicx}
\usepackage{ifthen}
\usepackage{hyperref}
\usepackage{tikz}
\usetikzlibrary{shapes, arrows, shapes.multipart,shapes.geometric,topaths,calc, positioning, decorations.pathreplacing, calligraphy}

\newcommand{\junk}[1]{}
\newcommand{\lcm}{\mathrm{lcm}}

\makeatletter
\newenvironment{proof}{{\bf Proof:  }}{\hfill\rule{2mm}{2mm}}

\newtheorem{theorem}{Theorem}

\newtheorem{proposition}[theorem]{Proposition}

\newcommand{\sq}{~\square~}
\newcommand{\F}{{\cal F}}

\title{Cycle decompositions of cartesian products of two cycles}

\author{
Moriah Aberle\thanks{University of Hawai'i at Manoa, \texttt{maberle@hawaii.edu}} 
\and Sarah Gold\thanks{Haverford College, \texttt{segold@haverford.edu}} 
\and Rivkah Moshe\thanks{Boston University, \texttt{rmoshe12@bu.edu}} 
\and David Offner\thanks{Carnegie Mellon University, \texttt{doffner@andrew.cmu.edu}}
}

\date{\today}

\begin{document}

\maketitle

\begin{abstract}
We say a graph $H$ decomposes a graph $G$ if there exists a partition of the edges of $G$ into subgraphs isomorphic to $H$.  We seek to characterize necessary and sufficient conditions for a cycle of length $k$, denoted $C_k$, to decompose the Cartesian product of two cycles $C_m ~\square~ C_n$.  We prove that if $m$ is a multiple of 3, then the Cartesian product of a cycle  $C_m$ and any other cycle can be decomposed into 3 cycles of equal length. This extends work of Kotzig, who proved in 1973 that the Cartesian product of two cycles can always be decomposed into two cycles of equal length. We also show that if $k$, $m$, and $n$ are positive, and $k$ divides $4mn$ then $C_{4k}$ decomposes $C_{4m} ~\square~ C_{4n}$. 
\end{abstract}

\section{Introduction}

Let $C_k$ represent the cycle of length $k$. We will use the convention that $V(C_k) = \{0,1, \ldots, k-1\}$ and $E(C_k) = \{\{i,i+1\} : 0 \le i \le k-1\}$, where addition is interpreted $\mod{k}$, and for convenience we will write edges $\{u,v\}$ as $uv$. The cartesian product of two graphs $G$ and $H$ is denoted as $G \sq H$ and is the graph where $V(G \sq H) = \{ (u,v): u \in V(G), v \in V(H)\}$ and $E(G \sq H) = \{(u, v)(u', v') : u = u', vv' \in E(H) \text{ or } v = v', uu' \in E(G)\}$.  Thus $V(C_m \sq C_n) = \{(i,j) : 0 \le i \le m-1, 0 \le j \le n-1\}$ and $E(C_m \sq C_n) = \{(i, j)(i+1,j) : 0 \le i \le m-1, 0 \le j \le n-1\} \cup \{(i, j)(i,j+1) : 0 \le i \le m-1, 0 \le j \le n-1\}$, where addition in the first coordinate is interpreted $\mod{m}$ and addition in the second coordinate is interpreted $\mod{n}$.  We refer to the cartesian product of two cycles as a torus.

\begin{figure}\begin{center}
\begin{tikzpicture}[scale=.75, transform shape]
%B1
    \draw[red, ultra thick, dashed] (-.5, -2) -- (0, -2) -- (0, 0) -- (1, 0) -- (1, -2) -- (2, -2); \draw[red, ultra thick, dashed] (-.5, -3) -- (2, -3);
    \draw[orange, ultra thick] (-.5, 0) -- (0, 0) -- (0, .5); \draw[orange, ultra thick] (2, 0) -- (1, 0) -- (1, .5); \draw[orange, ultra thick] (-.5, -1) -- (2, -1); \draw[orange, ultra thick] (0, -3.5) -- (0, -2) -- (1, -2) -- (1, -3.5);
%A1
    \draw[red, ultra thick, dashed] (2,.5) -- (2,-1) -- (3,-1) -- (3,0) -- (4,0) -- (4,.5);\draw[red, ultra thick, dashed] (2, -2) -- (4,-2) -- (4,-1) -- (5, -1);\draw[red, dashed, ultra thick] (2,-3) -- (2,-3) -- (2,-3.5);\draw[red, ultra thick, dashed] (4,-3.5) -- (4,-3) -- (4.5, -3);
    \draw[blue, densely dotted, ultra thick] (5, 0) -- (4,0) -- (4,-1) -- (3,-1) -- (3,-3) -- (4,-3) -- (4,-2) -- (5, -2);
    \draw[orange, ultra thick] (3, .5) -- (3,0) -- (2,0);\draw[orange, ultra thick] (2, -1) -- (2,-1) -- (2,-3) -- (3,-3) -- (3, -3.5);
%B2
    \draw[blue, ultra thick, densely dotted] (5, 0) -- (5, .5); \draw[blue, ultra thick, densely dotted] (6, .5) -- (6, 0) -- (7, 0); \draw[blue, ultra thick, densely dotted] (5, -2) -- (5, -1) -- (6, -1) -- (6, -2) -- (7, -2); \draw[blue, ultra thick, densely dotted] (5, -3.5) -- (5, -3) -- (6, -3) -- (6, -3.5); \draw[dashed, red, ultra thick] (5, -3) -- (5, -2) -- (6, -2) -- (6, -3) -- (7, -3);\draw[dashed, red, ultra thick] (5, -1) -- (5, 0) -- (6, 0) -- (6, -1) -- (7, -1);
%A2
    \draw[red, ultra thick, dashed] (7,.5) -- (7, 0) -- (8,0) -- (8, -2) -- (7,-2) -- (7,-1);\draw[red, ultra thick, dashed] (7,-3) -- (7,-3.5);
    \draw[blue, densely dotted, ultra thick] (7,0) -- (7,-1) -- (9,-1) -- (9,0) -- (10,0);\draw[blue, ultra thick, densely dotted] (7,-2) -- (7,-3) -- (9, -3) -- (9,-2) -- (10,-2);
    \draw[orange, ultra thick] (8,.5) -- (8,0) -- (9,0) -- (9,.5);\draw[orange, ultra thick] (8,-3.5) -- (8,-2) -- (9,-2) -- (9, -1) -- (10, -1);\draw[orange, ultra thick] (9, -3.5) -- (9,-3) -- (10, -3);
%B3
    \draw[blue, ultra thick, densely dotted] (10, 0) -- (10, .5); \draw[blue, ultra thick, densely dotted] (11, .5) -- (11, 0) -- (12, 0); \draw[blue, ultra thick, densely dotted] (10, -2) -- (10, -1) -- (11, -1) -- (11, -2) -- (12, -2); \draw[blue, ultra thick, densely dotted] (10, -3.5) -- (10, -3) -- (11, -3) -- (11, -3.5); \draw[orange, ultra thick] (10, -3) -- (10, -2) -- (11, -2) -- (11, -3) -- (12, -3);\draw[orange, ultra thick] (10, -1) -- (10, 0) -- (11, 0) -- (11, -1) -- (12, -1);
%A3
    \draw[red, ultra thick, dashed] (13, .5) -- (13,0) -- (14,0) -- (14,-2) -- (14.5,-2);\draw[red, dashed, ultra thick] (14.5,-3) -- (13,-3) -- (13,-3.5);
    \draw[blue, densely dotted, ultra thick] (12,0) -- (12,-1) -- (13,-1) -- (13, -3) -- (12,-3) -- (12,-2);
    \draw[orange, ultra thick] (12, .5) -- (12, 0) -- (13,0) -- (13,-1) -- (14,-1) -- (14.5, -1);\draw[orange, ultra thick] (14,.5) -- (14,0) -- (14.5, 0);\draw[orange, ultra thick] (12, -1) -- (12, -2) -- (14, -2) -- (14, -3.5);\draw[orange, ultra thick] (12,-3) -- (12, -3.5);
\foreach \x in {0,1,...,14}{\node at (\x,1){$\x$};}\foreach \y in {0,1,2,3}{\node at (-1,-\y){$\y$};}
\foreach \x in {0,1,...,14}{\foreach\y in {0,-1,-2,-3}{\filldraw[fill=white,draw=black] (\x,\y) circle (3pt);}}
\end{tikzpicture}\end{center}
\caption{A decomposition of $C_4\sq C_{15}$ into three cycles (dashed red, dotted blue, and solid yellow), each of length 40.}\label{Ex1}
\end{figure}
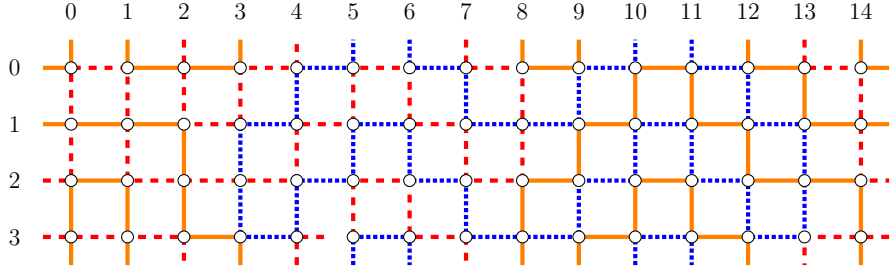

When drawing cartesian products of cycles $C_m \sq C_n$, we adopt the convention that each copy of $C_m$ is drawn vertically, and each copy of $C_n$ is drawn horizontally, with $(0,0)$ in the upper left.  For $0 \le i \le m-1$, the edge $(i,0)(i,n-1)$ is represented by a half-edge to the left of $(i,0)$ and a half-edge to the right of $(i,n-1)$.  Similarly, for $0 \le j \le n-1$, the edge $(0,j)(m-1,j)$ is represented by a half-edge above $(0,j)$ and a half-edge below $(m-1,j)$. See Figures~\ref{Ex1} and \ref{Ex2}.

%% Example: 4x15 

Given two graphs $G$ and $H$, define the union $G \cup H$ to be the graph where $V(G \cup H) = V(G) \cup V(H)$ and $E(G \cup H)= E(G) \cup E(H)$. We use the notation $G \sqcup H$ to denote an edge-disjoint union of graphs, that is, a union of graphs $G$ and $H$ where $E(G) \cap E(H) = \emptyset$. A decomposition of a graph $G$ is a set of subgraphs that partition the edges in $G$. That is, a set of graphs $\{G_1, G_2, ..., G_m \}$ is a decomposition of $G$ if and only if $E(G) = \bigsqcup\limits_{1 \leq i \leq m} E(G_i)$. Note $ E(G_i) \cap E(G_j) = \emptyset$ if $i \neq j$. If each graph in a decomposition is isomorphic to a single graph $H$, we call it an $H$-decomposition, and we say $H$ decomposes $G$.  In this paper we are interested in $C_k$-decompositions of cartesian products $C_m \sq C_n$.  Throughout the paper we assume $m$ and $n$ are both at least three.

Our goal is to determine necessary and sufficient conditions for cycle decompositions of cartesian products of cycles.  Since $C_m \sq C_n$ has $mn$ vertices and $2mn$ edges, if $C_k$ decomposes $C_m \sq C_n$, then it must be the case that $k \le mn$ and $k$ divides $2mn$.  We call the necessary condition that $k$ divides $2mn$ the divisibility condition.

There is a vast literature on cycle decompositions of graphs, particularly into Hamiltonian cycles, much of it summarized in \cite{alspach}. The most general result in this area on cartesian product graphs is due to  Stong~\cite{stong}, who proved that if $G$ and $H$ both have Hamiltonian decompositions, and satisfy some mild additional assumptions, then $G \sq H$ also has a Hamiltonian decomposition.    There has also been interest in decomposition of cartesian products of graphs into smaller cycles, particularly in the case of hypercube graphs.  See, for example, \cite{AOT21}, \cite{HSW}, \cite{MR}, \cite{TWB}.

Since the number of edges in $C_m \sq C_n$ is even, a Hamiltonian cycle always meets the divisibility condition, and in 1973, Kotzig~\cite{kotzig} proved that every cartesian product of two cycles can be decomposed into two (Hamiltonian) cycles of equal length.

In Section~\ref{c3m} we prove a theorem on decomposition into three cycles analogous to Kotzig's result on two cycles.  The divisibility condition implies that if $C_m \sq C_n$ can be decomposed into three cycles, then at least one of $m$ and $n$ must be a multiple of three.  We prove this condition is also sufficient.

\begin{theorem}\label{3cycles}
 If $m \ge 3$ and $n$ is a positive multiple of 3, then $C_m \sq C_n$ can be decomposed into 3 cycles of equal length.
\end{theorem}

In 2022, Gibson and Offner~\cite{GO22} proved that if $\ell \ge 1$, $n \ge 1$, and $n$ divides $4\ell$, then  $C_{4\ell} \sq C_4$ has a decomposition into cycles of length $4n$ . In Section~\ref{mult4} we generalize their technique to show that cartesian products of cycles whose lengths are multiples of four can be decomposed into a wide range of cycle lengths divisible by four.

\begin{theorem}\label{4m4n}
If $k$, $m$, and $n$ are positive integers and $k$ divides $4mn$, then $C_{4k}$ decomposes $C_{4m} \sq C_{4n}$.
\end{theorem}

This resolves the question of decomposability for all possible cycles $C_{4k}$ that meet the divisibility condition for $C_{4m} \sq C_{4n}$ except where $k$ divides $8mn$ but not $4mn$. For example, if $k=8$, $m=1$, and $n=5$ then the divisibility condition does not preclude a $C_{4k} = C_{32}$-decomposition of $C_{4m} \sq C_{4n} = C_{4} \sq C_{20}$, but this is not produced by our method.  Of course, it also leaves open the question of when $C_k$ decomposes $C_{4m} \sq C_{4n}$ when $k$ is congruent to $2 \mod{4}$. We note that it is not the case that all cycles $C_k$ where $k \le 16mn$ and $k$ divides $2(4m)(4n)= 32mn$ will decompose $C_{4m} \sq C_{4n}$. For example, Proposition~\ref{c6thm} implies that $C_{4m} \sq C_{4n}$ has no $C_6$-decomposition for any choices of $m$ and $n$.

In Section~\ref{odd} we present some other positive and negative results on decompositions into cycles of length six and cycles with odd lengths. 
The general case of necessary and sufficient conditions for $C_k$ to decompose $C_m \sq C_n$  remains open.

\section{Decomposition of $C_m \sq C_n$ into three cycles}\label{c3m}

%TODO: replace orange with black in 3 cycle figures?

%Example: 10x18
\begin{figure}\begin{center}
\begin{tikzpicture}[scale=.65,transform shape]
%A
    \draw[blue, ultra thick, densely dotted] (-.5, 0) -- (2, 0) -- (2, -1) -- (5, -1) -- (5, -2) -- (6, -2); \draw[blue, ultra thick, densely dotted] (-.5, -2) -- (0, -2) -- (0, -3) -- (5, -3) -- (5, -4) -- (6, -4); \draw[blue, ultra thick, densely dotted] (5, .5) -- (5, 0) -- (6, 0); \draw[blue, ultra thick, densely dotted] (5.5, -5) -- (5, -5) -- (5, -6); \draw[blue, ultra thick, densely dotted] (-.5, -4) -- (2, -4) -- (2, -5) -- (-.5, -5);\draw[orange, ultra thick] (0, .5) -- (0, -1) -- (1, -1) -- (1, .5); \draw[orange, ultra thick] (3, .5) -- (3, 0) -- (4, 0) -- (4, -2) -- (5, -2) -- (5, -3) -- (5.5, -3); \draw[orange, ultra thick] (-.5, -3) -- (0, -3) -- (0, -5.5); \draw[orange, ultra thick] (1, -5.5)  -- (1, -2) -- (3, -2) -- (3, -4) -- (4, -4) -- (4, -5) -- (3, -5) -- (3, -5.5);\draw[red, ultra thick, dashed] (-.5, -1) -- (0, -1) -- (0, -2) -- (1, -2) -- (1, -1) -- (2, -1) -- (2, -4) -- (3, -4) -- (3, -5) -- (2, -5) -- (2, -5); \draw[dashed, red, ultra thick] (2, .5) -- (2, 0) -- (3, 0) -- (3, -2) -- (4, -2) -- (4, -4) -- (5, -4) -- (5, -5) -- (4, -5) -- (4, -5); \draw[dashed, red, ultra thick] (4, .5) -- (4, 0) -- (5, 0) -- (5, -1) -- (6, -1);
%B(0)
    \draw[blue, ultra thick, densely dotted] (6, 0) -- (12, 0); \draw[blue, ultra thick, densely dotted] (6, -2) -- (12, -2); \draw[blue, ultra thick, densely dotted] (6, -4) -- (12, -4); \draw[blue, ultra thick, densely dotted] (5.55, -5) -- (12, -5); \draw[red, ultra thick, dashed] (6, -1) -- (6, .5); \draw[dashed, red, ultra thick] (8, .5) -- (8, -1) -- (7, -1) -- (7, -3) -- (6, -3) -- (6, -5); \draw[red, ultra thick, dashed] (10, .5) -- (10, -1) -- (9, -1) -- (9, -3) -- (8, -3) -- (8, -5); \draw[red, ultra thick, dashed] (12, -1) -- (11, -1) -- (11, -3) -- (10, -3) -- (10, -5);\draw[orange, ultra thick] (7, .5) -- (7, -1) -- (6, -1) -- (6, -3) -- (5.4, -3); \draw[orange, ultra thick] (9, .5) -- (9, -1) -- (8, -1) -- (8, -3) -- (7, -3) -- (7, -5.5); \draw[orange, ultra thick] (11, .5) -- (11, -1) -- (10, -1) -- (10, -3) -- (9, -3) -- (9, -5.5); \draw[orange, ultra thick] (11.6, -3) -- (11, -3) -- (11, -5.5);
%B(1)
    \draw[blue, ultra thick, densely dotted] (12, 0) -- (17.5, 0); \draw[blue, ultra thick, densely dotted] (12, -2) -- (17.5, -2); \draw[blue, ultra thick, densely dotted] (12, -4) -- (17.5, -4); \draw[blue, ultra thick, densely dotted] (12, -5) -- (17.5, -5);
    \draw[red, ultra thick, dashed] (12, -1) -- (12, .5); \draw[dashed, red, ultra thick] (14, .5) -- (14, -1) -- (13, -1) -- (13, -3) -- (12, -3) -- (12, -5); \draw[red, ultra thick, dashed] (16, .5) -- (16, -1) -- (15, -1) -- (15, -3) -- (14, -3) -- (14, -9.5); \draw[red, ultra thick, dashed] (17.5, -1) -- (17, -1) -- (17, -3) -- (16, -3) -- (16, -9.5);
    \draw[orange, ultra thick] (13, .5) -- (13, -1) -- (12, -1) -- (12, -3) -- (11.4, -3); \draw[orange, ultra thick] (15, .5) -- (15, -1) -- (14, -1) -- (14, -3) -- (13, -3) -- (13, -9.5); \draw[orange, ultra thick] (17, .5) -- (17, -1) -- (16, -1) -- (16, -3) -- (15, -3) -- (15, -9.5); \draw[orange, ultra thick] (17.5, -3) -- (17, -3) -- (17, -9.5);
%C(0)
    \draw[blue, ultra thick, densely dotted] (-.5, -6) -- (5, -6); \draw[blue, ultra thick, densely dotted] (-.5, -7) -- (5, -7); 
    \draw[blue, ultra thick, densely dotted] (17.5,-6) -- (12,-6) -- (12,-7) -- (17.5,-7);
    \draw[orange, ultra thick] (0, -5.4) -- (0, -7.5); \draw[orange, ultra thick] (1, -5.4) -- (1, -7.5); \draw[orange, ultra thick] (3, -5.4) -- (3, -7.5); 
    \draw[orange, ultra thick] (7, -5.4) -- (7, -6) -- (6, -6) -- (6, -7) -- (7, -7) -- (7, -7.5); \draw[orange, ultra thick] (9, -5.4) -- (9, -6) -- (8, -6) -- (8, -7) -- (9, -7) -- (9, -7.5); \draw[orange, ultra thick] (11, -5.4) -- (11, -6) -- (10, -6) -- (10, -7) -- (11, -7) -- (11, -7.5); 
    \draw[red, ultra thick, dashed] (2, -5) -- (2, -7);\draw[red, ultra thick, dashed] (4, -5) -- (4, -7); \draw[red, ultra thick, dashed] (6, -5) -- (6, -6) -- (5, -6) -- (5, -7) -- (6, -7) -- (6, -7); \draw[red, dashed, ultra thick] (12,-5) -- (12, -6) -- (11, -6) -- (11, -7) -- (12, -7) -- (12,-7);
    \draw[red, ultra thick, dashed] (8, -5) -- (8, -6) -- (7, -6) -- (7, -7) -- (8, -7) -- (8, -7);\draw[red, ultra thick, dashed] (10, -5) -- (10, -6) -- (9, -6) -- (9, -7) -- (10, -7) -- (10, -7);
%C(1)
    \draw[blue, ultra thick, densely dotted] (-.5, -8) -- (5, -8) -- (5, -7); \draw[blue, ultra thick, densely dotted] (-.5, -9) -- (5, -9) -- (5, -9.5); 
    \draw[blue, ultra thick, densely dotted] (17.5,-8) -- (12,-8) -- (12,-9) -- (17.5,-9);
    \draw[orange, ultra thick] (0, -7.4) -- (0, -9.5); \draw[orange, ultra thick] (1, -7.4) -- (1, -9.5); \draw[orange, ultra thick] (3, -7.4) -- (3, -9.5); 
    \draw[orange, ultra thick] (7, -7.4) -- (7, -8) -- (6, -8) -- (6, -9) -- (7, -9) -- (7, -9.5); \draw[orange, ultra thick] (9, -7.4) -- (9, -8) -- (8, -8) -- (8, -9) -- (9, -9) -- (9, -9.5); \draw[orange, ultra thick] (11, -7.4) -- (11, -8) -- (10, -8) -- (10, -9) -- (11, -9) -- (11, -9.5); 
    \draw[red, ultra thick, dashed] (2, -7) -- (2, -9.5);\draw[red, ultra thick, dashed] (4, -7) -- (4, -9.5); \draw[red, ultra thick, dashed] (6, -7) -- (6, -8) -- (5, -8) -- (5, -9) -- (6, -9) -- (6, -9.5); \draw[red, dashed, ultra thick] (12,-7) -- (12, -8) -- (11, -8) -- (11, -9) -- (12, -9) -- (12,-9.5);
    \draw[red, ultra thick, dashed] (8, -7) -- (8, -8) -- (7, -8) -- (7, -9) -- (8, -9) -- (8, -9.5);\draw[red, ultra thick, dashed] (10, -7) -- (10, -8) -- (9, -8) -- (9, -9) -- (10, -9) -- (10, -9.5);

    \foreach \x in {0,1,...,17}{\node at (\x,1){$\x$}; \foreach \y in {0,1,...,9}{\filldraw[fill=white, draw=black] (\x,-\y) circle (3pt);}}
    \foreach \y in {0,1,...,9}{\node at (-1,-\y){$\y$};}

\end{tikzpicture}
\end{center}
\caption{A decomposition of $C_{10} \sq C_{18}$ into three cycles, each of length 120.}\label{Ex2}
\end{figure}
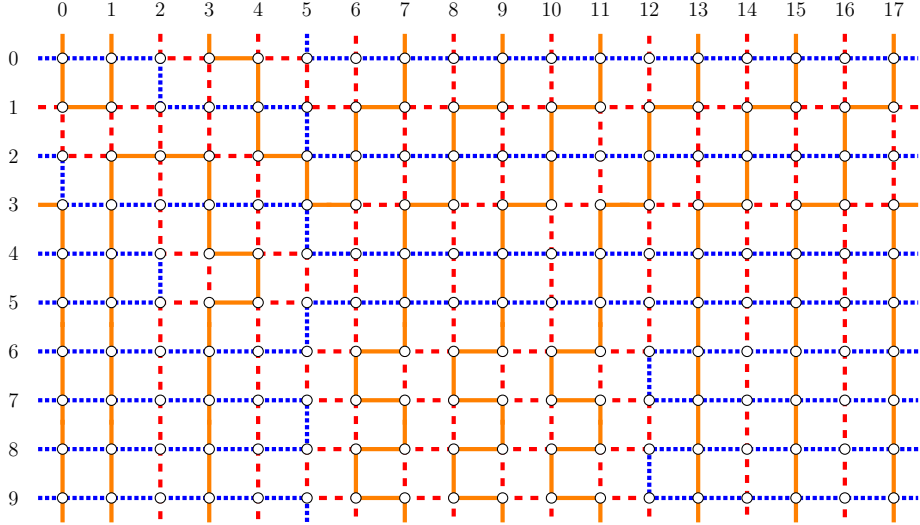

\begin{proof}[Proof of Theorem~\ref{3cycles}]
Let $m$ and $n$ be positive integers where $m \ge 3$ and $n$ is a multiple of 3. 

This proof is divided into seven cases:
\begin{enumerate}
    \item $m$ is odd,
    \item $m=4$ and $n$ is even, 
    \item $m=4$ and $n=3$, 
    \item $m=4$ and $n\ge 9$ is odd,
    \item $m\geq 6$ is even and $n = 6$,
    \item $m \geq 6$ is even and $n \ge 12$ is even, 
    \item $m \geq 6$ is even and $n$ is odd.
\end{enumerate}

The proof strategy in all cases is inductive:  We divide the edges in $C_m \sq C_n$ into blocks, where we use the term block to denote a set of edges incident with a set of vertices $(i,j)$ where $a \le i \le b$ and $c \le j \le d$ for some integers $a,b,c,d$ (These form rectangular blocks in the representations of the cartesian products in the figures).  We refer to such a block as a $(b-a+1) \times (d-c+1)$ block, and use different letters to distinguish different types of block, calling them $A$ blocks, $B$ blocks, and so forth.  We refer to  $C_m \sq C_n$ for the smallest values of $m$ and $n$ in a given case as a minimal block, and in all cases show how to partition the edges of this minimal block into three cycles of the same length, giving the base case of the induction.  In all cases we will describe one cycle as dashed red, one as solid yellow, and one as dotted blue.

For the inductive step, we suppose $m$ or $n$ are larger than in the minimal block, so we divide the remaining edges among some other blocks and show how to partition these edges into three colors, where the edges of a given color in each block form a set of paths, and there are the same number of edges of each color.  We ensure that edges that are in multiple blocks are colored with the same color, so the paths in the other blocks extend the cycles in the minimal block and preserve the fact that in $C_m \sq C_n$ all edges are contained in one of three equal length cycles. 

In each case, we describe the block structure of the decomposition, and explicitly list the edges of each color in each block. Throughout the proof, addition in the first coordinate is done modulo $m$ and addition in the second coordinate is done modulo $n$.

\paragraph{Case 1 ($m$ is odd):} Assume $m = 3+2M$ for some nonnegative integer $M$ and $n=3N$ for some positive integer $N$. The decomposition of $C_m \sq C_n$ is illustrated in Figure~\ref{ABcase1} (middle) and contains $N$ $3 \times 3$ $A$ blocks and $MN$ $2 \times 3$ $B$ blocks, which are defined below.  The $A$ block is the minimal block for this case.

For the $3 \times 3$ $A$ block define three edge sets $A^r$, $A^y$, and $A^b$, each of which forms a cycle with six edges, shown as dashed red, solid yellow, and dotted blue edges, respectively, in Figure~\ref{ABcase1} (top left). 

$A^r$ is the set of edges $\{(-1,0)(0, 0), (0,0)(0,1), (0,1)(1, 1), (1,1)(1,2), (1,2)(2,2),\\ (2,2)(2,3), (2,-1)(2, 0), (2,0)(3,0)\}$.

$A^y$ is the set of edges $\{(0,-1)(0,0), (0,0)(1,0), (1, 0)(1, 1), (1,1)(2, 1), (2, 1)(2, 2), \\(2, 2)(3,2), (-1,2)(0,2), (0,2)(0,3)\}.$

$A^b$ is the set of edges $\{(-1, 1)(0,1), (0,1)(0,2), (0, 2)(1,2), (1, 2)(1, 3), (1, -1)(1, 0), \\(1, 0)(2,0), (2,0)(2,1), (2,1)(3,1)\}.$

For $0 \le J \le N-1$ let the block $A(J)$ in $C_m \sq C_n$ be the set of edges incident to the vertices $(i,j)$ where $0 \le i \le 2$ and $3J \le j \le 2+3J$.  To color the edges in $A(J)$, define three edge sets $A^r(J)$, $A^y(J)$, and $A^b(J)$, where for $\ell \in \{r,y,b\}$, edge $(i,j)(i',j')$ is in $A^\ell$ if and only if $(i,j+3J)(i',j'+3J)$ is in $A^\ell(J)$. For example, the edge $(2,7)(2,8)$ is yellow (in $A^y(2)$) because $(2,1)(2,2)$ is yellow (in $A^y$).

For the $2 \times 3$ $B$ block define three edge sets $B^r$, $B^y$, and $B^b$, each containing four edges and shown as dashed red, solid yellow, and dotted blue, respectively, in Figure~\ref{ABcase1} (top right).

$B^r = \{(0, -1)(0, 0), (0,0)(1,0), (1, 0)(1, 1), (1, 1)(0, 1), (0, 1)(0,2)\}$

$B^y=\{(2, -1)(2, 0), (2, 0)(3, 0), (-1, 0)(0, 0), (0, 0)(0, 1), (0, 1)(-1, 1), (3, 1)(2, 1), (2, 1)(2, 2)\}$ 

$B^b =\{(1, -1)(1, 0), (1, 0)(2, 0), (2, 0)(2, 1), (2, 1)(1, 1), (1, 1)(1, 2)\}$

For $0 \le I \le M-1$ and $0 \le J \le N-1$, let the block $B(I,J)$ in $C_m \sq C_n$ be the set of edges incident with $(i,j)$ where $3+2I \le i \le 4+2I$ and $3J \le j \le 2+3J$. To color the edges in $B(I,J)$, define three edge sets $B^r(I,J)$, $B^y(I,J)$, and $B^b(I,J)$, where for $\ell \in \{r,y,b\}$, edge $(i,j)(i',j')$ is in $B^\ell$ if and only if $(i + 3+2I,j+3J)(i'+3+2I,j'+3J)$ is in $B^\ell(I,J)$. For example, $(7, 4)(8,4)$ is red (in $B^r(2,1)$) because $(0,1)(1,1)$ is red (in $B^r$).

For $\ell \in \{r,y,b\}$ each cycle in the decomposition is the union of the $A^\ell(J)$ and $B^\ell(I,J)$ sets, i.e.

\[ \left(\bigcup_{J=0}^{N-1}A^\ell(J) \right) \cup
\left(\bigcup_{I=0}^{M-1}\bigcup_{J=0}^{N-1}B^\ell(I,J) \right).\]

Since each of the $N$ $A$ blocks contains six edges of each color, and each of the $MN$ $B$ blocks contains four edges of each color, each cycle contains $6N + 4MN$ edges.  Since $m=3+2M$ and $n=3N$, the total number of edges in $C_m \sq C_n$ is $2mn= 2(3+2M)(3N) = 18N + 12MN$.  Thus each cycle contains one third of the total edges, as desired.

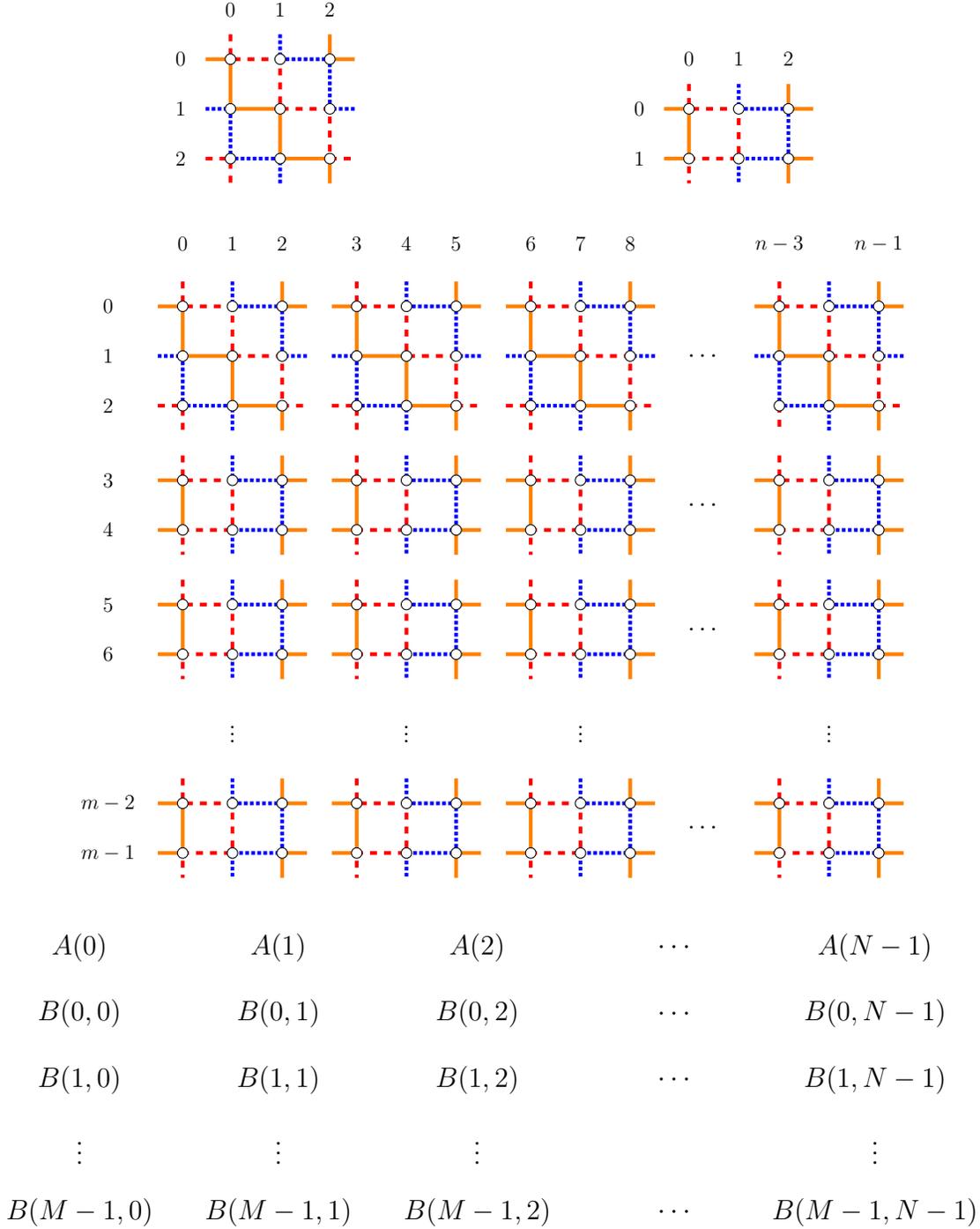
\begin{figure}
\begin{center}
%% A Block
\begin{tikzpicture}[scale=.75, transform shape]
\draw[red, ultra thick, dashed] (-.5, -2) -- (0, -2) -- (0, -2.5); \draw[blue, ultra thick, densely dotted] (-.5, -1) -- (0, -1) -- (0, -2) -- (1, -2) -- (1, -2.5); \draw[orange, ultra thick] (-.5, 0) -- (0, 0) -- (0, -1) -- (1, -1) -- (1, -2) -- (2, -2) -- (2, -2.5); \draw[red, ultra thick, dashed] (0, .5) -- (0, 0) -- (1, 0) -- (1, -1) -- (2, -1) -- (2, -2) -- (2.5, -2); \draw[blue, ultra thick, densely dotted] (1, .5) -- (1, 0) -- (2, 0) -- (2, -1) -- (2.5, -1); \draw[orange, ultra thick] (2, .5) -- (2, 0) -- (2.5, 0); \foreach \x in {0, 1, 2}{\foreach \y in {0, -1, -2}{\draw[draw = black, fill = white] (\x, \y) circle (3pt);}}
\foreach \x in {0,1,2}{\node at (\x,1){$\x$};}\foreach \y\i in {0/0,-1/1,-2/2}{\node at (-1,\y){$\i$};}
\end{tikzpicture}
\hspace{1.5in}
%B Block
\begin{tikzpicture}[scale=.75, transform shape]
\draw[red, ultra thick, dashed] (0, -3) -- (0, -3.5) -- (1, -3.5) -- (1, -4.5) -- (0, -4.5) -- (0, -5); \draw[blue, ultra thick, densely dotted] (1, -3) -- (1, -3.5) -- (2, -3.5) -- (2, -4.5) -- (1, -4.5) -- (1, -5); \draw[orange, ultra thick] (-.5, -3.5) -- (0, -3.5) -- (0, -4.5) -- (-.5, -4.5); \draw[orange, ultra thick] (2, -3) -- (2, -3.5) -- (2.5, -3.5); \draw[orange, ultra thick] (2.5, -4.5) -- (2, -4.5) -- (2, -5); \foreach \x in {0, 1, 2}{\foreach \y in {-3.5, -4.5}{\draw[draw = black, fill = white] (\x, \y) circle (3pt);}}
\foreach \x in {0,1,2}{\node at (\x,-2.5){$\x$};}\foreach \y\i in {-3.5/0,-4.5/1}{\node at (-1,\y){$\i$};}
\end{tikzpicture}

\vspace{.25in}

\begin{tikzpicture}[scale=.75, transform shape]
% Indices
\node at (0, 1.25){0};\node at (1, 1.25){1};\node at (2, 1.25){2};\node at (3.5, 1.25){3};\node at (4.5, 1.25){4};\node at (5.5, 1.25){5};\node at (7, 1.25){6};\node at (8, 1.25){7};\node at (9, 1.25){8};\node at (12, 1.25){$n-3$};
%\node at (13, 1.25){$n-2$};
\node at (14, 1.25){$n-1$};\node at (-1.5, 0){0};\node at (-1.5, -1){1};\node at (-1.5, -2){2};\node at (-1.5, -3.5){3};\node at (-1.5, -4.5){4};\node at (-1.5, -6){5};\node at (-1.5, -7){6};\node at (-1.5, -10){$m-2$};\node at (-1.5, -11){$m-1$};

\foreach \y in {-1, -4, -6.5, -10.5}{\node at (10.5,\y){\large$\cdots$};}
\foreach \x in {1, 4.5, 8, 13}{\node at (\x,-8.5){\large \vdots};}

% A(0)
    \draw[red, ultra thick, dashed] (-.5, -2) -- (0, -2) -- (0, -2.5); \draw[blue, ultra thick, densely dotted] (-.5, -1) -- (0, -1) -- (0, -2) -- (1, -2) -- (1, -2.5); \draw[orange, ultra thick] (-.5, 0) -- (0, 0) -- (0, -1) -- (1, -1) -- (1, -2) -- (2, -2) -- (2, -2.5); \draw[red, ultra thick, dashed] (0, .5) -- (0, 0) -- (1, 0) -- (1, -1) -- (2, -1) -- (2, -2) -- (2.5, -2); \draw[blue, ultra thick, densely dotted] (1, .5) -- (1, 0) -- (2, 0) -- (2, -1) -- (2.5, -1); \draw[orange, ultra thick] (2, .5) -- (2, 0) -- (2.5, 0); \foreach \x in {0, 1, 2}{\foreach \y in {0, -1, -2}{\draw[draw = black, fill = white] (\x, \y) circle (3pt);}}
    
% A(1)
    \draw[red, ultra thick, dashed] (3, -2) -- (3.5, -2) -- (3.5, -2.5); \draw[blue, ultra thick, densely dotted] (3, -1) -- (3.5, -1) -- (3.5, -2) -- (4.5, -2) -- (4.5, -2.5); \draw[orange, ultra thick] (3, 0) -- (3.5, 0) -- (3.5, -1) -- (4.5, -1) -- (4.5, -2) -- (5.5, -2) -- (5.5, -2.5); \draw[red, ultra thick, dashed] (3.5, .5) -- (3.5, 0) -- (4.5, 0) -- (4.5, -1) -- (5.5, -1) -- (5.5, -2) -- (6, -2); \draw[blue, ultra thick, densely dotted] (4.5, .5) -- (4.5, 0) -- (5.5, 0) -- (5.5, -1) -- (6, -1); \draw[orange, ultra thick] (5.5, .5) -- (5.5, 0) -- (6, 0); \foreach \x in {3.5, 4.5, 5.5}{\foreach \y in {0, -1, -2}{\draw[draw = black, fill = white] (\x, \y) circle (3pt);}}
    
% A(2)
    \draw[red, ultra thick, dashed] (6.5, -2) -- (7, -2) -- (7, -2.5); \draw[blue, ultra thick, densely dotted] (6.5, -1) -- (7, -1) -- (7, -2) -- (8, -2) -- (8, -2.5); \draw[orange, ultra thick] (6.5, 0) -- (7, 0) -- (7, -1) -- (8, -1) -- (8, -2) -- (9, -2) -- (9, -2.5); \draw[red, ultra thick, dashed] (7, .5) -- (7, 0) -- (8, 0) -- (8, -1) -- (9, -1) -- (9, -2) -- (9.5, -2); \draw[blue, ultra thick, densely dotted] (8, .5) -- (8, 0) -- (9, 0) -- (9, -1) -- (9.5, -1); \draw[orange, ultra thick] (9, .5) -- (9, 0) -- (9.5, 0); \foreach \x in {7, 8, 9}{\foreach \y in {0, -1, -2}{\draw[draw = black, fill = white] (\x, \y) circle (3pt);}}

% A(N-1)
    \draw[red, ultra thick, dashed] (12, -2) -- (12, -2.5); \draw[blue, ultra thick, densely dotted] (11.5, -1) -- (12, -1) -- (12, -2) -- (13, -2) -- (13, -2.5); \draw[orange, ultra thick] (11.5, 0) -- (12, 0) -- (12, -1) -- (13, -1) -- (13, -2) -- (14, -2) -- (14, -2.5); \draw[red, ultra thick, dashed] (12, .5) -- (12, 0) -- (13, 0) -- (13, -1) -- (14, -1) -- (14, -2) -- (14.5, -2); \draw[blue, ultra thick, densely dotted] (13, .5) -- (13, 0) -- (14, 0) -- (14, -1) -- (14.5, -1); \draw[orange, ultra thick] (14, .5) -- (14, 0) -- (14.5, 0); \foreach \x in {12, 13, 14}{\foreach \y in {0, -1, -2}{\draw[draw = black, fill = white] (\x, \y) circle (3pt);}}

% B(0, 0)
    \draw[red, ultra thick, dashed] (0, -3) -- (0, -3.5) -- (1, -3.5) -- (1, -4.5) -- (0, -4.5) -- (0, -5); \draw[blue, ultra thick, densely dotted] (1, -3) -- (1, -3.5) -- (2, -3.5) -- (2, -4.5) -- (1, -4.5) -- (1, -5); \draw[orange, ultra thick] (-.5, -3.5) -- (0, -3.5) -- (0, -4.5) -- (-.5, -4.5); \draw[orange, ultra thick] (2, -3) -- (2, -3.5) -- (2.5, -3.5); \draw[orange, ultra thick] (2.5, -4.5) -- (2, -4.5) -- (2, -5); \foreach \x in {0, 1, 2}{\foreach \y in {-3.5, -4.5}{\draw[draw = black, fill = white] (\x, \y) circle (3pt);}}
    
% B(0, 1)
   \draw[red, ultra thick, dashed] (3.5, -3) -- (3.5, -3.5) -- (4.5, -3.5) -- (4.5, -4.5) -- (3.5, -4.5) -- (3.5, -5); \draw[blue, ultra thick, densely dotted] (4.5, -3) -- (4.5, -3.5) -- (5.5, -3.5) -- (5.5, -4.5) -- (4.5, -4.5) -- (4.5, -5); \draw[orange, ultra thick] (3, -3.5) -- (3.5, -3.5) -- (3.5, -4.5) -- (3, -4.5); \draw[orange, ultra thick] (5.5, -3) -- (5.5, -3.5) -- (6, -3.5); \draw[orange, ultra thick] (6, -4.5) -- (5.5, -4.5) -- (5.5, -5); \foreach \x in {3.5, 4.5, 5.5}{\foreach \y in {-3.5, -4.5}{\draw[draw = black, fill = white] (\x, \y) circle (3pt);}}
   
% B(0, 2)
    \draw[red, ultra thick, dashed] (7, -3) -- (7, -3.5) -- (8, -3.5) -- (8, -4.5) -- (7, -4.5) -- (7, -5); \draw[blue, ultra thick, densely dotted] (8, -3) -- (8, -3.5) -- (9, -3.5) -- (9, -4.5) -- (8, -4.5) -- (8, -5); \draw[orange, ultra thick] (6.5, -3.5) -- (7, -3.5) -- (7, -4.5) -- (6.5, -4.5); \draw[orange, ultra thick] (9, -3) -- (9, -3.5) -- (9.5, -3.5); \draw[orange, ultra thick] (9.5, -4.5) -- (9, -4.5) -- (9, -5); \foreach \x in {7, 8, 9}{\foreach \y in {-3.5, -4.5}{\draw[draw = black, fill = white] (\x, \y) circle (3pt);}}
    
% B(0, N-1)
    \draw[red, ultra thick, dashed] (12, -3) -- (12, -3.5) -- (13, -3.5) -- (13, -4.5) -- (12, -4.5) -- (12, -5); \draw[blue, ultra thick, densely dotted] (13, -3) -- (13, -3.5) -- (14, -3.5) -- (14, -4.5) -- (13, -4.5) -- (13, -5); \draw[orange, ultra thick] (11.5, -3.5) -- (12, -3.5) -- (12, -4.5) -- (11.5, -4.5); \draw[orange, ultra thick] (14, -3) -- (14, -3.5) -- (14.5, -3.5); \draw[orange, ultra thick] (14.5, -4.5) -- (14, -4.5) -- (14, -5); \foreach \x in {12, 13, 14}{\foreach \y in {-3.5, -4.5}{\draw[draw = black, fill = white] (\x, \y) circle (3pt);}}
    
% B(1, 0)
    \draw[red, ultra thick, dashed] (0, -5.5) -- (0, -6) -- (1, -6) -- (1, -7) -- (0, -7) -- (0, -7.5); \draw[blue, ultra thick, densely dotted] (1, -5.5) -- (1, -6) -- (2, -6) -- (2, -7) -- (1, -7) -- (1, -7.5); \draw[orange, ultra thick] (-.5, -6) -- (0, -6) -- (0, -7) -- (-.5, -7); \draw[orange, ultra thick] (2, -5.5) -- (2, -6) -- (2.5, -6); \draw[orange, ultra thick] (2.5, -7) -- (2, -7) -- (2, -7.5); \foreach \x in {0, 1, 2}{\foreach \y in {-6, -7}{\draw[draw = black, fill = white] (\x, \y) circle (3pt);}}
    
% B(1, 1)
   \draw[red, ultra thick, dashed] (3.5, -5.5) -- (3.5, -6) -- (4.5, -6) -- (4.5, -7) -- (3.5, -7) -- (3.5, -7.5); \draw[blue, ultra thick, densely dotted] (4.5, -5.5) -- (4.5, -6) -- (5.5, -6) -- (5.5, -7) -- (4.5, -7) -- (4.5, -7.5); \draw[orange, ultra thick] (3, -6) -- (3.5, -6) -- (3.5, -7) -- (3, -7); \draw[orange, ultra thick] (5.5, -5.5) -- (5.5, -6) -- (6, -6); \draw[orange, ultra thick] (6, -7) -- (5.5, -7) -- (5.5, -7.5); \foreach \x in {3.5, 4.5, 5.5}{\foreach \y in {-6, -7}{\draw[draw = black, fill = white] (\x, \y) circle (3pt);}}
   
% B(1, 2)
    \draw[red, ultra thick, dashed] (7, -5.5) -- (7, -6) -- (8, -6) -- (8, -7) -- (7, -7) -- (7, -7.5); \draw[blue, ultra thick, densely dotted] (8, -5.5) -- (8, -6) -- (9, -6) -- (9, -7) -- (8, -7) -- (8, -7.5); \draw[orange, ultra thick] (6.5, -6) -- (7, -6) -- (7, -7) -- (6.5, -7); \draw[orange, ultra thick] (9, -5.5) -- (9, -6) -- (9.5, -6); \draw[orange, ultra thick] (9.5, -7) -- (9, -7) -- (9, -7.5); \foreach \x in {7, 8, 9}{\foreach \y in {-6, -7}{\draw[draw = black, fill = white] (\x, \y) circle (3pt);}}

% B(1, N-1)
    \draw[red, ultra thick, dashed] (12, -5.5) -- (12, -6) -- (13, -6) -- (13, -7) -- (12, -7) -- (12, -7.5); \draw[blue, ultra thick, densely dotted] (13, -5.5) -- (13, -6) -- (14, -6) -- (14, -7) -- (13, -7) -- (13, -7.5); \draw[orange, ultra thick] (11.5, -6) -- (12, -6) -- (12, -7) -- (11.5, -7); \draw[orange, ultra thick] (14, -5.5) -- (14, -6) -- (14.5, -6); \draw[orange, ultra thick] (14.5, -7) -- (14, -7) -- (14, -7.5); \foreach \x in {12, 13, 14}{\foreach \y in {-6, -7}{\draw[draw = black, fill = white] (\x, \y) circle (3pt);}}

% B(M-1, 0)
    \draw[red, ultra thick, dashed] (0, -9.5) -- (0, -10) -- (1, -10) -- (1, -11) -- (0, -11) -- (0, -11.5); \draw[blue, ultra thick, densely dotted] (1, -9.5) -- (1, -10) -- (2, -10) -- (2, -11) -- (1, -11) -- (1, -11.5); \draw[orange, ultra thick] (-.5, -10) -- (0, -10) -- (0, -11) -- (-.5, -11); \draw[orange, ultra thick] (2, -9.5) -- (2, -10) -- (2.5, -10); \draw[orange, ultra thick] (2.5, -11) -- (2, -11) -- (2, -11.5); \foreach \x in {0, 1, 2}{\foreach \y in {-10, -11}{\draw[draw = black, fill = white] (\x, \y) circle (3pt);}}

% B(M-1, 1)
    \draw[red, ultra thick, dashed] (3.5, -9.5) -- (3.5, -10) -- (4.5, -10) -- (4.5, -11) -- (3.5, -11) -- (3.5, -11.5); \draw[blue, ultra thick, densely dotted] (4.5, -9.5) -- (4.5, -10) -- (5.5, -10) -- (5.5, -11) -- (4.5, -11) -- (4.5, -11.5); \draw[orange, ultra thick] (3, -10) -- (3.5, -10) -- (3.5, -11) -- (3, -11); \draw[orange, ultra thick] (5.5, -9.5) -- (5.5, -10) -- (6, -10); \draw[orange, ultra thick] (6, -11) -- (5.5, -11) -- (5.5, -11.5); \foreach \x in {3.5, 4.5, 5.5}{\foreach \y in {-10, -11}{\draw[draw = black, fill = white] (\x, \y) circle (3pt);}}

% B(M-1, 2)
    \draw[red, ultra thick, dashed] (7, -9.5) -- (7, -10) -- (8, -10) -- (8, -11) -- (7, -11) -- (7, -11.5); \draw[blue, ultra thick, densely dotted] (8, -9.5) -- (8, -10) -- (9, -10) -- (9, -11) -- (8, -11) -- (8, -11.5); \draw[orange, ultra thick] (6.5, -10) -- (7, -10) -- (7, -11) -- (6.5, -11); \draw[orange, ultra thick] (9, -9.5) -- (9, -10) -- (9.5, -10); \draw[orange, ultra thick] (9.5, -11) -- (9, -11) -- (9, -11.5); \foreach \x in {7, 8, 9}{\foreach \y in {-10, -11}{\draw[draw = black, fill = white] (\x, \y) circle (3pt);}}

% B(M-1, N-1)
    \draw[red, ultra thick, dashed] (12, -9.5) -- (12, -10) -- (13, -10) -- (13, -11) -- (12, -11) -- (12, -11.5); \draw[blue, ultra thick, densely dotted] (13, -9.5) -- (13, -10) -- (14, -10) -- (14, -11) -- (13, -11) -- (13, -11.5); \draw[orange, ultra thick] (11.5, -10) -- (12, -10) -- (12, -11) -- (11.5, -11); \draw[orange, ultra thick] (14, -9.5) -- (14, -10) -- (14.5, -10); \draw[orange, ultra thick] (14.5, -11) -- (14, -11) -- (14, -11.5); \foreach \x in {12, 13, 14}{\foreach \y in {-10, -11}{\draw[draw = black, fill = white] (\x, \y) circle (3pt);}}
\end{tikzpicture}

\vspace{.25in}

%A/B Table
\begin{tikzpicture}[scale=1, transform shape]
    \node at (10, 0){$A(0)$}; \node at (10, -1){$B(0,0)$}; \node at (10, -2){$B(1,0)$}; \node at (10, -3){$\vdots$}; \node at (10, -4){$B(M-1,0)$};\node at (13, 0){$A(1)$}; \node at (13, -1){$B(0,1)$}; \node at (13, -2){$B(1,1)$}; \node at (13, -3){$\vdots$};\node at (13, -4){$B(M-1,1)$};\node at (16, 0){$A(2)$}; \node at (16, -1){$B(0,2)$}; \node at (16, -2){$B(1,2)$};\node at (16, -3){$\vdots$}; \node at (16, -4){$B(M-1,2)$};\node at (19, 0){$\cdots$}; \node at (19, -1){$\cdots$}; \node at (19, -2){$\cdots$};\node at (19, -4){$\cdots$};
    \node at (22, 0){$A(N-1)$}; \node at (22, -1){$B(0,N-1)$}; \node at (22, -2){$B(1,N-1)$}; \node at (22, -3){$\vdots$}; \node at (22, -4){$B(M-1,N-1)$};
\end{tikzpicture}
\end{center}
\caption{Top left: The $A$ block for Case 1, with $A^r$, $A^y$, and $A^b$ in dashed red, solid yellow, and dotted blue. Top right: The $B$ block for Case 1, with $B^r$, $B^y$, and $B^b$ in dashed red, solid yellow, and dotted blue.  Middle: The decomposition of $C_m \sq C_n$ into three cycles in Case 1 ($m$ odd). Bottom: The arrangement of $A$ and $B$ blocks in the decomposition of $C_m \sq C_n$ in Case 1.}\label{ABcase1}
\end{figure}

\paragraph{Case 2 ($m=4$ and $n$ is even):} Assume $m=4$ and $n=6+6N$ for a nonnegative integer $N$. The decomposition of $C_m \sq C_n$ is illustrated in Figure~\ref{ABcase2} (bottom) and contains one $4 \times 6$ $A$ block and $N$ $4 \times 6$ $B$ blocks, defined below.  The $A$ block is the minimal block for this case.

For the $ 4 \times 6$ $A$ block, we define three edge sets $A^r$, $A^y$, and $A^b$, each of which forms a cycle with 16 edges, shown as dashed red, solid yellow, and dotted blue, respectively, in Figure~\ref{ABcase2} (top left).

$A^r = \{(1, -1)(1, 0), (1, 0)(2, 0), (2, 0)(2, 1), (2, 1)(1, 1), (1,1)(1,2),(1,2)(2,2),\\(2,2)(3,2), (3,2)(4,2), (-1,2)(0,2), (0,2)(0,3), (0,3)(1,3), (1,3)(2,3), \\(2,3)(2,4), (2,4)(3,4), (3,4)(4,4), (-1,4)(0,4), (0,4)(0,5), (0,5)(1,5), (1,5)(1,6)\}$

$A^y = \{(-1, 0)(0,0), (0,0)(1,0), (1,0)(1,1), (1,1)(0,1), (0,1)(-1,1), (4, 1)(3,1), \\(3,1)(2,1), (2,1)(2,2), (2,2)(2,3), (2,3)(3,3), (3,3)(4,3), (-1, 3)(0,3), (0,3)(0,4),\\ (0,4)(1,4), (1,4)(2,4), (2,4)(2,5), (2,5)(3,5), (3,5)(3,6), (3, -1)(3,0), (3,0)(4,0)\}$

$A^b = \{(0,-1)(0,0), (0,0)(0,1), (0,1)(0,2), (0,2)(1,2), (1,2)(1,3), (1,3)(1,4),\\ (1,4)(1,5), (1,5)(2,5), (2,5)(2,6), (2,-1)(2,0), (2,0)(3,0), (3,0)(3,1), (3,1)(3,2), \\(3,2)(3,3), (3,3)(3,4), (3,4)(3,5), (3,5)(4,5), (-1,5)(0,5), (0,5)(0,6)\}$

For the $4 \times 6$ $B$ block define three edge sets $B^r$, $B^y$, and $B^b$, each containing 16 edges and shown as dashed red, solid yellow, and dotted blue, respectively, in Figure~\ref{ABcase2} (top right).

$B^r = \{(4,0)(3,0), (3,0)(3,1), (3,1)(3,2), (3,2)(2,2), (2,2)(1,2), (1,2)(1,3),\\ (1,3)(0,3), (0,3)(0,4), (0,4)(1,4), (1,4)(2,4), (2,4)(3,4), (3,4)(3,5), (3,5)(2,5),\\ (2,5)(1,5), (1,5)(1,6), (-1, 0)(0,0), (0,0)(1,0), (1,0)(1,-1)\}$

$B^y = \{(-1,1)(0,1), (0,1)(1,1), (1,1)(1,2), (1,2)(0,2), (0,2)(-1,2), (4,2)(3,2), \\(3,2)(3,3), (3,3)(2,3), (2,3)(1,3), (1,3)(1,4), (1,4)(1,5), (1,5)(0,5), (0,5)(-1,5),\\ (4,5)(3,5), (3,5)(3,6), (3, -1)(3,0), (3,0)(2,0), (2,0)(2,1), (2,1)(3,1), (3,1)(4,1)\}$

$B^b = \{(0, -1)(0,0), (0,0)(0,1), (0,1)(0,2), (0,2)(0,3), (0,3)(-1,3), (4,3)(3,3), \\(3,3)(3,4), (3,4)(4,4), (-1,4)(0,4), (0,4)(0,5), (0,5)(0,6), (2,-1)(2,0), (2,0)(1,0), \\(1,0)(1,1), (1,1)(2,1), (2,1)(2,2), (2,2)(2,3), (2,3)(2,4), (2,4)(2,5), (2,5)(2,6)\}$

For $0 \leq J \leq N-1$, let the block $B(J)$ in $C_m \sq C_n$ be the set of edges incident to the vertices $(i, j)$ where $0\leq i \leq 3$ and $6+6J\leq j\leq 11+6J$. To color the edges in $B(J)$, define three edge sets $B^r(J)$, $B^y(J)$, and $B^b(J)$, where for $\ell \in \{r,y,b\}$, edge $(i,j)(i',j')$ is in $B^\ell$ if and only if $(i,j+6+6J)(i',j'+6+6J)$ is in $B^\ell(J)$. For example, $(1,13)(1,14)$ is yellow (in $B^y(1)$) because $(1,1)(1,2)$ is yellow (in $B^y$).

For $\ell \in \{r,y,b\}$ each cycle in the decomposition is the union of the $A^\ell$ and $B^\ell(J)$ sets, i.e.

\[ A^\ell \cup
\left(\bigcup_{J=0}^{N-1}B^\ell(J) \right).\]

Since the $A$ block contains 16 edges of each color, and each of the $N$ $B$ blocks contain 16 edges of each color, each cycle contains $16+16N$ edges.  Since $m=4$ and $n=6+6N$, the total number of edges in $C_m \sq C_n$ is $2mn= 2(4)(6+6N) = 48+48N$.  Thus each cycle contains one third of the total edges, as desired.

\begin{figure}
\begin{center}
\begin{tikzpicture}
    \draw[blue, ultra thick, densely dotted] (-.5, 0) -- (2, 0) -- (2, -1) -- (5, -1) -- (5, -2) -- (5.5, -2); \draw[blue, ultra thick, densely dotted] (-.5, -2) -- (0, -2) -- (0, -3) -- (5, -3) -- (5, -3.5); \draw[blue, ultra thick, densely dotted] (5, .5) -- (5, 0) -- (5.5, 0);
    \draw[orange, ultra thick] (0, .5) -- (0, -1) -- (1, -1) -- (1, .5); \draw[orange, ultra thick] (3, .5) -- (3, 0) -- (4, 0) -- (4, -2) -- (5, -2) -- (5, -3) -- (5.5, -3); \draw[orange, ultra thick] (-.5, -3) -- (0, -3) -- (0, -3.5); \draw[orange, ultra thick] (1, -3.5) -- (1, -2) -- (3, -2) -- (3, -3.5);
    \draw[red, ultra thick, dashed] (-.5, -1) -- (0, -1) -- (0, -2) -- (1, -2) -- (1, -1) -- (2, -1) -- (2, -3.5); \draw[red, ultra thick, dashed] (2, .5) -- (2, 0) -- (3, 0) -- (3, -2) -- (4, -2) -- (4, -3.5); \draw[red, ultra thick, dashed] (4, .5) -- (4, 0) -- (5, 0) -- (5, -1) -- (5.5, -1);\foreach \x in {0, 1, ..., 5}{\foreach \y in {0, -1, -2, -3}{\filldraw[draw=black, fill=white] (\x, \y) circle (3pt);}}
    \foreach \x in {0,1,2,3,4,5}{\node at (\x,1){$\x$};}\foreach \y\i in {0/0,-1/1,-2/2,-3/3}{\node at (-1,\y){$\i$};}
\end{tikzpicture}
\hspace{.5in}
%B Block
\begin{tikzpicture}[scale=1, transform shape]
    \draw[blue, ultra thick, densely dotted] (6, 0) -- (9.5, 0) -- (9.5, .5); \draw[blue, ultra thick, densely dotted] (10.5, .5) -- (10.5, 0) -- (12, 0); \draw[blue, ultra thick, densely dotted] (6, -2) -- (6.5, -2) -- (6.5, -1) -- (7.5, -1) -- (7.5, -2) -- (12, -2); 
    \draw[blue, ultra thick, densely dotted] (9.5, -3.5) -- (9.5, -3) -- (10.5, -3) -- (10.5, -3.5);
    \draw[orange, ultra thick] (7.5, .5) -- (7.5, -1) -- (8.5, -1) -- (8.5, .5); \draw[orange, ultra thick] (6, -3) -- (6.5, -3) -- (6.5, -2) -- (7.5, -2) -- (7.5, -3.5); \draw[orange, ultra thick] (8.5, -3.5) -- (8.5, -3) -- (9.5, -3) -- (9.5, -1) -- (11.5, -1) -- (11.5, .5); \draw[orange, ultra thick] (11.5, -3.5) -- (11.5, -3) -- (12, -3);
    \draw[red, ultra thick, dashed] (6, -1) -- (6.5, -1) -- (6.5, .5); \draw[red, ultra thick, dashed] (6.5, -3.5) -- (6.5, -3) -- (8.5, -3) -- (8.5, -1) -- (9.5, -1) -- (9.5, 0) -- (10.5, 0) -- (10.5, -3) -- (11.5, -3) -- (11.5, -1) -- (12, -1);
    \foreach \x in {6.5, 7.5, 8.5, 9.5, 10.5, 11.5}{\foreach \y in {0, -1, -2, -3}{\filldraw[draw=black, fill=white] (\x, \y) circle (3pt);}}
    \foreach \x\j in {6.5/0,7.5/1,8.5/2,9.5/3,10.5/4,11.5/5}{\node at (\x,1){$\j$};}\foreach \y\i in {0/0,-1/1,-2/2,-3/3}{\node at (5.5,\y){$\i$};}
\end{tikzpicture}

\vspace{.25in}

\begin{tikzpicture}[scale=.55, transform shape]
% Vertex labels
\foreach \x\i in {0/0, 1/1, 2/2, 3/3, 4/4, 5/5, 6.5/6, 7.5/7, 8.5/8, 9.5/9, 10.5/10, 11.5/11, 13/12, 14/13, 15/14, 16/15, 17/16, 18/17, 19.5/ , 21.5/$n-5$, 22.5/ , 23.5/$n-3$, 24.5/ , 25.5/$n-1$}{\node at (\x, 1.25){\i};}\foreach \y in {0,1,2,3}{\node at (-1, -\y){\y};}

% A(0)
    \draw[blue, ultra thick, densely dotted] (-.5, 0) -- (2, 0) -- (2, -1) -- (5, -1) -- (5, -2) -- (5.5, -2); \draw[blue, ultra thick, densely dotted] (-.5, -2) -- (0, -2) -- (0, -3) -- (5, -3) -- (5, -3.5); \draw[blue, ultra thick, densely dotted] (5, .5) -- (5, 0) -- (5.5, 0);
    \draw[orange, ultra thick] (0, .5) -- (0, -1) -- (1, -1) -- (1, .5); \draw[orange, ultra thick] (3, .5) -- (3, 0) -- (4, 0) -- (4, -2) -- (5, -2) -- (5, -3) -- (5.5, -3); \draw[orange, ultra thick] (-.5, -3) -- (0, -3) -- (0, -3.5); \draw[orange, ultra thick] (1, -3.5) -- (1, -2) -- (3, -2) -- (3, -3.5);
    \draw[red, ultra thick, dashed] (-.5, -1) -- (0, -1) -- (0, -2) -- (1, -2) -- (1, -1) -- (2, -1) -- (2, -3.5); \draw[red, ultra thick, dashed] (2, .5) -- (2, 0) -- (3, 0) -- (3, -2) -- (4, -2) -- (4, -3.5); \draw[red, ultra thick, dashed] (4, .5) -- (4, 0) -- (5, 0) -- (5, -1) -- (5.5, -1);\foreach \x in {0, 1, ..., 5}{\foreach \y in {0, -1, -2, -3}{\filldraw[draw=black, fill=white] (\x, \y) circle (3pt);}}
    
    \node at (2.5, -4.5){\huge{$A$}}; 
% B(0)
    \draw[blue, ultra thick, densely dotted] (6, 0) -- (9.5, 0) -- (9.5, .5); \draw[blue, ultra thick, densely dotted] (10.5, .5) -- (10.5, 0) -- (12, 0); \draw[blue, ultra thick, densely dotted] (6, -2) -- (6.5, -2) -- (6.5, -1) -- (7.5, -1) -- (7.5, -2) -- (12, -2); \draw[blue, ultra thick, densely dotted] (9.5, -3.5) -- (9.5, -3) -- (10.5, -3) -- (10.5, -3.5);
    \draw[orange, ultra thick] (7.5, .5) -- (7.5, -1) -- (8.5, -1) -- (8.5, .5); \draw[orange, ultra thick] (6, -3) -- (6.5, -3) -- (6.5, -2) -- (7.5, -2) -- (7.5, -3.5); \draw[orange, ultra thick] (8.5, -3.5) -- (8.5, -3) -- (9.5, -3) -- (9.5, -1) -- (11.5, -1) -- (11.5, .5); \draw[orange, ultra thick] (11.5, -3.5) -- (11.5, -3) -- (12, -3);
    \draw[red, ultra thick, dashed] (6, -1) -- (6.5, -1) -- (6.5, .5); \draw[red, ultra thick, dashed] (6.5, -3.5) -- (6.5, -3) -- (8.5, -3) -- (8.5, -1) -- (9.5, -1) -- (9.5, 0) -- (10.5, 0) -- (10.5, -3) -- (11.5, -3) -- (11.5, -1) -- (12, -1);
    \foreach \x in {6.5, 7.5, 8.5, 9.5, 10.5, 11.5}{\foreach \y in {0, -1, -2, -3}{\filldraw[draw=black, fill=white] (\x, \y) circle (3pt);}}
    
    \node at (9, -4.5){\huge{$B(0)$}};
% B(1)
    \draw[blue, ultra thick, densely dotted] (12.5, 0) -- (16, 0) -- (16, .5); \draw[blue, ultra thick, densely dotted] (17, .5) -- (17, 0) -- (18.5, 0); \draw[blue, ultra thick, densely dotted] (12.5, -2) -- (13, -2) -- (13, -1) -- (14, -1) -- (14, -2) -- (18.5, -2); \draw[blue, ultra thick, densely dotted] (16, -3.5) -- (16, -3) -- (17, -3) -- (17, -3.5);
    \draw[orange, ultra thick] (14, .5) -- (14, -1) -- (15, -1) -- (15, .5); \draw[orange, ultra thick] (12.5, -3) -- (13, -3) -- (13, -2) -- (14, -2) -- (14, -3.5); \draw[orange, ultra thick] (15, -3.5) -- (15, -3) -- (16, -3) -- (16, -1) -- (18, -1) -- (18, .5); \draw[orange, ultra thick] (18, -3.5) -- (18, -3) -- (18.5, -3);
    \draw[red, ultra thick, dashed] (12.5, -1) -- (13, -1) -- (13, .5); \draw[red, ultra thick, dashed] (13, -3.5) -- (13, -3) -- (15, -3) -- (15, -1) -- (16, -1) -- (16, 0) -- (17, 0) -- (17, -3) -- (18, -3) -- (18, -1) -- (18.5, -1);
    \foreach \x in {13, 14, 15, 16, 17, 18}{\foreach \y in {0, -1, -2, -3}{\filldraw[draw=black, fill=white] (\x, \y) circle (3pt);}}
    
    \node at (15.5, -4.5){\huge{$B(1)$}}; 
\node at (19.25, -1.5){$\cdots$};
\node at (19.25, -4.5){\huge{$\cdots$}};
% B(N-1)
    \draw[blue, ultra thick, densely dotted] (20, 0) -- (23.5, 0) -- (23.5, .5); \draw[blue, ultra thick, densely dotted] (24.5, .5) -- (24.5, 0) -- (26, 0); \draw[blue, ultra thick, densely dotted] (20, -2) -- (20.5, -2) -- (20.5, -1) -- (21.5, -1) -- (21.5, -2) -- (26, -2); \draw[blue, ultra thick, densely dotted] (23.5, -3.5) -- (23.5, -3) -- (24.5, -3) -- (24.5, -3.5);
    \draw[orange, ultra thick] (21.5, .5) -- (21.5, -1) -- (22.5, -1) -- (22.5, .5); \draw[orange, ultra thick] (20, -3) -- (20.5, -3) -- (20.5, -2) -- (21.5, -2) -- (21.5, -3.5); \draw[orange, ultra thick] (22.5, -3.5) -- (22.5, -3) -- (23.5, -3) -- (23.5, -1) -- (25.5, -1) -- (25.5, .5); \draw[orange, ultra thick] (25.5, -3.5) -- (25.5, -3) -- (26, -3);
    \draw[red, ultra thick, dashed] (20, -1) -- (20.5, -1) -- (20.5, .5); \draw[red, ultra thick, dashed] (20.5, -3.5) -- (20.5, -3) -- (22.5, -3) -- (22.5, -1) -- (23.5, -1) -- (23.5, 0) -- (24.5, 0) -- (24.5, -3) -- (25.5, -3) -- (25.5, -1) -- (26, -1);
    \foreach \x in {20.5, 21.5, 22.5, 23.5, 24.5, 25.5}{\foreach \y in {0, -1, -2, -3}{\filldraw[draw=black, fill=white] (\x, \y) circle (3pt);}}
    
    \node at (23, -4.5){\huge{$B(N-1)$}};
\end{tikzpicture}
\end{center}
    \caption{Above left: The $A$ block for Case 2, with $A^r$, $A^y$, and $A^b$ in dashed red, solid yellow, and dotted blue. Above right: The $B$ block for Case 2, with $B^r$, $B^y$, and $B^b$ in dashed red, solid yellow, and dotted blue. Below: The decomposition of $C_m \sq C_n$ into three cycles, and the arrangement of $A$ and $B$ blocks in $C_m \sq C_n$ in Case 2.}\label{ABcase2}
\end{figure}
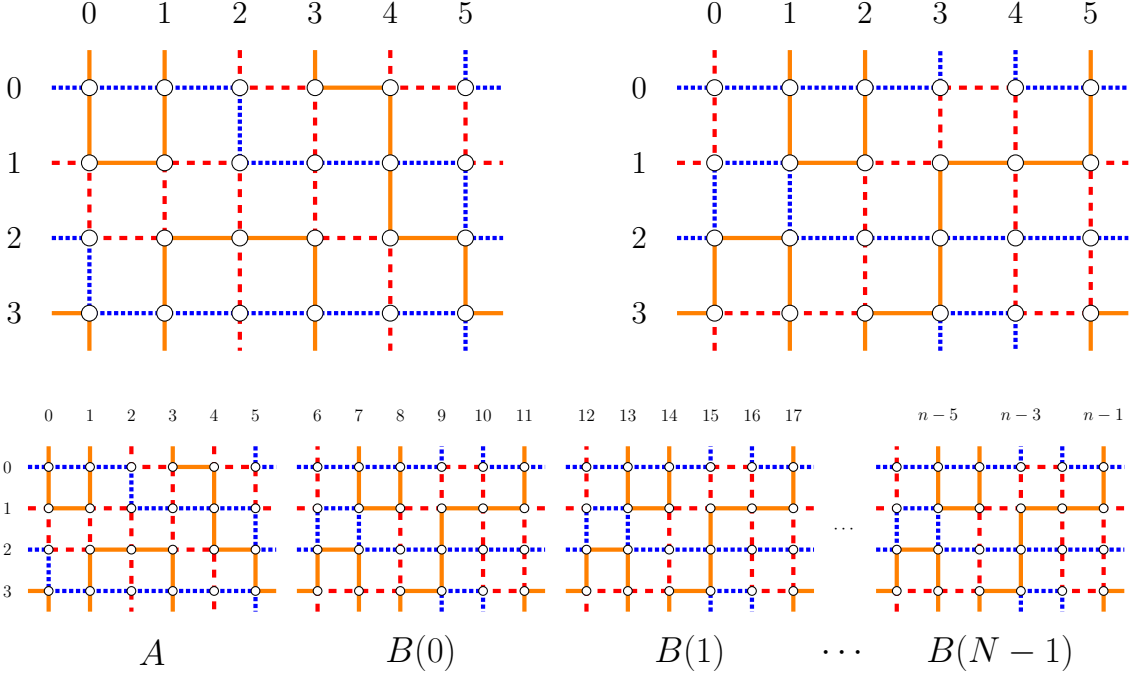

\paragraph{Case 3 ($m=4$ and $n=3$):} Assume $m=4$ and $n=3$. For the single $4 \times 3$ $A$ block, define three edge sets $A^r$, $A^y$, and $A^b$, each of which forms a cycle with eight edges, shown as dashed red, solid yellow, and dotted blue, respectively, in Figure~\ref{Acase3}.

$A^r = \{(-1, 0)(0, 0), (0,0)(0,1), (0, 1)(1, 1), (1, 1)(1, 0), (1, 0)(2,0), (2, 0)(2,1),\\ (2,1)(3,1), (3,1)(3,0), (3,0)(4,0)\}$

$A^y = \{(0,-1)(0,0), (0,0)(1,0), (1,0)(1,-1), (1,3)(1,2), (1,2)(2,2),(2,2)(2,3),\\ (2,-1)(2,0), (2,0)(3,0), (3,0)(3,-1), (3,3)(3,2), (3,2)(4,2), (-1,2)(0,2), (0,2)(0,3)\}$

$A^b = \{(-1, 1)(0, 1), (0,1)(0,2), (0, 2)(1, 2), (1, 2)(1, 1), (1, 1)(2,1), (2, 1)(2,2),\\ (2,2)(3,2), (3,2)(3,1), (3,1)(4,1)\}$

%The decomposition of $C_m \sq C_n$ is illustrated in Figure~\ref{Case2}.

\begin{figure}
\begin{center}
\begin{tikzpicture}[scale=.7]
    \draw[red, ultra thick, dashed] (0, -3) -- (0, -3.5) -- (1, -3.5) -- (1, -4.5) -- (0, -4.5) -- (0, -5.5) -- (1,-5.5) -- (1,-6.5) -- (0,-6.5) -- (0,-7); \draw[blue, ultra thick, densely dotted] (1, -3) -- (1, -3.5) -- (2, -3.5) -- (2, -4.5) -- (1, -4.5) -- (1, -5.5) -- (2, -5.5) -- (2, -6.5) -- (1,-6.5) -- (1, -7); 
    \draw[orange, ultra thick] (-.5, -3.5) -- (0, -3.5) -- (0, -4.5) -- (-.5, -4.5); 
    \draw[orange, ultra thick] (-.5, -5.5) -- (0, -5.5) -- (0, -6.5) -- (-.5, -6.5); 
    \draw[orange, ultra thick] (2, -3) -- (2, -3.5) -- (2.5, -3.5); 
    \draw[orange, ultra thick] (2.5, -4.5) -- (2, -4.5) -- (2, -5.5) -- (2.5, -5.5); 
    \draw[orange, ultra thick] (2.5, -6.5) -- (2, -6.5) -- (2, -7); 
    \foreach \x in {0, 1, 2}{\foreach \y in {-3.5, -4.5, -5.5, -6.5}{\draw[draw = black, fill = white] (\x, \y) circle (3pt);}}
    \foreach \x in {0,1,2}{\node at (\x,-2.5){$\x$};}\foreach \y\i in {-3.5/0,-4.5/1,-5.5/2,-6.5/3}{\node at (-1,\y){$\i$};}
\end{tikzpicture}
\end{center}
    \caption{The $A$ block for Case 3, with $A^r$, $A^y$, and $A^b$ in dashed red, solid yellow, and dotted blue.}
    \label{Acase3}
\end{figure}
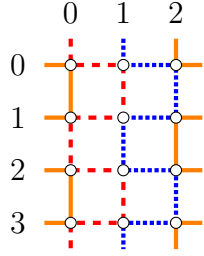

\paragraph{Case 4 ($m=4$ and $n\ge 9$ is odd):} Assume $m=4$ and $n=9+6N$ for a nonnegative integer $N$. For this case, our minimal $4 \times 9$ $A$ block is the union of three $4 \times 3$ blocks $A_1$, $A_2$, and $A_3$, shown in Figure~\ref{ACase4} (top).  The decomposition of $C_m \sq C_n$ contains one of each of the blocks $A_1$, $A_2$, and $A_3$, and $N$ of each of three types of $4 \times 2$ blocks $B_1$, $B_2$, and $B_3$, defined below and shown in Figure~\ref{ACase4} (middle below).

In the decomposition of $C_m \sq C_n$, $A_1$ contains edges incident to the vertices $(i,j)$ where $0\leq i\leq 3$ and  $2N \leq j \leq 2+2N$, $A_2$ contains edges incident to the vertices $(i,j)$ where $0\leq i\leq 3$ and $3+4N \leq j \leq 5+4N$, and $A_3$ contains edges incident to the vertices $(i,j)$ where  $0\leq i\leq 3$ and $6+6N \leq j \leq 8+6N$.
For $k \in \{1,2,3\}$ and $\ell \in \{r,y,b\}$, we define the sets of edges $A_k^r$, $A_k^y$, and $A_k^b$, which are depicted as dashed red, solid yellow, and dotted blue, respectively, in Figure~\ref{ACase4} (top).

Note that when $N=0$, the $A_1$, $A_2$, and $A_3$ blocks form the minimal $4 \times 9$ $A$ block depicted in Figure~\ref{ACase4} (middle above), and in this case, for $\ell \in \{r,y,b\}$ the union $A_1^\ell \cup A_2^\ell \cup A_3^\ell$ is a cycle with 24 edges.

$A_1^r = \{(-1, 0)(0,0), (0,0)(1,0), (1,0)(1,1), (1,1)(0,1), (0,1)(0,2), (0,2)(-1,2),\\ (4,2)(3,2), (3,2)(3,3), (3,-1)(3,0), (3,0)(4,0)\}.$

$A_1^y = \{(0,-1)(0,0), (0,0)(0,1), (0,1)(-1,1), (4,1)(3,1), (3,1)(3,0), (3,0)(2,0),\\ (2,0)(1,0), (1,0)(1,-1)\}.$

$A_1^b = \{(0,3)(0,2), (0,2)(1,2), (1,2)(1,1), (1,1)(2,1), (2,1)(3,1), (3,1)(3,2),\\ (3,2)(2,2), (2,2)(2,3)\}.$

$A_2^r = \{(-1,0)(0,0), (0,0)(0,1), (0,1)(1,1), (1,1)(2,1), (2,1)(2,0), (2,0)(1,0), \\(1,0)(1,-1), (4,0)(3,0), (3,0)(3,-1)\}.$

$A_2^y = \{(-1,1)(0,1), (0,1)(0,2), (0,2)(-1,2), (4,2)(3,2), (3,2)(3,3), (4,1)(3,1),\\ (3,1)(2,1), (2,1)(2,2), (2,2)(1,2), (1,2)(1,3)\}.$

$A_2^b = \{(0,-1)(0,0), (0,0)(1,0), (1,0)(1,1), (1,1)(1,2), (1,2)(0,2), (0,2)(0,3), \\(2,-1)(2,0), (2,0)(3,0), (3,0)(3,1), (3,1)(3,2), (3,2)(2,2), (2,2)(2,3)\}.$

$A_3^r = \{(-1,1)(0,1), (0,1)(0,2), (0,2)(1,2), (1,2)(2,2), (2,2)(2,3), (4,1)(3,1), \\(3,1)(3,2), (3,2)(3,3)\}.$

$A_3^y = \{(-1,0)(0,0), (0,0)(0,1), (0,1)(1,1), (1,1)(1,2), (1,2)(1,3), (1,-1)(1,0),\\ (1,0)(2,0), (2,0)(2,1), (2,1)(2,2), (2,2)(3,2), (3,2)(4,2), (3,-1)(3,0), (3,0)(4,0)\}.$

$A_3^b = \{(0,-1)(0,0), (0,0)(1,0), (1,0)(1,1), (1,1)(2,1), (2,1)(3,1), (3,1)(3,0), \\(3,0)(2,0), (2,0)(2,-1)\}.$

For the inductive step, we define three $4 \times 2$ blocks $B_1$, $B_2$, $B_3$, and for $k \in \{1,2,3\}$ and $\ell \in \{r,y,b\}$, we define the sets of edges $B_k^r$, $B_k^y$, and $B_k^b$, shown as dashed red, solid yellow, and dotted blue, respectively, in Figure~\ref{ACase4} (middle below). Note that each $B_1$, $B_2$, and $B_3$ block contains eight edges each of two of the colors, and no edges of the third color.

$B_1^r = \{(2,-1)(2,0), (2,0)(1,0), (1,0)(0,0), (0,0)(0,1), (0,1)(1,1), (1,1)(2,1), \\(2,1)(2,3), (3,-1)(3,0), (3,0)(3,1), (3,1)(3,2)\}.$

$B_1^y = \{(0,-1)(0,0), (0,0)(-1,0), (4,0)(3,0), (3,0)(2,0), (2,0)(2,1), (2,1)(3,1),\\ (3,1)(4,1), (-1,1)(0,1), (0,1)(0,2), (1,-1)(1,0), (1,0)(1,1), (1,1)(1,2)\}.$

$B_1^b = \{ \}$.

$B_2^r = \{(1,-1)(1,0), (1,0)(0,0), (0,0)(0,1), (0,1)(1,1), (1,1)(1,2), (3,-1)(3,0), \\(3,0)(2,0), (2,0)(2,1), (2,1)(3,1), (3,1)(3,2)\}.$

$B_2^y = \{ \}$.

$B_2^b = \{(2,-1)(2,0), (2,0)(1,0), (1,0)(1,1), (1,1)(2,1), (2,1)(2,2), (0,-1)(0,0),\\ (0,0)(-1,0), (-1,1)(0,1), (0,1)(0,2), (4,0)(3,0), (3,0)(3,1), (3,1)(4,1)\}.$

$B_3^r = \{ \}$.

$B_3^y = \{(1,-1)(1,0), (1,0)(0,0), (0,0)(0,1), (0,1)(1,1), (1,1)(1,2), (3,-1)(3,0), \\(3,0)(2,0), (2,0)(2,1), (2,1)(3,1), (3,1)(3,2)\}.$

$B_3^b = \{(2,-1)(2,0), (2,0)(1,0), (1,0)(1,1), (1,1)(2,1), (2,1)(2,2), (0,-1)(0,0),\\ (0,0)(-1,0), (-1,1)(0,1), (0,1)(0,2), (4,0)(0,0), (3,0)(3,1), (3,1)(4,1)\}.$

Unlike the other cases, rather than being added next to the $A$ block in the inductive step, the blocks $B_1$, $B_2$, and $B_3$ are spliced in between the blocks $A_1$, $A_2$, and $A_3$ in a way that the cycles of the $A$ block are extended, as shown in Figure~\ref{ACase4} (bottom).  See also Figure~\ref{Ex1}. The $B_1$ blocks fit between $A_3$ and $A_1$, the $B_2$ blocks fit between $A_1$ and $A_2$, and the $B_3$ blocks fit between $A_2$ and $A_3$. 

For $0\leq J\leq N-1,$ let
\begin{itemize}
\item $B_1(J)$ be the set of edges incident to vertices $(i, j)$ where $0\leq i\leq 3$ and $2J\leq j \leq 1+2J$.
\item $B_2(J)$ be the set of edges incident to vertices $(i, j)$ where $0\leq i\leq 3$ and $3+2N+2J\leq j \leq 4+ 2N+2J.$ 
\item  $B_3(J)$ be the set of edges incident to vertices $(i, j)$ where $0\leq i \leq 3$ and $6+4N +2J\leq j \leq 7+4N+2J$.
\end{itemize}

For $k \in \{1,2,3\}$, to color the edges in $B_k(J)$, define three edge sets $B_k^r(J)$, $B_k^y(J)$, and $B_k^b(J)$, where for $\ell \in \{r,y,b\}$, edge $(i,j)(i',j')$ is in $B_k^\ell$ if and only if $(i,j+2J)(i',j'+2J)$ is in $B_k^\ell(J)$. For example, $(1,4N+1)(2,4N+2)$ is blue (in $B_2^b(N-1)$) because $(1,2N+3)(2,2N+3)$ is blue (in $B_2^b$).

The decomposition of $C_m \sq C_n$ is illustrated in Figure~\ref{ACase4} (bottom). 
For $\ell \in \{r,y,b\}$ each cycle in the decomposition is the union of the $A_1^\ell$, $A_2^\ell$, $A_3^\ell$, $B_1^\ell(J)$, $B_2^\ell(J)$, and $B_3^\ell(J)$ sets, i.e.

\[ A_1^\ell \cup A_2^\ell \cup A_3^\ell \cup
\left(\bigcup_{J=0}^{N-1} (B_1^\ell(J) \cup B_2^\ell(J) \cup B_3^\ell(J) )\right).\]

Since the $A$ block contains 24 edges of each color, and $B_1 \cup B_2 \cup B_3$ contains 16 edges of each color,  each cycle contains $24 + 16N$ edges.  Since $m=4$ and $n=9+6N$, the total number of edges in $C_m \sq C_n$ is $2mn= 2(4)(9+6N) = 72 + 48N$.  Thus each cycle contains one third of the total edges, as desired.

\begin{figure}
\begin{center}
\begin{tikzpicture}[scale=1, transform shape]
% A_1
    \draw[red, ultra thick, dashed] (0,.5) -- (0,-1) -- (1,-1) -- (1,0) -- (2,0) -- (2,.5);\draw[red, ultra thick, dashed] (-.5, -2) -- (2,-2) -- (2,-1) -- (2.5, -1);\draw[red, dashed, ultra thick] (-.5,-3) -- (0,-3) -- (0,-3.5);\draw[red, ultra thick, dashed] (2,-3.5) -- (2,-3) -- (2.5, -3);
    \draw[blue, densely dotted, ultra thick] (2.5, 0) -- (2,0) -- (2,-1) -- (1,-1) -- (1,-3) -- (2,-3) -- (2,-2) -- (2.5, -2);
    \draw[orange, ultra thick] (1, .5) -- (1,0) -- (-.5,0);\draw[orange, ultra thick] (-.5, -1) -- (0,-1) -- (0,-3) -- (1,-3) -- (1, -3.5);
    \foreach \x in {0,1,2}{\foreach\y in {0,-1,-2,-3}{\filldraw[fill=white,draw=black](\x,\y) circle (3pt);}}
    \foreach \y\i in {0,-1/1,-2/2,-3/3}{\node at (-1,\y){$\i$};}
    \node at (0,1){$2N$};
    \end{tikzpicture}\hspace{1cm}\begin{tikzpicture}
% A_2
    \draw[red, ultra thick, dashed] (3.5,.5) -- (3.5, 0) -- (4.5,0) -- (4.5, -2) -- (3.5,-2) -- (3.5,-1) -- (3,-1);\draw[red, ultra thick, dashed] (3,-3)--(3.5,-3) -- (3.5,-3.5);
    \draw[blue, densely dotted, ultra thick] (3,0) -- (3.5,0) -- (3.5,-1) -- (5.5,-1) -- (5.5,0) -- (6,0);\draw[blue, ultra thick, densely dotted] (3,-2) -- (3.5,-2) -- (3.5,-3) -- (5.5, -3) -- (5.5,-2) -- (6,-2);
    \draw[orange, ultra thick] (4.5,.5) -- (4.5,0) -- (5.5,0) -- (5.5,.5);\draw[orange, ultra thick] (4.5,-3.5) -- (4.5,-2) -- (5.5,-2) -- (5.5, -1) -- (6, -1);\draw[orange, ultra thick] (5.5, -3.5) -- (5.5,-3) -- (6, -3);
    \foreach \x in {3.5,4.5,5.5}{\foreach\y in {0,-1,-2,-3}{\filldraw[fill=white,draw=black](\x,\y) circle (3pt);}}
    \foreach \y\i in {0,-1/1,-2/2,-3/3}{\node at (2.5,\y){$\i$};}
    \node at (3.5,1){$4N+3$};
    \end{tikzpicture}\hspace{1cm}\begin{tikzpicture}
% A_3
    \draw[red, ultra thick, dashed] (8, .5) -- (8,0) -- (9,0) -- (9,-2) -- (9.5,-2);\draw[red, dashed, ultra thick] (9.5,-3) -- (8,-3) -- (8,-3.5);
    \draw[blue, densely dotted, ultra thick] (6.5, 0) -- (7,0) -- (7,-1) -- (8,-1) -- (8, -3) -- (7,-3) -- (7,-2) -- (6.5, -2);
    \draw[orange, ultra thick] (7, .5) -- (7, 0) -- (8,0) -- (8,-1) -- (9,-1) -- (9.5, -1);\draw[orange, ultra thick] (9,.5) -- (9,0) -- (9.5, 0);\draw[orange, ultra thick] (6.5, -1) -- (7, -1) -- (7, -2) -- (9, -2) -- (9, -3.5);\draw[orange, ultra thick] (6.5, -3) -- (7,-3) -- (7, -3.5);
    \foreach \x in {7,8,9}{\foreach\y in {0,-1,-2,-3}{\filldraw[fill=white,draw=black](\x,\y) circle (3pt);}}
    \foreach \y\i in {0,-1/1,-2/2,-3/3}{\node at (6,\y){$\i$};}
    \node at (7,1){$6N+6$};
\end{tikzpicture}

\vspace{.25in}

\begin{tikzpicture}[scale=1, transform shape]
% A_1
    \draw[red, ultra thick, dashed] (0,.5) -- (0,-1) -- (1,-1) -- (1,0) -- (2,0) -- (2,.5);\draw[red, ultra thick, dashed] (-.5, -2) -- (2,-2) -- (2,-1) -- (2.5, -1);\draw[red, dashed, ultra thick] (-.5,-3) -- (0,-3) -- (0,-3.5);\draw[red, ultra thick, dashed] (2,-3.5) -- (2,-3) -- (2.5, -3);
    \draw[blue, densely dotted, ultra thick] (2.5, 0) -- (2,0) -- (2,-1) -- (1,-1) -- (1,-3) -- (2,-3) -- (2,-2) -- (2.5, -2);
    \draw[orange, ultra thick] (1, .5) -- (1,0) -- (-.5,0);\draw[orange, ultra thick] (-.5, -1) -- (0,-1) -- (0,-3) -- (1,-3) -- (1, -3.5);

    \draw[red, ultra thick, dashed] (3,.5) -- (3, 0) -- (4,0) -- (4, -2) -- (3,-2) -- (3,-1) -- (2.5,-1);\draw[red, ultra thick, dashed] (2.5,-3)--(3,-3) -- (3,-3.5); \draw[blue, densely dotted, ultra thick] (2.5,0) -- (3,0) -- (3,-1) -- (5,-1) -- (5,0) -- (5.5,0);\draw[blue, ultra thick, densely dotted] (2.5,-2) -- (3,-2) -- (3,-3) -- (5, -3) -- (5,-2) -- (5.5,-2); \draw[orange, ultra thick] (4,.5) -- (4,0) -- (5,0) -- (5,.5);\draw[orange, ultra thick] (4,-3.5) -- (4,-2) -- (5,-2) -- (5, -1) -- (5.5, -1);\draw[orange, ultra thick] (5, -3.5) -- (5,-3) -- (5.5, -3);

    \draw[red, ultra thick, dashed] (7, .5) -- (7,0) -- (8,0) -- (8,-2) -- (8.5,-2);\draw[red, dashed, ultra thick] (8.5,-3) -- (7,-3) -- (7,-3.5); \draw[blue, densely dotted, ultra thick] (5.5, 0) -- (6,0) -- (6,-1) -- (7,-1) -- (7, -3) -- (6,-3) -- (6,-2) -- (5.5, -2); \draw[orange, ultra thick] (6, .5) -- (6, 0) -- (7,0) -- (7,-1) -- (8,-1) -- (8.5, -1);\draw[orange, ultra thick] (8,.5) -- (8,0) -- (8.5, 0);\draw[orange, ultra thick] (5.5, -1) -- (6, -1) -- (6, -2) -- (8, -2) -- (8, -3.5);\draw[orange, ultra thick] (5.5, -3) -- (6,-3) -- (6, -3.5);
    
    \foreach\x in {0, 1, ..., 8}{\foreach\y in {0, -1, -2, -3}{\filldraw[fill=white, draw=black] (\x, \y) circle (3pt);}}
    \foreach \x in {0,1,...,8}{\node at (\x,1){$\x$};}\foreach \y\i in {0/0,-1/1,-2/2,-3/3}{\node at (-1,\y){$\i$};}
        
\end{tikzpicture}

\vspace{.25in}

\begin{tikzpicture}[scale=.85, transform shape]
%B1
    \draw[red, ultra thick, dashed] (-1.5, -4) -- (-1, -4) -- (-1, -2) -- (0, -2) -- (0, -4) -- (.5, -4); \draw[red, ultra thick, dashed] (-1.5, -5) -- (.5, -5);
    \draw[orange, ultra thick] (-1.5, -2) -- (-1, -2) -- (-1, -1.5); \draw[orange, ultra thick] (.5, -2) -- (0, -2) -- (0, -1.5); \draw[orange, ultra thick] (-1.5, -3) -- (.5, -3); \draw[orange, ultra thick] (-1, -5.5) -- (-1, -4) -- (0, -4) -- (0, -5.5); \foreach \i in {-1, 0}{\foreach \j in {-2, -3, -4, -5}{\filldraw[fill=white, draw=black] (\i, \j) circle (3pt);}}
    \foreach \y\i in {-2/0,-3/1,-4/2,-5/3}{\node at (-2,\y){$\i$};}
    \foreach \x\j in {-1/0,0/1}{\node at (\x,-1){$\j$};}
    \end{tikzpicture}\hspace{1cm}\begin{tikzpicture}[scale=.85, transform shape]
%B2
    \draw[blue, ultra thick, densely dotted] (1.5, -2) -- (2, -2) -- (2, -1.5); \draw[blue, ultra thick, densely dotted] (3, -1.5) -- (3, -2) -- (3.5, -2); \draw[blue, ultra thick, densely dotted] (1.5, -4) -- (2, -4) -- (2, -3) -- (3, -3) -- (3, -4) -- (3.5, -4); \draw[blue, ultra thick, densely dotted] (2, -5.5) -- (2, -5) -- (3, -5) -- (3, -5.5); \draw[dashed, red, ultra thick] (1.5, -5) -- (2, -5) -- (2, -4) -- (3, -4) -- (3, -5) -- (3.5, -5);\draw[dashed, red, ultra thick] (1.5, -3) -- (2, -3) -- (2, -2) -- (3, -2) -- (3, -3) -- (3.5, -3);\foreach \i in {2, 3}{\foreach \j in {-2, -3,..., -5}{\filldraw[fill=white, draw=black] (\i, \j) circle (3pt);}}
    \foreach \y\i in {-2/0,-3/1,-4/2,-5/3}{\node at (1,\y){$\i$};}
    \node at (2,-1){$2N+3$};
    \end{tikzpicture}\hspace{1cm}\begin{tikzpicture}[scale=.85, transform shape]
%B3
    \draw[blue, ultra thick, densely dotted] (4.5, -2) -- (5, -2) -- (5, -1.5); \draw[blue, ultra thick, densely dotted] (6, -1.5) -- (6, -2) -- (6.5, -2); \draw[blue, ultra thick, densely dotted] (4.5, -4) -- (5, -4) -- (5, -3) -- (6, -3) -- (6, -4) -- (6.5, -4); \draw[blue, ultra thick, densely dotted] (5, -5.5) -- (5, -5) -- (6, -5) -- (6, -5.5); \draw[orange, ultra thick] (4.5, -5) -- (5, -5) -- (5, -4) -- (6, -4) -- (6, -5) -- (6.5, -5);\draw[orange, ultra thick] (4.5, -3) -- (5, -3) -- (5, -2) -- (6, -2) -- (6, -3) -- (6.5, -3);\foreach \i in {5, 6}{\foreach \j in {-2, -3,..., -5}{\filldraw[fill=white, draw=black] (\i, \j) circle (3pt);}}
    \foreach \y\i in {-2/0,-3/1,-4/2,-5/3}{\node at (4,\y){$\i$};}
    \node at (5,-1){$4N+6$};
\end{tikzpicture}

\vspace{.25in}

\begin{tikzpicture}[scale=.57, transform shape]
%B1(0)
    \draw[red, ultra thick, dashed] (-1.5, -2) -- (-1, -2) -- (-1, 0) -- (0, 0) -- (0, -2) -- (.5, -2); \draw[red, ultra thick, dashed] (-1.5, -3) -- (.5, -3);
    \draw[orange, ultra thick] (-1.5, 0) -- (-1, 0) -- (-1, .5); \draw[orange, ultra thick] (.5, 0) -- (0, 0) -- (0, .5); \draw[orange, ultra thick] (-1.5, -1) -- (.5, -1); \draw[orange, ultra thick] (-1, -3.5) -- (-1, -2) -- (0, -2) -- (0, -3.5); \foreach \i in {-1, 0}{\foreach \j in {0, -1, -2, -3}{\filldraw[fill=white, draw=black] (\i, \j) circle (3pt);}}
\node at (1, -1.5){$\cdots$};
%B1(N-1)
    \draw[red, ultra thick, dashed] (1.5, -2) -- (2, -2) -- (2, 0) -- (3, 0) -- (3, -2) -- (3.5, -2); \draw[red, ultra thick, dashed] (1.5, -3) -- (3.5, -3);
    \draw[orange, ultra thick] (1.5, 0) -- (2, 0) -- (2, .5); \draw[orange, ultra thick] (3.5, 0) -- (3, 0) -- (3, .5); \draw[orange, ultra thick] (1.5, -1) -- (3.5, -1); \draw[orange, ultra thick] (2, -3.5) -- (2, -2) -- (3, -2) -- (3, -3.5); \foreach \i in {2, 3}{\foreach \j in {0, -1, -2, -3}{\filldraw[fill=white, draw=black] (\i, \j) circle (3pt);}}
% A_1
    \draw[red, ultra thick, dashed] (4.5,.5) -- (4.5,-1) -- (5.5,-1) -- (5.5,0) -- (6.5,0) -- (6.5,.5);\draw[red, ultra thick, dashed] (4, -2) -- (6.5,-2) -- (6.5,-1) -- (7, -1);\draw[red, dashed, ultra thick] (4,-3) -- (4.5,-3) -- (4.5,-3.5);\draw[red, ultra thick, dashed] (6.5,-3.5) -- (6.5,-3) -- (7, -3);
    \draw[blue, densely dotted, ultra thick] (7, 0) -- (6.5,0) -- (6.5,-1) -- (5.5,-1) -- (5.5,-3) -- (6.5,-3) -- (6.5,-2) -- (7, -2);
    \draw[orange, ultra thick] (5.5, .5) -- (5.5,0) -- (4,0);\draw[orange, ultra thick] (4, -1) -- (4.5,-1) -- (4.5,-3) -- (5.5,-3) -- (5.5, -3.5);
    \foreach \i in {4.5, 5.5, 6.5}{\foreach \j in {0, -1, -2, -3}{\filldraw[fill=white, draw=black] (\i, \j) circle (3pt);}}
%B2(0)
    \draw[blue, ultra thick, densely dotted] (7.5, 0) -- (8, 0) -- (8, .5); \draw[blue, ultra thick, densely dotted] (9, .5) -- (9, 0) -- (9.5, 0); \draw[blue, ultra thick, densely dotted] (7.5, -2) -- (8, -2) -- (8, -1) -- (9, -1) -- (9, -2) -- (9.5, -2); \draw[blue, ultra thick, densely dotted] (8, -3.5) -- (8, -3) -- (9, -3) -- (9, -3.5); \draw[dashed, red, ultra thick] (7.5, -3) -- (8, -3) -- (8, -2) -- (9, -2) -- (9, -3) -- (9.5, -3);\draw[dashed, red, ultra thick] (7.5, -1) -- (8, -1) -- (8, 0) -- (9, 0) -- (9, -1) -- (9.5, -1);\foreach \i in {8, 9}{\foreach \j in {0, -1, -2, -3}{\filldraw[fill=white, draw=black] (\i, \j) circle (3pt);}}
\node at (10, -1.5){$\cdots$};
%B2(N-1)
    \draw[blue, ultra thick, densely dotted] (10.5, 0) -- (11, 0) -- (11, .5); \draw[blue, ultra thick, densely dotted] (12, .5) -- (12, 0) -- (12.5, 0); \draw[blue, ultra thick, densely dotted] (10.5, -2) -- (11, -2) -- (11, -1) -- (12, -1) -- (12, -2) -- (12.5, -2); \draw[blue, ultra thick, densely dotted] (11, -3.5) -- (11, -3) -- (12, -3) -- (12, -3.5); \draw[dashed, red, ultra thick] (10.5, -3) -- (11, -3) -- (11, -2) -- (12, -2) -- (12, -3) -- (12.5, -3);\draw[dashed, red, ultra thick] (10.5, -1) -- (11, -1) -- (11, 0) -- (12, 0) -- (12, -1) -- (12.5, -1);\foreach \i in {11, 12}{\foreach \j in {0, -1, -2, -3}{\filldraw[fill=white, draw=black] (\i, \j) circle (3pt);}}
% A_2
    \draw[red, ultra thick, dashed] (13.5,.5) -- (13.5, 0) -- (14.5,0) -- (14.5, -2) -- (13.5,-2) -- (13.5,-1) -- (13,-1);\draw[red, ultra thick, dashed] (13,-3)--(13.5,-3) -- (13.5,-3.5);
    \draw[blue, densely dotted, ultra thick] (13,0) -- (13.5,0) -- (13.5,-1) -- (15.5,-1) -- (15.5,0) -- (16,0);\draw[blue, ultra thick, densely dotted] (13,-2) -- (13.5,-2) -- (13.5,-3) -- (15.5, -3) -- (15.5,-2) -- (16,-2); \draw[orange, ultra thick] (14.5,.5) -- (14.5,0) -- (15.5,0) -- (15.5,.5);\draw[orange, ultra thick] (14.5,-3.5) -- (14.5,-2) -- (15.5,-2) -- (15.5, -1) -- (16, -1);\draw[orange, ultra thick] (15.5, -3.5) -- (15.5,-3) -- (16, -3);\foreach \i in {13.5, 14.5, 15.5}{\foreach \j in {0, -1, -2, -3}{\filldraw[fill=white, draw=black] (\i, \j) circle (3pt);}}
%B3(0)
    \draw[blue, ultra thick, densely dotted] (16.5, 0) -- (17, 0) -- (17, .5); \draw[blue, ultra thick, densely dotted] (18, .5) -- (18, 0) -- (18.5, 0); \draw[blue, ultra thick, densely dotted] (16.5, -2) -- (17, -2) -- (17, -1) -- (18, -1) -- (18, -2) -- (18.5, -2); \draw[blue, ultra thick, densely dotted] (17, -3.5) -- (17, -3) -- (18, -3) -- (18, -3.5); \draw[orange, ultra thick] (16.5, -3) -- (17, -3) -- (17, -2) -- (18, -2) -- (18, -3) -- (18.5, -3);\draw[orange, ultra thick] (16.5, -1) -- (17, -1) -- (17, 0) -- (18, 0) -- (18, -1) -- (18.5, -1);\foreach \i in {17, 18}{\foreach \j in {0, -1, -2, -3}{\filldraw[fill=white, draw=black] (\i, \j) circle (3pt);}}
\node at (19, -1.5){$\cdots$};
%B3(N-1)
    \draw[blue, ultra thick, densely dotted] (19.5, 0) -- (20, 0) -- (20, .5); \draw[blue, ultra thick, densely dotted] (21, .5) -- (21, 0) -- (21.5, 0); \draw[blue, ultra thick, densely dotted] (19.5, -2) -- (20, -2) -- (20, -1) -- (21, -1) -- (21, -2) -- (21.5, -2); \draw[blue, ultra thick, densely dotted] (20, -3.5) -- (20, -3) -- (21, -3) -- (21, -3.5); \draw[orange, ultra thick] (19.5, -3) -- (20, -3) -- (20, -2) -- (21, -2) -- (21, -3) -- (21.5, -3);\draw[orange, ultra thick] (19.5, -1) -- (20, -1) -- (20, 0) -- (21, 0) -- (21, -1) -- (21.5, -1);\foreach \i in {20, 21}{\foreach \j in {0, -1, -2, -3}{\filldraw[fill=white, draw=black] (\i, \j) circle (3pt);}}
% A_3
    \draw[red, ultra thick, dashed] (23.5, .5) -- (23.5,0) -- (24.5,0) -- (24.5,-2) -- (25,-2);\draw[red, dashed, ultra thick] (25,-3) -- (23.5,-3) -- (23.5,-3.5);
    \draw[blue, densely dotted, ultra thick] (22, 0) -- (22.5,0) -- (22.5,-1) -- (23.5,-1) -- (23.5, -3) -- (22.5,-3) -- (22.5,-2) -- (22, -2);
    \draw[orange, ultra thick] (22.5, .5) -- (22.5, 0) -- (23.5,0) -- (23.5,-1) -- (24.5,-1) -- (25, -1);\draw[orange, ultra thick] (24.5,.5) -- (24.5,0) -- (25, 0);\draw[orange, ultra thick] (22, -1) -- (22.5, -1) -- (22.5, -2) -- (24.5, -2) -- (24.5, -3.5);\draw[orange, ultra thick] (22, -3) -- (22.5,-3) -- (22.5, -3.5);
    \foreach \x in {22.5, 23.5, 24.5}{\foreach \y in {0, -1, -2, -3}{\filldraw[fill=white, draw = black] (\x, \y) circle (3pt);}}

\large{\node at (-.75, -5){$B_1(0)$}; \node at (1, -5){$\cdots$}; \node at (2.75, -5){$B_1(N-1)$}; \node at (5.5, -5){$A_1$};\node at (8.25, -5){$B_2(0)$};\node at (10, -5){$\cdots$};\node at (11.75, -5){$B_2(N-1)$};\node at (14.5, -5){$A_2$};\node at (17.25, -5){$B_3(0)$};\node at (19, -5){$\cdots$};\node at (20.75, -5){$B_3(N-1)$};\node at (23.5, -5){$A_3$};}

\foreach\x\n in {-1/0, 0/1, 2/$2(N-1)$, 4.5/$2N$, 8/$2N+3$, 11/$4N+1$, 13.5/$4N+3$, 17/$4N+6$, 20/$6N+4$, 22.5/$6N+6$, 24.5/$n-1$}{\node at (\x,1){\n};}
\foreach\y\m in {0/0, -1/1, -2/2, -3/3}{\node at (-2,\y){\m};}

\node at (24.5, 1.5){$6N+8$};

\end{tikzpicture}

\end{center}
\caption{Top (left to right): The blocks $A_1$, $A_2$, and $A_3,$ for Case 4. Above middle: The (minimal) $A$ block for Case 4. Below middle: The blocks $B_1$, $B_2$, and $B_3,$ for Case 4. Bottom: Decomposition of $C_m\sq C_n$, with the arrangements of $A_1$, $A_2$, $A_3$, $B_1$, $B_2$, and $B_3$ blocks in $C_m \sq C_n$ in Case 4.
}\label{ACase4}
\end{figure}
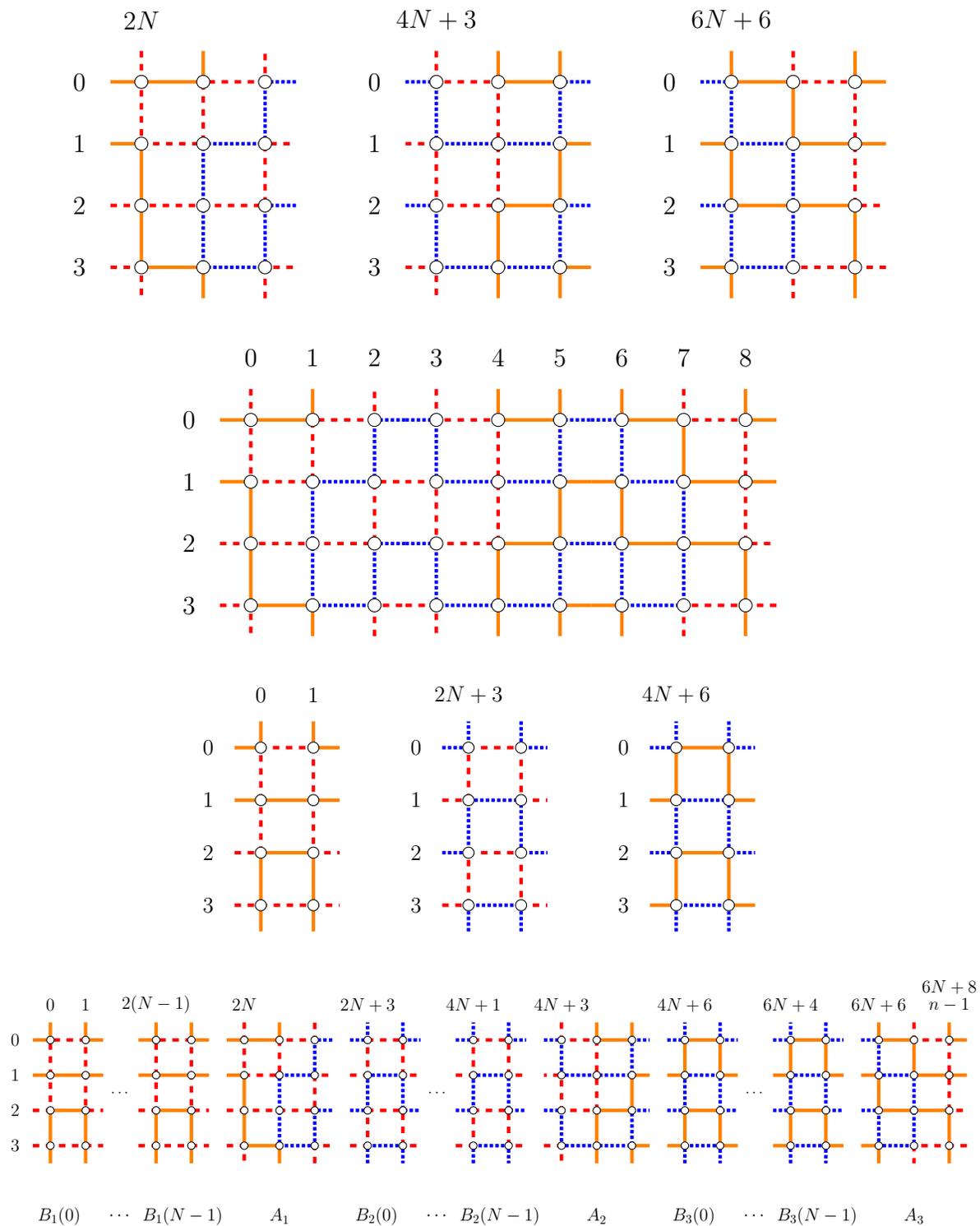

\paragraph{Case 5 ($m\geq 6$ is even and $n = 6$):} Assume $m = 6+2M$ for a nonnegative integer $M$ and $n=6$. 
The decomposition of $C_m \sq C_n$ is illustrated in Figure~\ref{Case5} and contains one $6 \times 6$ $A$ block and $M$ $2 \times 6$ $B$ blocks, shown in Figure~\ref{Case5} (left) and defined below.  The $A$ block is the minimal block for this case.

For the $6 \times 6$ $A$ block define three edge sets $A^r$, $A^y$, and $A^b$, each of which forms a cycle with 24 edges, shown as dashed red, solid yellow, and dotted blue, respectively, in Figure~\ref{Case5} (left above).

$A^r = \{(1,-1)(1,0), (1,0)(2,0), (2,0)(2,1), (2,1)(1,1), (1,1)(1,2), (1,2)(2,2),\\ (2,2)(3,2), (3,2)(4,2), (4,2)(4,3), (4,3)(5,3), (5,3)(5,2), (5,2)(6,2), (-1,2)(0,2),\\ (0,2)(0,3), (0,3)(1,3), (1,3)(2,3), (2,3)(2,4), (2,4)(3,4), (3,4)(4,4), (4,4)(4,5), \\(4,5)(5,5), (5,5)(5,4), (5,4)(6,4), (-1,4)(0,4), (0,4)(0,5), (0,5)(1,5), (1,5)(1,6)\}$

$A^y = \{(3,-1)(3,0), (3,0)(4,0), (4,0)(5,0), (5,0)(6,0), (-1,0)(0,0), (0,0)(1,0), \\(1,0)(1,1), (1,1)(0,1), (0,1)(-1,1), (6,1)(5,1), (5,1)(4,1), (4,1)(3,1), (3,1)(2,1),\\ (2,1)(2,2), (2,2)(2,3), (2,3)(3,3), (3,3)(4,3), (4,3)(4,4), (4,4)(5,4), (5,4)(5,3),\\ (5,3)(6,3), (-1,3)(0,3), (0,3)(0,4), (0,4)(1,4), (1,4)(2,4), (2,4)(2,5), (2,5)(3,5),\\ (3,5)(3,6)\}$

$A^b = \{(0,-1)(0,0), (0,0)(0,1), (0,1)(0,2), (0,2)(1,2), (1,2)(1,3), (1,3)(1,4), \\(1,4)(1,5), (1,5)(2,5), (2,5)(2,6), (2,-1)(2,0), (2,0)(3,0), (3,0)(3,1), (3,1)(3,2),\\ (3,2)(3,3), (3,3)(3,4), (3,4)(3,5), (3,5)(4,5), (4,5)(4,6), (4,-1)(4,0), (4,0)(4,1),\\ (4,1)(4,2), (4,2)(5,2), (5,2)(5,1), (5,1)(5,0), (5,0)(5,-1), (5,6)(5,5), (5,5)(6,5)\}$

For the $2 \times 6$ $B$ block define edge sets $B^r$, $B^y$, and $B^b$, each containing eight edges, shown as dashed red, solid yellow, and dotted blue, respectively in Figure~\ref{Case5} (left below). 

$B^r = \{(-1,2)(0,2), (0,2)(0,3), (0,3)(1,3), (1,3)(1,2), (1,2)(2,2), (-1,4)(0,4),\\ (0,4)(0,5), (0,5)(1,5), (1,5)(1,4), (1,4)(2,4)\}$

$B^y = \{(-1,0)(0,0), (0,0)(1,0), (1,0)(2,0), (-1,1)(0,1), (0,1)(1,1), (1,1)(2,1), \\(-1,3)(0,3), (0,3)(0,4), (0,4)(1,4), (1,4)(1,3), (1,3)(2,3)\}$

$B^b = \{(-1,5)(0,5), (0,5)(0,6), (0,-1)(0,0), (0,0)(0,1), (0,1)(0,2), (0,2)(1,2), \\(1,2)(1,1), (1,1)(1,0), (1,0)(1,-1), (1,6)(1,5), (1,5)(2,5)\}$

For $0\leq I \leq M-1$, let the block $B(I)$ to be the set of edges incident to the vertices $(i, j)$ where $6+2I \leq i\leq 7+2I$ and $0\leq j \leq 5$.
To color the edges in $B(I)$, define three edge sets $B^r(I)$, $B^y(I)$, and $B^b(I)$, where for $\ell \in \{r,y,b\}$, edge $(i,j)(i',j')$ is in $B^\ell$ if and only if $(i+6 + 2I,j)(i'+6+2I,j')$ is in $B^\ell(I)$. For example, $(8,2)(9,2)$ is blue (in $B^b(1)$) because $(0,2)(1,2)$ is blue (in $B^b$).

For $\ell \in \{r,y,b\}$ each cycle in the decomposition is the union of the $A^\ell$ and $B^\ell(I)$ sets, i.e.

\[ A^\ell \cup
\left(\bigcup_{I=0}^{M-1}B^\ell(I) \right).\]

Since the $A$ block contains 24 edges of each color, and each of the $M$ $B$ blocks contain eight edges of each color, each cycle contains $24 + 8M$ edges.  Since $m=6+2M$ and $n=6$, the total number of edges in $C_m \sq C_n$ is $2mn= 2(6+2M)(6) = 72 + 24M$.  Thus each cycle contains one third of the total edges, as desired.

\begin{figure}
\begin{center}
\begin{tikzpicture}[scale=1, transform shape]
% A 
    \draw[blue, ultra thick, densely dotted] (-.5, 0) -- (2, 0) -- (2, -1) -- (5, -1) -- (5, -2) -- (5.5, -2); \draw[blue, ultra thick, densely dotted] (-.5, -2) -- (0, -2) -- (0, -3) -- (5, -3) -- (5, -4) -- (5.5, -4); \draw[blue, ultra thick, densely dotted] (5, .5) -- (5, 0) -- (5.5, 0); \draw[blue, ultra thick, densely dotted] (5.5, -5) -- (5, -5) -- (5, -5.5); \draw[blue, ultra thick, densely dotted] (-.5, -4) -- (2, -4) -- (2, -5) -- (-.5, -5);
    \draw[orange, ultra thick] (0, .5) -- (0, -1) -- (1, -1) -- (1, .5); \draw[orange, ultra thick] (3, .5) -- (3, 0) -- (4, 0) -- (4, -2) -- (5, -2) -- (5, -3) -- (5.5, -3); \draw[orange, ultra thick] (-.5, -3) -- (0, -3) -- (0, -5.5); \draw[orange, ultra thick] (1, -5.5)  -- (1, -2) -- (3, -2) -- (3, -4) -- (4, -4) -- (4, -5) -- (3, -5) -- (3, -5.5);
    \draw[red, ultra thick, dashed] (-.5, -1) -- (0, -1) -- (0, -2) -- (1, -2) -- (1, -1) -- (2, -1) -- (2, -4) -- (3, -4) -- (3, -5) -- (2, -5) -- (2, -5.5); \draw[dashed, red, ultra thick] (2, .5) -- (2, 0) -- (3, 0) -- (3, -2) -- (4, -2) -- (4, -4) -- (5, -4) -- (5, -5) -- (4, -5) -- (4, -5.5); \draw[dashed, red, ultra thick] (4, .5) -- (4, 0) -- (5, 0) -- (5, -1) -- (5.5, -1);
    \foreach \x in {0, 1, ..., 5}{\foreach \y in {0, -1, ..., -5}{\filldraw[draw=black, fill=white] (\x, \y) circle (3pt);}}
% C_0(0)
    \draw[red, dashed, ultra thick] (2,-7) -- (2,-7.5) -- (3,-7.5) -- (3,-8.5) -- (2,-8.5) -- (2,-9); \draw[red, dashed, ultra thick] (4,-7) -- (4,-7.5) -- (5,-7.5) -- (5,-8.5) -- (4,-8.5) -- (4,-9); 
    \draw[blue, densely dotted, ultra thick] (5,-7) -- (5,-7.5) -- (5.5,-7.5); \draw[blue, densely dotted, ultra thick] (-.5,-7.5) -- (2,-7.5) -- (2,-8.5) -- (-.5,-8.5); \draw[blue, densely dotted, ultra thick] (5.5,-8.5) -- (5,-8.5) -- (5,-9); 
    \draw[orange, ultra thick] (0,-7) -- (0,-9); \draw[orange, ultra thick] (1,-7) -- (1,-9);
    \draw[orange, ultra thick] (3,-7) -- (3,-7.5) -- (4,-7.5) -- (4,-8.5) -- (3,-8.5) -- (3,-9);
    \foreach \x in {0,1,...,5}{\foreach\y in {-7.5,-8.5}{\filldraw[draw=black, fill=white] (\x, \y) circle (3pt);}}
    \foreach \x in {0,1,2,3,4,5}{\foreach\y in {1,-6.5}{\node at (\x,\y){$\x$};}}\foreach \y\i in {0/0,-1/1,-2/2,-3/3, -4/4, -5/5,-7.5/0,-8.5/1}{\node at (-1,\y){$\i$};}
\end{tikzpicture}
\hfill
\begin{tikzpicture}[scale=.6, transform shape]
% A 
    \draw[blue, ultra thick, densely dotted] (-.5, 0) -- (2, 0) -- (2, -1) -- (5, -1) -- (5, -2) -- (5.5, -2); \draw[blue, ultra thick, densely dotted] (-.5, -2) -- (0, -2) -- (0, -3) -- (5, -3) -- (5, -4) -- (5.5, -4); \draw[blue, ultra thick, densely dotted] (5, .5) -- (5, 0) -- (5.5, 0); \draw[blue, ultra thick, densely dotted] (5.5, -5) -- (5, -5) -- (5, -5.5); \draw[blue, ultra thick, densely dotted] (-.5, -4) -- (2, -4) -- (2, -5) -- (-.5, -5);
    \draw[orange, ultra thick] (0, .5) -- (0, -1) -- (1, -1) -- (1, .5); \draw[orange, ultra thick] (3, .5) -- (3, 0) -- (4, 0) -- (4, -2) -- (5, -2) -- (5, -3) -- (5.5, -3); \draw[orange, ultra thick] (-.5, -3) -- (0, -3) -- (0, -5.5); \draw[orange, ultra thick] (1, -5.5)  -- (1, -2) -- (3, -2) -- (3, -4) -- (4, -4) -- (4, -5) -- (3, -5) -- (3, -5.5);
    \draw[red, ultra thick, dashed] (-.5, -1) -- (0, -1) -- (0, -2) -- (1, -2) -- (1, -1) -- (2, -1) -- (2, -4) -- (3, -4) -- (3, -5) -- (2, -5) -- (2, -5.5); \draw[dashed, red, ultra thick] (2, .5) -- (2, 0) -- (3, 0) -- (3, -2) -- (4, -2) -- (4, -4) -- (5, -4) -- (5, -5) -- (4, -5) -- (4, -5.5); \draw[dashed, red, ultra thick] (4, .5) -- (4, 0) -- (5, 0) -- (5, -1) -- (5.5, -1);
    \foreach \x in {0, 1, ..., 5}{\foreach \y in {0, -1, ..., -5}{\filldraw[draw=black, fill=white] (\x, \y) circle (3pt);}}
% C_0(0)
    \draw[red, dashed, ultra thick] (2,-6) -- (2,-6.5) -- (3,-6.5) -- (3,-7.5) -- (2,-7.5) -- (2,-8); \draw[red, dashed, ultra thick] (4,-6) -- (4,-6.5) -- (5,-6.5) -- (5,-7.5) -- (4,-7.5) -- (4,-8); 
    \draw[blue, densely dotted, ultra thick] (5,-6) -- (5,-6.5) -- (5.5,-6.5); \draw[blue, densely dotted, ultra thick] (-.5,-6.5) -- (2,-6.5) -- (2,-7.5) -- (-.5,-7.5); \draw[blue, densely dotted, ultra thick] (5.5,-7.5) -- (5,-7.5) -- (5,-8); \draw[orange, ultra thick] (0,-6) -- (0,-8); \draw[orange, ultra thick] (1,-6) -- (1,-8);\draw[orange, ultra thick] (3,-6) -- (3,-6.5) -- (4,-6.5) -- (4,-7.5) -- (3,-7.5) -- (3,-8);\foreach \x in {0,1,...,5}{\foreach\y in {-6.5,-7.5}{\filldraw[draw=black, fill=white] (\x, \y) circle (3pt);}}
% C_0(1)
    \draw[red, dashed, ultra thick] (2,-8.5) -- (2,-9) -- (3,-9) -- (3,-10) -- (2,-10) -- (2,-10.5); \draw[red, dashed, ultra thick] (4,-8.5) -- (4,-9) -- (5,-9) -- (5,-10) -- (4,-10) -- (4,-10.5); 
    \draw[blue, densely dotted, ultra thick] (5,-8.5) -- (5,-9) -- (5.5,-9); \draw[blue, densely dotted, ultra thick] (-.5,-9) -- (2,-9) -- (2,-10) -- (-.5,-10); \draw[blue, densely dotted, ultra thick] (5.5,-10) -- (5,-10) -- (5,-10.5); \draw[orange, ultra thick] (0,-8.5) -- (0,-10.5); \draw[orange, ultra thick] (1,-8.5) -- (1,-10.5);\draw[orange, ultra thick] (3,-8.5) -- (3,-9) -- (4,-9) -- (4,-10) -- (3,-10) -- (3,-10.5);\foreach \x in {0,1,...,5}{\foreach\y in {-9,-10}{\filldraw[draw=black, fill=white] (\x, \y) circle (3pt);}}
\node at (2.5,-11.5){$\vdots$};
% C_0(M-1)
    \draw[red, dashed, ultra thick] (2,-12.5) -- (2,-13) -- (3,-13) -- (3,-14) -- (2,-14) -- (2,-14.5); \draw[red, dashed, ultra thick] (4,-12.5) -- (4,-13) -- (5,-13) -- (5,-14) -- (4,-14) -- (4,-14.5); 
    \draw[blue, densely dotted, ultra thick] (5,-12.5) -- (5,-13) -- (5.5,-13); \draw[blue, densely dotted, ultra thick] (-.5,-13) -- (2,-13) -- (2,-14) -- (-.5,-14); \draw[blue, densely dotted, ultra thick] (5.5,-14) -- (5,-14) -- (5,-14.5); \draw[orange, ultra thick] (0,-12.5) -- (0,-14.5); \draw[orange, ultra thick] (1,-12.5) -- (1,-14.5);\draw[orange, ultra thick] (3,-12.5) -- (3,-13) -- (4,-13) -- (4,-14) -- (3,-14) -- (3,-14.5);\foreach \x in {0,1,...,5}{\foreach\y in {-13,-14}{\filldraw[draw=black, fill=white] (\x, \y) circle (3pt);}}
\node at (8, -2.5){\Large$A$}; \node at (8,-7) {\Large$B(0)$}; \node at (8,-9.5){\Large$B(1)$};\node at (8,-11.5){$\vdots$};\node at (8, -13.5){\Large$B(M-1)$};
% vertices
\foreach \x in {0, 1, ..., 5}{\node at (\x, 1.5){\large\x};}\foreach \y\j in {0/0, -1/1, -2/2, -3/3, -4/4, -5/5, -6.5/6, -7.5/7, -9/8, -10/9, -13/$m-2$, -14/$m-1$}{\node at (-1.5, \y){\large\j};}
\end{tikzpicture}
\end{center}
\caption{Left above: The $A$ block for Cases 5 and 6, with $A^r, A^y, A^b$ in dashed red, solid yellow, and dotted blue. Left below: The $B$ block for Case 5 with $B^r, B^y, B^b$ in dashed red, solid yellow, and dotted blue. Right: The decomposition of $C_m \sq C_6$ into three cycles, and the arrangement of $A$ and $B$ blocks in $C_m\sq C_6$ in Case 5.}\label{Case5}
\end{figure}
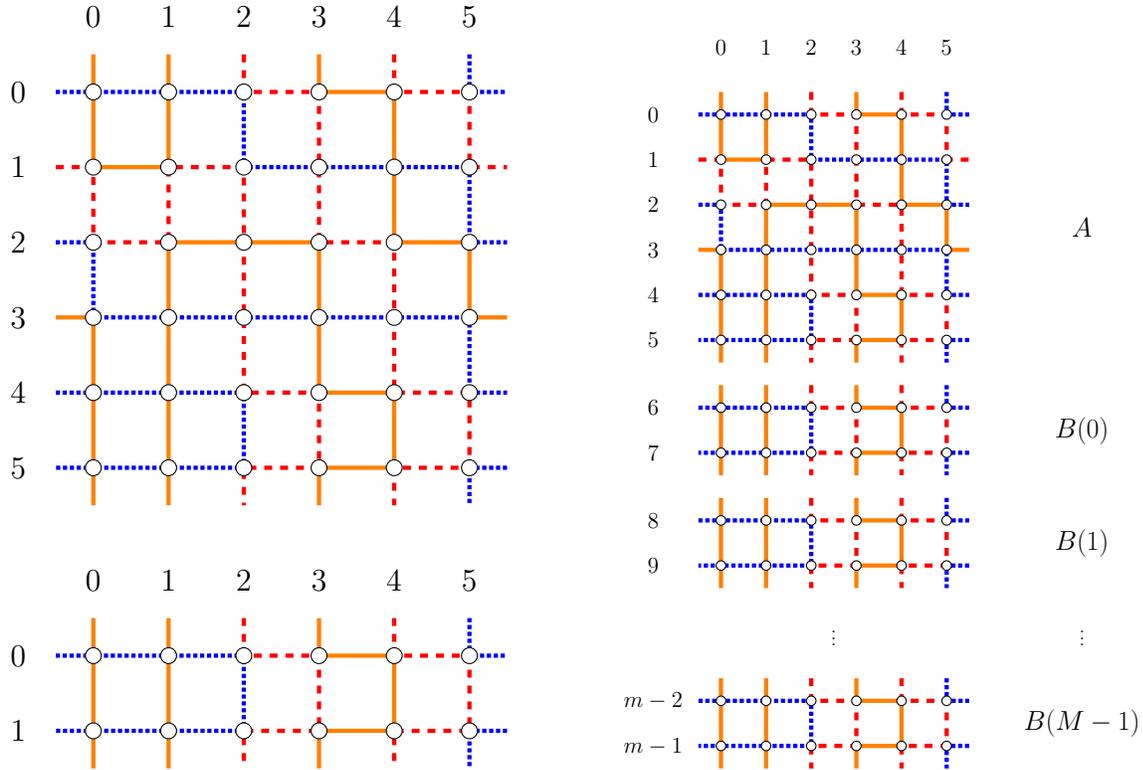

\paragraph{Case 6 ($m \geq 6$ is even and $n \ge 12$ is even):} Assume $m = 6+2M$ and $n=6+6N$ for a nonnegative integer $M$ and a positive integer $N$.  The decomposition of $C_m \sq C_n$ is illustrated in Figure~\ref{ABCCase6} and contains one $6 \times 6$ $A$ block, $N$ $6 \times 6$ $B$ blocks, and $M$ $2 \times n$ $C$ blocks, shown in Figures~\ref{BCase6} and \ref{CCase6} and defined below.  See also Figure~\ref{Ex2}.

The $A$ block and the three edge sets $A^r$, $A^y$, and $A^b$ in this case are identical to the $A$ block from Case 5, so the three edge sets $A^r$, $A^y$, and $A^b$ each form a cycle with 24 edges.

For the $6 \times 6$ $B$ block define three edge sets $B^r, B^y, B^b$, each containing 24 edges, shown as dashed red, solid yellow, and dotted blue respectively in Figure~\ref{BCase6}. 

$B^r = \{(-1,0)(0,0), (0,0)(1,0), (1,0)(1,-1), (1,6)(1,5), (1,5)(2,5), (2,5)(3,5),\\ (3,5)(3,4), (3,4)(4,4), (4,4)(5,4), (5,4)(6,4), (-1,4)(0,4), (0,4)(1,4), (1,4)(1,3),\\ (1,3)(2,3), (2,3)(3,3), (3,3)(3,2), (3,2)(4,2), (4,2)(5,2), (5,2)(6,2), (-1,2)(0,2), \\(0,2)(1,2), (1,2)(1,1), (1,1)(2,1), (2,1)(3,1), (3,1)(3,0), (3,0)(4,0), (4,0)(5,0),\\ (5,0)(6,0)\}.$

$B^y = \{(-1,1)(0,1), (0,1)(1,1), (1,1)(1,0), (1,0)(2,0), (2,0)(3,0), (3,0)(3,-1),\\ (3,6)(3,5), (3,5)(4,5), (4,5)(5,5), (5,5)(6,5), (-1,5)(0,5), (0,5)(1,5), (1,5)(1,4),\\ (1,4)(2,4), (2,4)(3,4), (3,4)(3,3), (3,3)(4,3), (4,3)(5,3), (5,3)(6,3), (-1,3)(0,3), \\(0,3)(1,3), (1,3)(1,2), (1,2)(2,2), (2,2)(3,2), (3,2)(3,1), (3,1)(4,1), (4,1)(5,1),\\ (5,1)(6,1)\}.$

$B^b = \{(0,-1)(0,0), (0,0)(0,1), (0,1)(0,2), (0,2)(0,3), (0,3)(0,4), (0,4)(0,5),\\ (0,5)(0,6), (2,-1)(2,0), (2,0)(2,1), (2,1)(2,2), (2,2)(2,3), (2,3)(2,4), (2,4)(2,5),\\ (2,5)(2,6),(4,-1)(4,0), (4,0)(4,1), (4,1)(4,2), (4,2)(4,3), (4,3)(4,4), (4,4)(4,5), \\(4,5)(4,6),(5,-1)(5,0), (5,0)(5,1), (5,1)(5,2), (5,2)(5,3), (5,3)(5,4), (5,4)(5,5),\\ (5,5)(5,6)\}.$

For $0\leq J\leq N-1,$, let the block $B(J)$ in $C_m \sq C_n$ be the set of edges incident to the vertices $(i,j)$ where $0 \leq i \leq 6$ and $6+6J\leq j\leq 11+6J$. 
To color the edges in $B(J)$, define three edge sets $B^r(J)$, $B^y(J)$, and $B^b(J)$, where for $\ell \in \{r,y,b\}$, edge $(i,j)(i',j')$ is in $B^\ell$ if and only if $(i,j+ 6 +6J)(i',j'+6 +6J)$ is in $B^\ell(J)$. For example, $(3,10)(3,11)$ is red (in $B^r(0)$) because $(3,4)(3,5)$ is red (in $B^r$).

The minimal block for this case has $m=6$ and $n=12,$ so it is the union of $A$ and $B(0)$.

For the $2 \times n$ $C$ block define edge sets $C^r$, $C^y$, and $C^b$, each containing $4n/3 = 8N+8$ edges, shown as dashed red, solid yellow, and dotted blue respectively in Figure~\ref{CCase6}. Since the blocks have variable size, we describe each edge set using set-builder notation rather than listing each edge explicitly.

$C^r = \{(-1,2)(0,2), (0,2)(1,2), (1,2)(2,2), (-1,4)(0,4), (0,4)(1,4), (1,4)(2,4), \\(-1,6+2x_1)(0, 6+2x_1), (0,6+2x_1)(0, 5+2x_1), (0,5+2x_1)(1, 5+2x_1), \\(1,5+2x_1)(1,6+2x_1), (1,6+2x_1)(2, 6+2x_1), (-1,n-2x_2)(0,n-2x_2), \\(0,n-2x_2)(1, n-2x_2),(1,n-2x_2)(2,n-2x_2) : 0\leq x_1\leq N+1, 1\leq x_2\leq 2N-2\}$.

$C^y = \{(-1, 0)(0,0), (0,0)(1,0), (1,0)(2,0), (-1, 1)(0,1), (0,1)(1,1), (1,1)(2,1),\\(-1, 3)(0,3), (0,3)(1,3), (1,3)(2,3), (-1,7+2x_3)(0,7+2x_3), (0,7+2x_3)(0,6+2x_3)\\ (0,6+2x_3)(1,6+2x_3),(1,6+2x_3)(1,7+2x_3), (1,7+2x_3)(2,7+2x_3),\\ (-1,n-(2x_4+1))(0,n-(2x_4+1)), (0,n-(2x_4+1))(1,n-(2x_4+1)), \\(1,n-(2x_4+1))(2,n-(2x_4+1)) : 0 \leq x_3\leq N, 0\leq x_4\leq 2N-2\}$.

$C^b = \{(-1,5)(0,5), (0,5)(0,4), (0,4)(0,3), (0,3)(0,2), (0,2)(0,1), (0,1)(0,0),\\ (0,0)(0,-1), (0,n-x_5)(0,n-(x_5+1)), (0,n-(4N-2))(1,n-(4N-2)),\\ (1,n-(x_6-1))(1,n-x_6),(2,5)(1,5), (1,5)(1,4), (1,4)(1,3), (1,3)(1,2), \\(1,2)(1,1), (1,1)(1,0), (1,0)(1,-1): 0\leq x_5\leq 4N-2, 1\leq x_6\leq 4N-2\}$

For $0\leq I\leq M-1$, let the block $C(I)$ in $C_m \sq C_n$ be the set of edges incident to the vertices $(i,j)$ where $6+2I \leq i \leq 7+2I$ and $0\leq j\leq n-1$. 
To color the edges in $C(I)$, define three edge sets $C^r(I)$, $C^y(I)$, and $C^b(I)$, where for $\ell \in \{r,y,b\}$, edge $(i,j)(i',j')$ is in $C^\ell$ if and only if $(i + 6 + 2I,j)(i'+6 + 2I,j')$ is in $C^\ell(J)$. For example, $(8,3)(9,3)$ is yellow (in $C^y(1)$) because $(0,3)(1,3)$ is yellow (in $C^y$).

For $\ell \in \{r,y,b\}$ each cycle in the decomposition is the union of the $A^\ell$, $B^\ell(J)$, and $C^\ell(I)$ sets, i.e.

\[ A^\ell \cup
\left(\bigcup_{J=0}^{N-1}B^\ell(J) \right)
\cup
\left(\bigcup_{I=0}^{M-1}C^\ell(I) \right).\]

Since the $A$ block contains 24 edges of each color, each of the $N$ $B$ blocks contain 24 edges of each color, and each of the $M$ $C$ blocks contain $8N+8$ edges of each color, each cycle contains $24 + 24N + 8M + 8MN$ edges.  Since $m=6+2M$ and $n=6+6N$, the total number of edges in $C_m \sq C_n$ is $2mn= 2(6+2M)(6+6N) = 72 + 72N + 24M + 24MN$.  Thus each cycle contains one third of the total edges, as desired.

\begin{figure}
\begin{center}
\begin{tikzpicture}
    \draw[blue, ultra thick, densely dotted] (-.5, 0) -- (2, 0) -- (2, -1) -- (5, -1) -- (5, -2) -- (5.5, -2); \draw[blue, ultra thick, densely dotted] (-.5, -2) -- (0, -2) -- (0, -3) -- (5, -3) -- (5, -4) -- (5.5, -4); \draw[blue, ultra thick, densely dotted] (5, .5) -- (5, 0) -- (5.5, 0); \draw[blue, ultra thick, densely dotted] (5.5, -5) -- (5, -5) -- (5, -5.5); \draw[blue, ultra thick, densely dotted] (-.5, -4) -- (2, -4) -- (2, -5) -- (-.5, -5);
    \draw[orange, ultra thick] (0, .5) -- (0, -1) -- (1, -1) -- (1, .5); \draw[orange, ultra thick] (3, .5) -- (3, 0) -- (4, 0) -- (4, -2) -- (5, -2) -- (5, -3) -- (5.5, -3); \draw[orange, ultra thick] (-.5, -3) -- (0, -3) -- (0, -5.5); \draw[orange, ultra thick] (1, -5.5)  -- (1, -2) -- (3, -2) -- (3, -4) -- (4, -4) -- (4, -5) -- (3, -5) -- (3, -5.5);
    \draw[red, ultra thick, dashed] (-.5, -1) -- (0, -1) -- (0, -2) -- (1, -2) -- (1, -1) -- (2, -1) -- (2, -4) -- (3, -4) -- (3, -5) -- (2, -5) -- (2, -5.5); \draw[dashed, red, ultra thick] (2, .5) -- (2, 0) -- (3, 0) -- (3, -2) -- (4, -2) -- (4, -4) -- (5, -4) -- (5, -5) -- (4, -5) -- (4, -5.5); \draw[dashed, red, ultra thick] (4, .5) -- (4, 0) -- (5, 0) -- (5, -1) -- (5.5, -1);
    \foreach \x in {0, 1, ..., 5}{\foreach \y in {0, -1, ..., -5}{\filldraw[draw=black, fill=white] (\x, \y) circle (3pt);}}
    \foreach \x in {0,1,2,3,4,5}{\node at (\x,1){$\x$};}\foreach \y\i in {0/0,-1/1,-2/2,-3/3,-4/4,-5/5}{\node at (-1,\y){$\i$};}
\end{tikzpicture}
\hfill
    \begin{tikzpicture}
    \draw[blue, ultra thick, densely dotted] (6, 0) -- (12, 0); \draw[blue, ultra thick, densely dotted] (6, -2) -- (12, -2); \draw[blue, ultra thick, densely dotted] (6, -4) -- (12, -4); \draw[blue, ultra thick, densely dotted] (6, -5) -- (12, -5);
    \draw[red, ultra thick, dashed] (6, -1) -- (6.5, -1) -- (6.5, .5); \draw[dashed, red, ultra thick] (8.5, .5) -- (8.5, -1) -- (7.5, -1) -- (7.5, -3) -- (6.5, -3) -- (6.5, -5.5); \draw[red, ultra thick, dashed] (10.5, .5) -- (10.5, -1) -- (9.5, -1) -- (9.5, -3) -- (8.5, -3) -- (8.5, -5.5); \draw[red, ultra thick, dashed] (12, -1) -- (11.5, -1) -- (11.5, -3) -- (10.5, -3) -- (10.5, -5.5);
    \draw[orange, ultra thick] (7.5, .5) -- (7.5, -1) -- (6.5, -1) -- (6.5, -3) -- (6, -3); \draw[orange, ultra thick] (9.5, .5) -- (9.5, -1) -- (8.5, -1) -- (8.5, -3) -- (7.5, -3) -- (7.5, -5.5); \draw[orange, ultra thick] (11.5, .5) -- (11.5, -1) -- (10.5, -1) -- (10.5, -3) -- (9.5, -3) -- (9.5, -5.5); \draw[orange, ultra thick] (12, -3) -- (11.5, -3) -- (11.5, -5.5);
    \foreach \x in {6.5, 7.5, 8.5, 9.5, 10.5, 11.5}{\foreach \y in {0, -1, ..., -5}{\filldraw[draw=black, fill=white] (\x, \y) circle (3pt);}}
    \foreach \x\j in {6.5/0,7.5/1,8.5/2,9.5/3,10.5/4,11.5/5}{\node at (\x,1){$\j$};}\foreach \y\i in {0/0,-1/1,-2/2,-3/3,-4/4,-5/5}{\node at (5.5,\y){$\i$};}
\end{tikzpicture}
\end{center}
\caption{Left: The $A$ block for Case 6, with $A^r$, $A^y$, and $A^b$ in dashed red, solid yellow, and dotted blue.
Right: The $B$ block for Case 6, with $B^r$, $B^y$, and $B^b$ in dashed red, solid yellow, and dotted blue.}
\label{BCase6}
\end{figure}
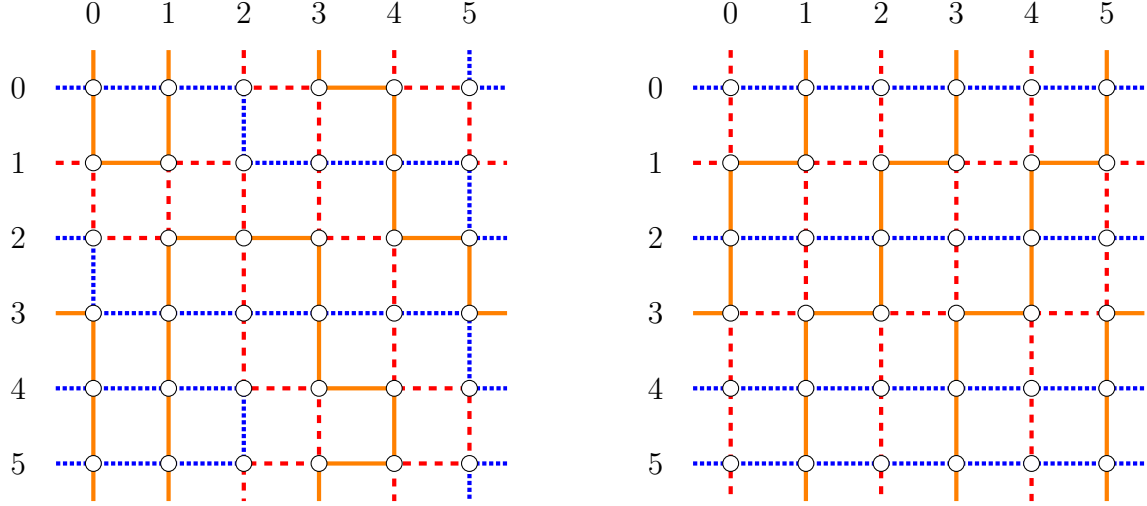
\begin{figure}
\begin{center}
\begin{tikzpicture}[scale=.7, transform shape]
    \draw[blue, ultra thick, densely dotted] (-.5, 0) -- (5, 0) -- (5, .5); \draw[blue, ultra thick, densely dotted] (-.5, -1) -- (5, -1) -- (5, -1.5); 
    \draw[orange, ultra thick] (0, .5) -- (0, -1.5); \draw[orange, ultra thick] (1, .5) -- (1, -1.5); \draw[orange, ultra thick] (3, .5) -- (3, -1.5); \draw[orange, ultra thick] (7, .5) -- (7, 0) -- (6, 0) -- (6, -1) -- (7, -1) -- (7, -1.5); \draw[orange, ultra thick] (9, .5) -- (9, 0) -- (8, 0) -- (8, -1) -- (9, -1) -- (9, -1.5);
    \draw[red, ultra thick, dashed] (2, .5) -- (2, -1.5);\draw[red, ultra thick, dashed] (4, .5) -- (4, -1.5); \draw[red, ultra thick, dashed] (6, .5) -- (6, 0) -- (5, 0) -- (5, -1) -- (6, -1) -- (6, -1.5); \draw[red, ultra thick, dashed] (8, .5) -- (8, 0) -- (7, 0) -- (7, -1) -- (8, -1) -- (8, -1.5); \draw[red, dashed, ultra thick] (9.5, 0) -- (9,0) -- (9,-1) -- (9.5, -1);
    \node at (10, -.5){$\cdots$};
    \draw[blue, densely dotted, ultra thick] (14.5, 0) -- (13,0) -- (13,-1) -- (14.5, -1);
    \draw[orange, ultra thick] (12, .5) -- (12, 0) -- (11, 0) -- (11, -1) -- (12, -1) -- (12, -1.5);\draw[orange, ultra thick] (14,.5) -- (14, -1.5);
    \draw[red, ultra thick, dashed] (10.5, 0) -- (11,0) -- (11,.5); \draw[red, ultra thick, dashed] (10.5, -1) -- (11,-1) -- (11,-1.5);\draw[red, ultra thick, dashed] (13, .5) -- (13, 0) -- (12, 0) -- (12, -1) -- (13, -1) -- (13, -1.5);
    \node at (15,-.5){$\cdots$};
    \draw[blue, densely dotted, ultra thick] (15.5, 0) -- (21.5, 0);\draw[blue, densely dotted, ultra thick] (15.5, -1) -- (21.5, -1);
    \draw[orange, ultra thick] (17, .5) -- (17, -1.5);\draw[orange, ultra thick] (19, .5) -- (19, -1.5);\draw[orange, ultra thick] (21, .5) -- (21, -1.5);
    \draw[red, dashed, ultra thick] (16, .5) -- (16, -1.5);\draw[red, dashed, ultra thick] (18, .5) -- (18, -1.5);\draw[red, dashed, ultra thick] (20, .5) -- (20, -1.5);
    \foreach\x in {0, 1, 2, 3, 4, 5, 6, 7, 8, 9, 11, 12, 13, 14, 16, 17, 18, 19, 20, 21}{\foreach \y in {0, -1}{\filldraw[fill=white, draw=black] (\x, \y) circle (3pt);}}
    \foreach \x\j in {0,1,2,3,4,5,6,7,8,9,17/$n-5$, 19/$n-3$, 21/$n-1$}{\node at (\x,1){\j};}\foreach \y\i in {0/0,-1/1}{\node at (-1,\y){$\i$};}
    %\node at (13,1.3){$\underbrace{n-(4N-2)}$}; 
    \draw [thick, decorate, decoration = {calligraphic brace, raise=.5pt, amplitude=5pt}] (14.25, 1.1) --  (11.75, 1.1) node[pos=.5, above]{$n-(4N-2)$}; 
    
\end{tikzpicture}
\caption{The $C$ block for Case 6, with $C^r$, $C^y$, and $C^b$ in dashed red, solid yellow, and dotted blue.}
\label{CCase6}
\end{center}
\end{figure}
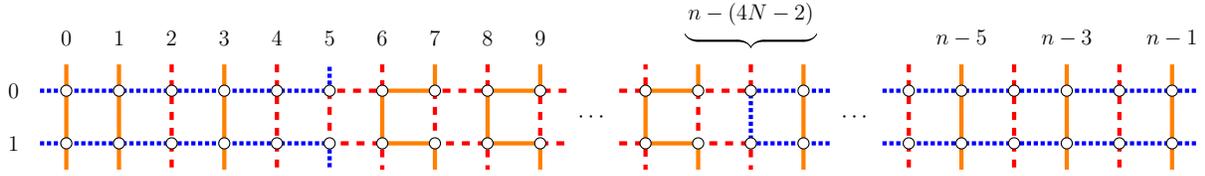

\begin{figure}
\begin{center}
\begin{tikzpicture}[scale=.5, transform shape]
%A
    \draw[blue, ultra thick, densely dotted] (-.5, 0) -- (2, 0) -- (2, -1) -- (5, -1) -- (5, -2) -- (5.5, -2); \draw[blue, ultra thick, densely dotted] (-.5, -2) -- (0, -2) -- (0, -3) -- (5, -3) -- (5, -4) -- (5.5, -4); \draw[blue, ultra thick, densely dotted] (5, .5) -- (5, 0) -- (5.5, 0); \draw[blue, ultra thick, densely dotted] (5.5, -5) -- (5, -5) -- (5, -5.5); \draw[blue, ultra thick, densely dotted] (-.5, -4) -- (2, -4) -- (2, -5) -- (-.5, -5);
    \draw[orange, ultra thick] (0, .5) -- (0, -1) -- (1, -1) -- (1, .5); \draw[orange, ultra thick] (3, .5) -- (3, 0) -- (4, 0) -- (4, -2) -- (5, -2) -- (5, -3) -- (5.5, -3); \draw[orange, ultra thick] (-.5, -3) -- (0, -3) -- (0, -5.5); \draw[orange, ultra thick] (1, -5.5)  -- (1, -2) -- (3, -2) -- (3, -4) -- (4, -4) -- (4, -5) -- (3, -5) -- (3, -5.5);
    \draw[red, ultra thick, dashed] (-.5, -1) -- (0, -1) -- (0, -2) -- (1, -2) -- (1, -1) -- (2, -1) -- (2, -4) -- (3, -4) -- (3, -5) -- (2, -5) -- (2, -5.5); \draw[dashed, red, ultra thick] (2, .5) -- (2, 0) -- (3, 0) -- (3, -2) -- (4, -2) -- (4, -4) -- (5, -4) -- (5, -5) -- (4, -5) -- (4, -5.5); \draw[dashed, red, ultra thick] (4, .5) -- (4, 0) -- (5, 0) -- (5, -1) -- (5.5, -1);
    \foreach \x in {0, 1, ..., 5}{\foreach \y in {0, -1, ..., -5}{\filldraw[draw=black, fill=white] (\x, \y) circle (3pt);}}
%B(0)
    \draw[blue, ultra thick, densely dotted] (6, 0) -- (12, 0); \draw[blue, ultra thick, densely dotted] (6, -2) -- (12, -2); \draw[blue, ultra thick, densely dotted] (6, -4) -- (12, -4); \draw[blue, ultra thick, densely dotted] (6, -5) -- (12, -5);
    \draw[red, ultra thick, dashed] (6, -1) -- (6.5, -1) -- (6.5, .5); \draw[dashed, red, ultra thick] (8.5, .5) -- (8.5, -1) -- (7.5, -1) -- (7.5, -3) -- (6.5, -3) -- (6.5, -5.5); \draw[red, ultra thick, dashed] (10.5, .5) -- (10.5, -1) -- (9.5, -1) -- (9.5, -3) -- (8.5, -3) -- (8.5, -5.5); \draw[red, ultra thick, dashed] (12, -1) -- (11.5, -1) -- (11.5, -3) -- (10.5, -3) -- (10.5, -5.5);
    \draw[orange, ultra thick] (7.5, .5) -- (7.5, -1) -- (6.5, -1) -- (6.5, -3) -- (6, -3); \draw[orange, ultra thick] (9.5, .5) -- (9.5, -1) -- (8.5, -1) -- (8.5, -3) -- (7.5, -3) -- (7.5, -5.5); \draw[orange, ultra thick] (11.5, .5) -- (11.5, -1) -- (10.5, -1) -- (10.5, -3) -- (9.5, -3) -- (9.5, -5.5); \draw[orange, ultra thick] (12, -3) -- (11.5, -3) -- (11.5, -5.5);
    \foreach \x in {6.5, 7.5, 8.5, 9.5, 10.5, 11.5}{\foreach \y in {0, -1, ..., -5}{\filldraw[draw=black, fill=white] (\x, \y) circle (3pt);}}
%B(1)
    \draw[blue, ultra thick, densely dotted] (12.5, 0) -- (18.5, 0); \draw[blue, ultra thick, densely dotted] (12.5, -2) -- (18.5, -2); \draw[blue, ultra thick, densely dotted] (12.5, -4) -- (18.5, -4); \draw[blue, ultra thick, densely dotted] (12.5, -5) -- (18.5, -5);
    \draw[red, ultra thick, dashed] (12.5, -1) -- (13, -1) -- (13, .5); \draw[dashed, red, ultra thick] (15, .5) -- (15, -1) -- (14, -1) -- (14, -3) -- (13, -3) -- (13, -5.5); \draw[red, ultra thick, dashed] (17, .5) -- (17, -1) -- (16, -1) -- (16, -3) -- (15, -3) -- (15, -5.5); \draw[red, ultra thick, dashed] (18.5, -1) -- (18, -1) -- (18, -3) -- (17, -3) -- (17, -5.5);
    \draw[orange, ultra thick] (14, .5) -- (14, -1) -- (13, -1) -- (13, -3) -- (12.5, -3); \draw[orange, ultra thick] (16, .5) -- (16, -1) -- (15, -1) -- (15, -3) -- (14, -3) -- (14, -5.5); \draw[orange, ultra thick] (18, .5) -- (18, -1) -- (17, -1) -- (17, -3) -- (16, -3) -- (16, -5.5); \draw[orange, ultra thick] (18.5, -3) -- (18, -3) -- (18, -5.5);
    \foreach \x in {13, 14, 15, 16, 17, 18}{\foreach \y in {0, -1, ..., -5}{\filldraw[draw=black, fill=white] (\x, \y) circle (3pt);}}
\node at (21.25, -3){$\cdots$};
%B(N-1)
    \draw[blue, ultra thick, densely dotted] (24, 0) -- (30, 0); \draw[blue, ultra thick, densely dotted] (24, -2) -- (30, -2); \draw[blue, ultra thick, densely dotted] (24, -4) -- (30, -4); \draw[blue, ultra thick, densely dotted] (24, -5) -- (30, -5);
    \draw[red, ultra thick, dashed] (24, -1) -- (24.5, -1) -- (24.5, .5); \draw[dashed, red, ultra thick] (26.5, .5) -- (26.5, -1) -- (25.5, -1) -- (25.5, -3) -- (24.5, -3) -- (24.5, -5.5); \draw[red, ultra thick, dashed] (28.5, .5) -- (28.5, -1) -- (27.5, -1) -- (27.5, -3) -- (26.5, -3) -- (26.5, -5.5); \draw[red, ultra thick, dashed] (30, -1) -- (29.5, -1) -- (29.5, -3) -- (28.5, -3) -- (28.5, -5.5);
    \draw[orange, ultra thick] (25.5, .5) -- (25.5, -1) -- (24.5, -1) -- (24.5, -3) -- (24, -3); \draw[orange, ultra thick] (27.5, .5) -- (27.5, -1) -- (26.5, -1) -- (26.5, -3) -- (25.5, -3) -- (25.5, -5.5); \draw[orange, ultra thick] (29.5, .5) -- (29.5, -1) -- (28.5, -1) -- (28.5, -3) -- (27.5, -3) -- (27.5, -5.5); \draw[orange, ultra thick] (30, -3) -- (29.5, -3) -- (29.5, -5.5);
    \foreach \x in {24.5, 25.5, 26.5, 27.5, 28.5, 29.5}{\foreach \y in {0, -1, ..., -5}{\filldraw[draw=black, fill=white] (\x, \y) circle (3pt);}}
%C(0)
    \draw[blue, ultra thick, densely dotted] (-.5, -6.5) -- (5, -6.5) -- (5, -6); \draw[blue, ultra thick, densely dotted] (-.5, -7.5) -- (5, -7.5) -- (5, -8); 
    \draw[orange, ultra thick] (0, -6) -- (0, -8); \draw[orange, ultra thick] (1, -6) -- (1, -8); \draw[orange, ultra thick] (3, -6) -- (3, -8); 
    \draw[orange, ultra thick] (7.5, -6) -- (7.5, -6.5) -- (6.5, -6.5) -- (6.5, -7.5) -- (7.5, -7.5) -- (7.5, -8); \draw[orange, ultra thick] (9.5, -6) -- (9.5, -6.5) -- (8.5, -6.5) -- (8.5, -7.5) -- (9.5, -7.5) -- (9.5, -8); \draw[orange, ultra thick] (11.5, -6) -- (11.5, -6.5) -- (10.5, -6.5) -- (10.5, -7.5) -- (11.5, -7.5) -- (11.5, -8); 
    \draw[red, ultra thick, dashed] (2, -6) -- (2, -8);\draw[red, ultra thick, dashed] (4, -6) -- (4, -8); \draw[red, ultra thick, dashed] (6.5, -6) -- (6.5, -6.5) -- (5, -6.5) -- (5, -7.5) -- (6.5, -7.5) -- (6.5, -8); \draw[red, dashed, ultra thick] (12, -6.5) -- (11.5, -6.5) -- (11.5, -7.5) -- (12, -7.5);
    \draw[red, ultra thick, dashed] (8.5, -6) -- (8.5, -6.5) -- (7.5, -6.5) -- (7.5, -7.5) -- (8.5, -7.5) -- (8.5, -8);\draw[red, ultra thick, dashed] (10.5, -6) -- (10.5, -6.5) -- (9.5, -6.5) -- (9.5, -7.5) -- (10.5, -7.5) -- (10.5, -8);
    \node at (15, -7){$\cdots$};
    \draw[red, dashed, ultra thick] (19.5, -6) -- (19.5, -6.5) -- (19, -6.5); \draw[red, dashed, ultra thick] (19.5, -8) -- (19.5, -7.5) -- (19, -7.5);\draw[red, dashed, ultra thick] (21.5, -6) -- (21.5, -6.5) -- (20.5, -6.5) -- (20.5, -7.5) -- (21.5, -7.5) -- (21.5, -8); 
    \draw[orange, ultra thick] (20.5, -6) -- (20.5, -6.5) -- (19.5, -6.5) -- (19.5, -7.5) -- (20.5, -7.5) -- (20.5, -8); \draw[orange, ultra thick] (22.5, -6) -- (22.5, -8);
    \draw[blue, ultra thick, densely dotted] (23, -6.5) -- (21.5, -6.5) -- (21.5, -7.5) -- (23, -7.5);
    \node at (23.5, -7){$\cdots$};
    \draw[blue, densely dotted, ultra thick] (30, -6.5) -- (24, -6.5);\draw[blue, densely dotted, ultra thick] (30, -7.5) -- (24, -7.5);\draw[red, ultra thick, dashed] (24.5, -6) -- (24.5, -8);\draw[red, ultra thick, dashed] (26.5, -6) -- (26.5, -8);\draw[red, ultra thick, dashed] (28.5, -6) -- (28.5, -8);
    \draw[orange, ultra thick] (25.5, -6) -- (25.5, -8);\draw[orange, ultra thick] (27.5, -6) -- (27.5, -8);\draw[orange, ultra thick] (29.5, -6) -- (29.5, -8);
    \foreach \x in {0, 1, 2, 3, 4, 5, 6.5, 7.5, 8.5, 9.5, 10.5, 11.5, 19.5, 20.5, 21.5, 22.5, 24.5, 25.5, 26.5, 27.5, 28.5, 29.5}{\foreach \y in {-6.5, -7.5}{\filldraw[draw=black, fill=white] (\x, \y) circle (3pt);}}
%C(1)
    \draw[blue, ultra thick, densely dotted] (-.5, -9) -- (5, -9) -- (5, -8.5); \draw[blue, ultra thick, densely dotted] (-.5, -10) -- (5, -10) -- (5, -10.5); 
    \draw[orange, ultra thick] (0, -8.5) -- (0, -10.5); \draw[orange, ultra thick] (1, -8.5) -- (1, -10.5); \draw[orange, ultra thick] (3, -8.5) -- (3, -10.5); 
    \draw[orange, ultra thick] (7.5, -8.5) -- (7.5, -9) -- (6.5, -9) -- (6.5, -10) -- (7.5, -10) -- (7.5, -10.5); \draw[orange, ultra thick] (9.5, -8.5) -- (9.5, -9) -- (8.5, -9) -- (8.5, -10) -- (9.5, -10) -- (9.5, -10.5); \draw[orange, ultra thick] (11.5, -8.5) -- (11.5, -9) -- (10.5, -9) -- (10.5, -10) -- (11.5, -10) -- (11.5, -10.5); 
    \draw[red, ultra thick, dashed] (2, -8.5) -- (2, -10.5);\draw[red, ultra thick, dashed] (4, -8.5) -- (4, -10.5); \draw[red, ultra thick, dashed] (6.5, -8.5) -- (6.5, -9) -- (5, -9) -- (5, -10) -- (6.5, -10) -- (6.5, -10.5); \draw[red, dashed, ultra thick] (12, -9) -- (11.5, -9) -- (11.5, -10) -- (12, -10);
    \draw[red, ultra thick, dashed] (8.5, -8.5) -- (8.5, -9) -- (7.5, -9) -- (7.5, -10) -- (8.5, -10) -- (8.5, -10.5);\draw[red, ultra thick, dashed] (10.5, -8.5) -- (10.5, -9) -- (9.5, -9) -- (9.5, -10) -- (10.5, -10) -- (10.5, -10.5);
    \node at (15, -9.5){$\cdots$};
    \draw[red, dashed, ultra thick] (19.5, -8.5) -- (19.5, -9) -- (19, -9); \draw[red, dashed, ultra thick] (19.5, -10.5) -- (19.5, -10) -- (19, -10);\draw[red, dashed, ultra thick] (21.5, -8.5) -- (21.5, -9) -- (20.5, -9) -- (20.5, -10) -- (21.5, -10) -- (21.5, -10.5); 
    \draw[orange, ultra thick] (20.5, -8.5) -- (20.5, -9) -- (19.5, -9) -- (19.5, -10) -- (20.5, -10) -- (20.5, -10.5); \draw[orange, ultra thick] (22.5, -8.5) -- (22.5, -10.5);
    \draw[blue, ultra thick, densely dotted] (23, -9) -- (21.5, -9) -- (21.5, -10) -- (23, -10);
    \node at (23.5, -9.5){$\cdots$};
    \draw[blue, densely dotted, ultra thick] (30, -9) -- (24, -9);\draw[blue, densely dotted, ultra thick] (30, -10) -- (24, -10);\draw[red, ultra thick, dashed] (24.5, -8.5) -- (24.5, -10.5);\draw[red, ultra thick, dashed] (26.5, -8.5) -- (26.5, -10.5);\draw[red, ultra thick, dashed] (28.5, -8.5) -- (28.5, -10.5);
    \draw[orange, ultra thick] (25.5, -8.5) -- (25.5, -10.5);\draw[orange, ultra thick] (27.5, -8.5) -- (27.5, -10.5);\draw[orange, ultra thick] (29.5, -8.5) -- (29.5, -10.5);
    \foreach \x in {0, 1, 2, 3, 4, 5, 6.5, 7.5, 8.5, 9.5, 10.5, 11.5, 19.5, 20.5, 21.5, 22.5, 24.5, 25.5, 26.5, 27.5, 28.5, 29.5}{\foreach \y in {-9, -10}{\filldraw[draw=black, fill=white] (\x, \y) circle (3pt);}}
\node at (15, -11.25){$\vdots$};

%C(M-1)
    \draw[blue, ultra thick, densely dotted] (-.5, -12.5) -- (5, -12.5) -- (5, -12); \draw[blue, ultra thick, densely dotted] (-.5, -13.5) -- (5, -13.5) -- (5, -14); 
    \draw[orange, ultra thick] (0, -12) -- (0, -14); \draw[orange, ultra thick] (1, -12) -- (1, -14); \draw[orange, ultra thick] (3, -12) -- (3, -14); 
    \draw[orange, ultra thick] (7.5, -12) -- (7.5, -12.5) -- (6.5, -12.5) -- (6.5, -13.5) -- (7.5, -13.5) -- (7.5, -14); \draw[orange, ultra thick] (9.5, -12) -- (9.5, -12.5) -- (8.5, -12.5) -- (8.5, -13.5) -- (9.5, -13.5) -- (9.5, -14); \draw[orange, ultra thick] (11.5, -12) -- (11.5, -12.5) -- (10.5, -12.5) -- (10.5, -13.5) -- (11.5, -13.5) -- (11.5, -14); 
    \draw[red, ultra thick, dashed] (2, -12) -- (2, -14);\draw[red, ultra thick, dashed] (4, -12) -- (4, -14); \draw[red, ultra thick, dashed] (6.5, -12) -- (6.5, -12.5) -- (5, -12.5) -- (5, -13.5) -- (6.5, -13.5) -- (6.5, -14); \draw[red, dashed, ultra thick] (12, -12.5) -- (11.5, -12.5) -- (11.5, -13.5) -- (12, -13.5);
    \draw[red, ultra thick, dashed] (8.5, -12) -- (8.5, -12.5) -- (7.5, -12.5) -- (7.5, -13.5) -- (8.5, -13.5) -- (8.5, -14);\draw[red, ultra thick, dashed] (10.5, -12) -- (10.5, -12.5) -- (9.5, -12.5) -- (9.5, -13.5) -- (10.5, -13.5) -- (10.5, -14);
    \node at (15, -13){$\cdots$};
    \draw[red, dashed, ultra thick] (19.5, -12) -- (19.5, -12.5) -- (19, -12.5); \draw[red, dashed, ultra thick] (19.5, -14) -- (19.5, -13.5) -- (19, -13.5);\draw[red, dashed, ultra thick] (21.5, -12) -- (21.5, -12.5) -- (20.5, -12.5) -- (20.5, -13.5) -- (21.5, -13.5) -- (21.5, -14); 
    \draw[orange, ultra thick] (20.5, -12) -- (20.5, -12.5) -- (19.5, -12.5) -- (19.5, -13.5) -- (20.5, -13.5) -- (20.5, -14); \draw[orange, ultra thick] (22.5, -12) -- (22.5, -14);
    \draw[blue, ultra thick, densely dotted] (23, -12.5) -- (21.5, -12.5) -- (21.5, -13.5) -- (23, -13.5);
    \node at (23.5, -13){$\cdots$};
    \draw[blue, densely dotted, ultra thick] (30, -12.5) -- (24, -12.5);\draw[blue, densely dotted, ultra thick] (30, -13.5) -- (24, -13.5);\draw[red, ultra thick, dashed] (24.5, -12) -- (24.5, -14);\draw[red, ultra thick, dashed] (26.5, -12) -- (26.5, -14);\draw[red, ultra thick, dashed] (28.5, -12) -- (28.5, -14);
    \draw[orange, ultra thick] (25.5, -12) -- (25.5, -14);\draw[orange, ultra thick] (27.5, -12) -- (27.5, -14);\draw[orange, ultra thick] (29.5, -12) -- (29.5, -14);
    \foreach \x in {0, 1, 2, 3, 4, 5, 6.5, 7.5, 8.5, 9.5, 10.5, 11.5, 19.5, 20.5, 21.5, 22.5, 24.5, 25.5, 26.5, 27.5, 28.5, 29.5}{\foreach \y in {-12.5, -13.5}{\filldraw[draw=black, fill=white] (\x, \y) circle (3pt);}}

    \foreach \x\i in {0/0, 1/1, 2/2, 3/3, 4/4, 5/5, 6.5/6, 7.5/7, 8.5/8, 9.5/9, 10.5/10, 11.5/11, 13/12, 14/13, 15/14, 16/15, 17/16, 18/17, 24.5/, 25.5/$n-5$, 26.5/, 27.5/$n-3$, 28.5/, 29.5/$n-1$}{\node at (\x, 1.5){\i};}
    \foreach \y\j in {0/0, -1/1, -2/2, -3/3, -4/4, -5/5, -6.5/6, -7.5/7, -9/8, -10/9, -12.5/$m-2$, -13.5/$m-1$}{\node at (-1.5, \y){\j};}
\end{tikzpicture} 
\begin{tikzpicture}[scale=.5, transform shape]
 \huge{   \node at (-1.5, 3){};\node at (2.5, 0){$A$};\node at (9, 0){$B(0)$};\node at (15.5, 0){$B(1)$};\node at (21.25, 0){$\cdots$};\node at (28, 0){$B(N-1)$};\node at (15.5, -2){$C(0)$};\node at (15.5, -4){$C(1)$};\node at (15.5, -5){$\vdots$};\node at (15.5, -7){$C(M-1)$};}
\end{tikzpicture}
\end{center}
\caption{Above: The decomposition of $C_m\sq C_n$ for Case 6 in dashed red, solid yellow, and dotted blue. Below: The arrangement of $A$, $B$, and $C$ blocks in $C_m\sq C_n$ for Case 6.}\label{ABCCase6}
\end{figure}
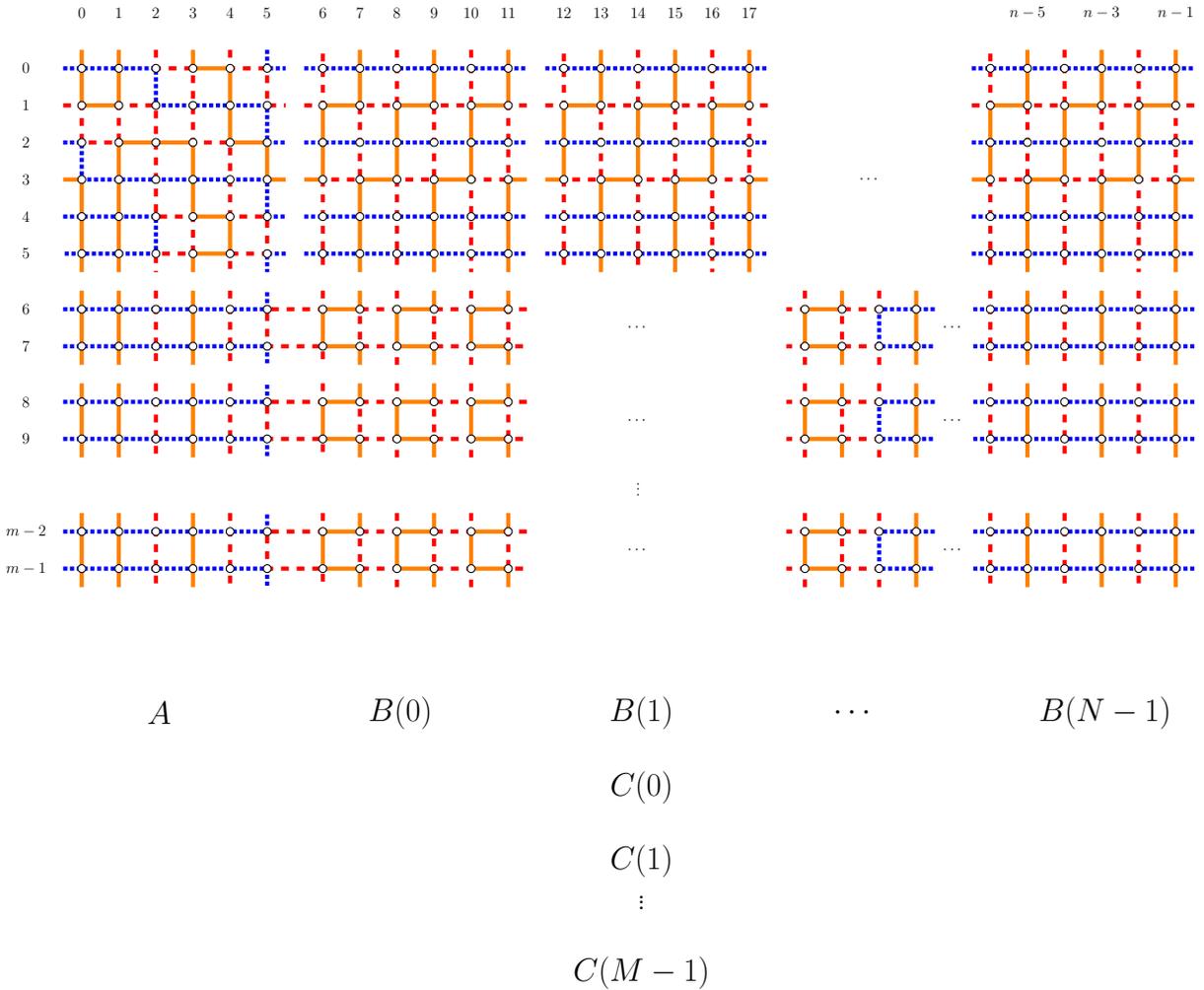

\paragraph{Case 7 ($m \geq 6$ is even and $n$ is odd):} Assume $m=6+2M$ and $n=3+6N$ for nonnegative integers $M$ and $N$. 
The decomposition of $C_m \sq C_n$ is illustrated in Figure~\ref{Case7} and contains one $6 \times 3$ $A$ block, $N$ $6 \times 6$ $B$ blocks, $M$ $2 \times 3$ $C$ blocks, and $MN$ $2 \times 6$ $D$ blocks, shown in Figure~\ref{ABCDCase7} and defined below.  The $A$ block is the minimal block for this case.

For the $6 \times 3$ $A$ block, define three edge sets $A^r$, $A^y$, and $A^b$, each of which forms a cycle with 12 edges, shown as dashed red, solid yellow, and dotted blue, respectively, in Figure~\ref{ABCDCase7} (above left).

$A^r = \{(1,-1)(1,0), (1,0)(0,0), (0,0)(0,1), (0,1)(-1,1), (6,1)(5,1), (5,1)(5,2),\\ (5,2)(4,2), (4,2)(4,3), (4,-1)(4,0), (4,0)(3,0), (3,0)(3,1), (3,1)(2,1), (2,1)(2,2), \\(2,2)(1,2), (1,2)(1,3)\}.$

$A^y = \{(2,-1)(2,0), (2,0)(1,0), (1,0)(1,1), (1,1)(0,1), (0,1)(0,2), (0,2)(-1,2),\\ (6,2)(5,2), (5,2)(5,3), (5,-1)(5,0), (5,0)(4,0), (4,0)(4,1), (4,1)(3,1), (3,1)(3,2), \\(3,2)(2,2), (2,2)(2,3)\}.$

$A^b = \{(0,-1)(0,0), (0,0)(-1,0), (6,0)(5,0), (5,0)(5,1), (5,1)(4,1), (4,1)(4,2),\\ (4,2)(3,2), (3,2)(3,3), (3,-1)(3,0), (3,0)(2,0), (2,0)(2,1), (2,1)(1,1), (1,1)(1,2),\\ (1,2)(0,2), (0,2)(0,3)\}.$

For the $6 \times 6$ block $B$ define edge sets $B^r$, $B^y$, and $B^b$, each containing 24 edges, shown as dashed red, solid yellow, and dotted blue, respectively in Figure~\ref{ABCDCase7} (above right).

$B^r = \{(1,-1)(1,0), (1,0)(0,0), (0,0)(0,1), (0,1)(-1,1), (6,1)(5,1), (5,1)(5,2),\\ (5,2)(4,2), (4,2)(4,3), (4,3)(3,3), (3,3)(3,4), (3,4)(2,4), (2,4)(2,5), (2,5)(1,5),\\ (1,5)(1,6), (4,-1)(4,0), (4,0)(3,0), (3,0)(3,1), (3,1)(2,1), (2,1)(2,2), (2,2)(1,2), \\(1,2)(1,3), (1,3)(0,3), (0,3)(0,4), (0,4)(-1,4), (6,4)(5,4), (5,4)(5,5), (5,5)(4,5),\\ (4,5)(4,6)\}$

$B^y = \{(2,-1)(2,0), (2,0)(1,0), (1,0)(1,1), (1,1)(0,1), (0,1)(0,2), (0,2)(-1,2),\\ (6,2)(5,2), (5,2)(5,3), (5,3)(4,3), (4,3)(4,4), (4,4)(3,4), (3,4)(3,5), (3,5)(2,5),\\ (2,5)(2,6), (5,-1)(5,0), (5,0)(4,0), (4,0)(4,1), (4,1)(3,1), (3,1)(3,2), (3,2)(2,2),\\ (2,2)(2,3), (2,3)(1,3), (1,3)(1,4), (1,4)(0,4), (0,4)(0,5), (0,5)(-1,5), (6,5)(5,5),\\ (5,5)(5,6)\}$

$B^b = \{(0,-1)(0,0), (0,0)(-1,0), (6,0)(5,0), (5,0)(5,1), (5,1)(4,1), (4,1)(4,2),\\ (4,2)(3,2), (3,2)(3,3), (3,3)(2,3), (2,3)(2,4), (2,4)(1,4), (1,4)(1,5), (1,5)(0,5), \\(0,5)(0,6), (3,-1)(3,0), (3,0)(2,0), (2,0)(2,1), (2,1)(1,1), (1,1)(1,2), (1,2)(0,2),\\ (0,2)(0,3), (0,3)(-1,3), (6,3)(5,3), (5,3)(5,4), (5,4)(4,4), (4,4)(4,5), (4,5)(3,5), \\(3,5)(3,6)\}$

For $0\leq J\leq N-1$, let the block $B(J)$ in $C_m \sq C_n$ be the set of edges incident to the vertices $(i,j)$ where $0\leq i\leq 5$ and $3+6J \leq j \leq 8+6J$.
To color the edges in $B(J)$, define three edge sets $B^r(J)$, $B^y(J)$, and $B^b(J)$, where for $\ell \in \{r,y,b\}$, edge $(i,j)(i',j')$ is in $B^\ell$ if and only if $(i,j+ 3 +6J)(i',j'+3 +6J)$ is in $B^\ell(J)$. For example, $(0,13)(0,14)$ is yellow (in $B^y(1)$) because $(0,4)(0,5)$ is yellow (in $B^y$).

For the $2 \times 3$ $C$ block, define three edge sets $C^r$, $C^y$, and $C^b$, each with four edges, shownT as dashed red, solid yellow, and dotted blue, respectively in Figure~\ref{ABCDCase7} (below left).

$C^r =\{(-1, 1)(0, 1), (0, 1)(0, 2), (0, 2)(1, 2), (1, 2)(1, 1), (1, 1)(2, 1)\}$

$C^y=\{(-1, 2)(0, 2), (0, 2)(0, 3), (0, -1)(0, 0), (0, 0)(1, 0), (1, 0)(1, -1), (1, 3)(1, 2),\\ (1, 2)(2, 2)\}$

$C^b = \{(-1, 0)(0, 0), (0,0)(0,1), (0, 1)(1, 1), (1, 1)(1, 0), (1, 0)(2,0)\}$

For $0 \le I \le M-1$, let the block $C(I)$ in $C_m \sq c_n$ be the set of edges incident to the vertices $(i,j)$ where $6+2I \leq i\leq 7+2I$ and $0\leq j\leq 2$.
To color the edges in $C(I)$, define three edge sets $C^r(I)$, $C^y(I)$, and $C^b(I)$, where for $\ell \in \{r,y,b\}$, edge $(i,j)(i',j')$ is in $C^\ell$ if and only if $(i + 6 + 2I,j)(i'+6 +2I,j')$ is in $C^\ell(I)$. For example, $(8,0)(9,0)$ is red (in $C^r(1)$) because $(0,0)(1,0)$ is red (in $C^r$).

Finally, for the $2 \times 6$ $D$ block define three edge sets $D^r$, $D^y$, and $D^b$, each containing eight edges, shown as dashed red, solid yellow, and dotted blue, respectively in Figure~\ref{ABCDCase7} (below right).

$D^r =\{(-1, 1)(0, 1), (0, 1)(0, 2), (0, 2)(1, 2), (1, 2)(1, 1), (1, 1)(2, 1), (-1, 4)(0, 4), \\(0, 4)(0, 5), (0, 5)(1, 5), (1, 5)(1, 4), (1, 4)(2, 4)\}$

$D^y=\{(-1, 2)(0, 2), (0, 2)(0, 3), (0, 3)(1, 3), (1, 3)(1, 2), (1, 2)(2, 2), (0, -1)(0, 0), \\(0, 0)(1, 0), (1, 0)(1, -1), (-1,5)(0,5), (0,5)(0,6), (2,5)(1,5), (1,5)(1,6)\}$

$D^b = \{(-1, 0)(0, 0), (0,0)(0,1), (0, 1)(1, 1), (1, 1)(1, 0), (1, 0)(2,0), (-1, 3)(0, 3),\\ (0,3)(0,4), (0, 4)(1, 4), (1, 4)(1, 3), (1, 3)(2,3)\}$

For $0 \le I \le M-1$ and $0 \le J \le N-1$, let the block $D(I,J)$ in $C_m \sq C_n$ be the set of edges incident to the vertices $(i,j)$ where $6+2I\leq i\leq 7+2I$ and $3+6J\leq j \leq 8+6J$.
To color the edges in $D(I,J)$, define three edge sets $D^r(I,J)$, $D^y(I,J)$, and $D^b(I,J)$, where for $\ell \in \{r,y,b\}$, edge $(i,j)(i',j')$ is in $D^\ell$ if and only if $(i + 6 + 2I,j + 3 + 6J)(i'+6 +2I,j' + 3 + 6J)$ is in $D^\ell(I,J)$. For example, $(6,12)(6,13)$ is blue (in $D^b(1,0)$) because $(0,3)(0,4)$ is blue (in $D^b$).

For $\ell \in \{r,y,b\}$ each cycle in the decomposition is the union of the $A^\ell$, $B^\ell(J)$, $C^\ell(I)$, and $D^\ell(I,J)$ sets, i.e.

\[ A^\ell \cup
\left(\bigcup_{J=0}^{N-1}B^\ell(J) \right)
\cup
\left(\bigcup_{I=0}^{M-1}C^\ell(I) \right)
\cup
\left(\bigcup_{I=0}^{M-1}\bigcup_{J=0}^{N-1}D^\ell(I,J) \right).\]

Since the $A$ block contains 12 edges of each color, each of the $N$ $B$ blocks contains 24 edges of each color, each of the $M$ $C$ blocks contains four edges of each color, and each of the $MN$ $D$ blocks contains eight edges of each color, each cycle contains $12 + 24N + 4M + 8MN$ edges.  Since $m=6+2M$ and $n=3+6N$, the total number of edges in $C_m \sq C_n$ is $2mn= 2(6+2M)(3+6N) = 36 + 72N + 12M + 24MN$.  Thus each cycle contains one third of the total edges, as desired.

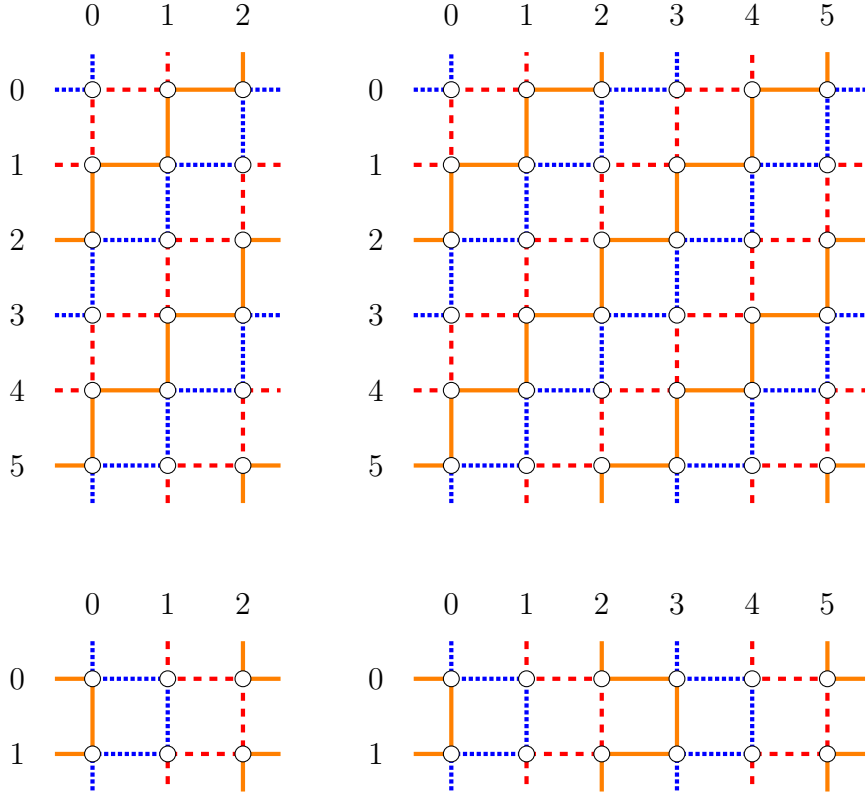
\begin{figure}
\begin{center}
\begin{tikzpicture}
%A Block
    \draw[red, ultra thick, dashed] (-.5, -1) -- (0, -1) -- (0, 0) -- (1, 0) -- (1, .5); \draw[red, ultra thick, dashed] (-.5, -4) -- (0, -4) -- (0, -3) -- (1, -3) -- (1, -2) -- (2, -2) -- (2, -1) -- (2.5, -1); \draw[red, ultra thick, dashed] (1, -5.5) -- (1, -5) -- (2, -5) -- (2, -4) -- (2.5, -4); 
    \draw[blue, ultra thick, densely dotted] (-.5, 0) -- (0, 0) -- (0, .5); \draw[blue, ultra thick, densely dotted] (-.5, -3) -- (0, -3) -- (0, -2) -- (1, -2) -- (1, -1) -- (2, -1) -- (2, 0) -- (2.5,0); \draw[blue, ultra thick, densely dotted] (0, -5.5) -- (0, -5) -- (1, -5) -- (1, -4) -- (2, -4) -- (2, -3) -- (2.5, -3);
    \draw[orange, ultra thick] (-.5, -2) -- (0, -2) -- (0, -1) -- (1, -1) -- (1, 0) -- (2, 0) -- (2, .5); \draw[orange, ultra thick] (-.5, -5) -- (0, -5) -- (0, -4) -- (1, -4) -- (1, -3) -- (2, -3) -- (2, -2) -- (2.5, -2); \draw[orange, ultra thick] (2, -5.5) -- (2, -5) -- (2.5, -5);
    \foreach \x in {0, 1, 2}{\node at (\x, 1){\x};\foreach \y in {0, -1, ..., -5}{\draw[draw = black, fill = white] (\x, \y) circle (3pt);}}
    \foreach \y\i in {0,-1/1,-2/2,-3/3,-4/4,-5/5}{\node at (-1,\y){$\i$};}
    \end{tikzpicture}\hspace{1cm}\begin{tikzpicture}
%B Block
    \draw[red, ultra thick, dashed] (3, -1) -- (3.5, -1) -- (3.5, 0) -- (4.5, 0) -- (4.5, .5); \draw[red, ultra thick, dashed] (3, -4) -- (3.5, -4) -- (3.5, -3) -- (4.5, -3) -- (4.5, -2) -- (5.5, -2) -- (5.5, -1) -- (6.5, -1) -- (6.5, 0) -- (7.5, 0) -- (7.5, .5); \draw[red, ultra thick, dashed] (4.5, -5.5) -- (4.5, -5) -- (5.5, -5) -- (5.5, -4) -- (6.5, -4) -- (6.5, -3) -- (7.5, -3) -- (7.5, -2) -- (8.5, -2) -- (8.5, -1) -- (9, -1); \draw[red, ultra thick, dashed] (7.5, -5.5) -- (7.5, -5) -- (8.5, -5) -- (8.5, -4) -- (9, -4);
    \draw[blue, ultra thick, densely dotted] (3, 0) -- (3.5, 0) -- (3.5, .5); \draw[blue, ultra thick, densely dotted] (3, -3) -- (3.5, -3) -- (3.5, -2) -- (4.5, -2) -- (4.5, -1) -- (5.5, -1) -- (5.5, 0) -- (6.5, 0) -- (6.5, .5); \draw[blue, ultra thick, densely dotted] (3.5, -5.5) -- (3.5, -5) -- (4.5, -5) -- (4.5, -4) -- (5.5, -4) -- (5.5, -3) -- (6.5, -3) -- (6.5, -2) -- (7.5, -2) -- (7.5, -1) -- (8.5, -1) -- (8.5, 0) -- (9, 0);\draw[blue, ultra thick, densely dotted] (6.5, -5.5) -- (6.5, -5) -- (7.5, -5) -- (7.5, -4) -- (8.5, -4) -- (8.5, -3) -- (9, -3);
    \draw[orange, ultra thick] (3, -2) -- (3.5, -2) -- (3.5, -1) -- (4.5, -1) -- (4.5, 0) -- (5.5, 0) -- (5.5, .5); \draw[orange, ultra thick] (3, -5) -- (3.5, -5) -- (3.5, -4) -- (4.5, -4) -- (4.5, -3) -- (5.5, -3) -- (5.5, -2) -- (6.5, -2) -- (6.5, -1) -- (7.5, -1) -- (7.5, 0) -- (8.5, 0) -- (8.5, .5); \draw[orange, ultra thick] (5.5, -5.5) -- (5.5, -5) -- (6.5, -5) -- (6.5, -4) -- (7.5, -4) -- (7.5, -3) -- (8.5, -3) -- (8.5, -2) -- (9, -2); \draw[orange, ultra thick] (8.5, -5.5) -- (8.5, -5) -- (9, -5);
    \foreach \x\j in {3.5/0, 4.5/1, 5.5/2, 6.5/3, 7.5/4, 8.5/5}{\node at (\x,1){\j};\foreach \y in {0, -1, ..., -5}{\draw[draw = black, fill = white] (\x, \y) circle (3pt);}}
    \foreach \y\i in {0,-1/1,-2/2,-3/3,-4/4,-5/5}{\node at (2.5,\y){$\i$};} 
    \end{tikzpicture}
\vskip 1cm    
    \begin{tikzpicture}
%C Block
    \draw[blue, ultra thick, densely dotted] (0, -6) -- (0, -6.5) -- (1, -6.5) -- (1, -7.5) -- (0, -7.5) -- (0, -8); \draw[red, ultra thick, dashed] (1, -6) -- (1, -6.5) -- (2, -6.5) -- (2, -7.5) -- (1, -7.5) -- (1, -8); \draw[orange, ultra thick] (-.5, -6.5) -- (0, -6.5) -- (0, -7.5) -- (-.5, -7.5); \draw[orange, ultra thick] (2, -6) -- (2, -6.5) -- (2.5, -6.5); \draw[orange, ultra thick] (2.5, -7.5) -- (2, -7.5) -- (2, -8); \foreach \x in {0, 1, 2}{\node at (\x,-5.5){\x};\foreach \y in {-6.5, -7.5}{\draw[draw = black, fill = white] (\x, \y) circle (3pt);}}
    \foreach \y\i in {-6.5/0,-7.5/1}{\node at (-1,\y){$\i$};}
    \end{tikzpicture}\hspace{1cm}\begin{tikzpicture}
%D Block
    \draw[blue, ultra thick, densely dotted] (3.5, -6) -- (3.5, -6.5) -- (4.5, -6.5) -- (4.5, -7.5) -- (3.5, -7.5) -- (3.5, -8); \draw[blue, ultra thick, densely dotted] (6.5, -6) -- (6.5, -6.5) -- (7.5, -6.5) -- (7.5, -7.5) -- (6.5, -7.5) -- (6.5, -8); \draw[red, ultra thick, dashed] (4.5, -6) -- (4.5, -6.5) -- (5.5, -6.5) -- (5.5, -7.5) -- (4.5, -7.5) -- (4.5, -8); \draw[red, ultra thick, dashed] (7.5, -6) -- (7.5, -6.5) -- (8.5, -6.5) -- (8.5, -7.5) -- (7.5, -7.5) -- (7.5, -8); \draw[orange, ultra thick] (5.5, -6) -- (5.5, -6.5) -- (6.5, -6.5) -- (6.5, -7.5) -- (5.5, -7.5) -- (5.5, -8); \draw[orange, ultra thick] (3, -6.5) -- (3.5, -6.5) -- (3.5, -7.5) -- (3, -7.5); \draw[orange, ultra thick] (8.5, -6) -- (8.5, -6.5) -- (9, -6.5); \draw[orange, ultra thick] (9, -7.5) -- (8.5, -7.5) -- (8.5, -8); \foreach \x\j in {3.5/0, 4.5/1, 5.5/2, 6.5/3, 7.5/4, 8.5/5}{\node at (\x,-5.5){\j};\foreach \y in {-6.5, -7.5}{\draw[draw = black, fill = white] (\x, \y) circle (3pt);}}
    \foreach \y\i in {-6.5/0,-7.5/1}{\node at (2.5,\y){$\i$};}
\end{tikzpicture}
\end{center}
\caption{For Case 7. Above left: The $A$ block. Above right: The $B$ block. Below left: The $C$ block. Below right: The $D$ block.}\label{ABCDCase7}
\end{figure}

\begin{figure}
\begin{center}
\begin{tikzpicture}[scale=.6, transform shape]
%A
    \draw[red, ultra thick, dashed] (-.5, -1) -- (0, -1) -- (0, 0) -- (1, 0) -- (1, .5); \draw[red, ultra thick, dashed] (-.5, -4) -- (0, -4) -- (0, -3) -- (1, -3) -- (1, -2) -- (2, -2) -- (2, -1) -- (2.5, -1); \draw[red, ultra thick, dashed] (1, -5.5) -- (1, -5) -- (2, -5) -- (2, -4) -- (2.5, -4); 
    \draw[blue, ultra thick, densely dotted] (-.5, 0) -- (0, 0) -- (0, .5); \draw[blue, ultra thick, densely dotted] (-.5, -3) -- (0, -3) -- (0, -2) -- (1, -2) -- (1, -1) -- (2, -1) -- (2, 0) -- (2.5,0); \draw[blue, ultra thick, densely dotted] (0, -5.5) -- (0, -5) -- (1, -5) -- (1, -4) -- (2, -4) -- (2, -3) -- (2.5, -3);
    \draw[orange, ultra thick] (-.5, -2) -- (0, -2) -- (0, -1) -- (1, -1) -- (1, 0) -- (2, 0) -- (2, .5); \draw[orange, ultra thick] (-.5, -5) -- (0, -5) -- (0, -4) -- (1, -4) -- (1, -3) -- (2, -3) -- (2, -2) -- (2.5, -2); \draw[orange, ultra thick] (2, -5.5) -- (2, -5) -- (2.5, -5);
    \foreach \x in {0, 1, 2}{\foreach \y in {0, -1, ..., -5}{\draw[draw = black, fill = white] (\x, \y) circle (3pt);}}
    
%B(0)
    \draw[red, ultra thick, dashed] (3, -1) -- (3.5, -1) -- (3.5, 0) -- (4.5, 0) -- (4.5, .5); \draw[red, ultra thick, dashed] (3, -4) -- (3.5, -4) -- (3.5, -3) -- (4.5, -3) -- (4.5, -2) -- (5.5, -2) -- (5.5, -1) -- (6.5, -1) -- (6.5, 0) -- (7.5, 0) -- (7.5, .5); \draw[red, ultra thick, dashed] (4.5, -5.5) -- (4.5, -5) -- (5.5, -5) -- (5.5, -4) -- (6.5, -4) -- (6.5, -3) -- (7.5, -3) -- (7.5, -2) -- (8.5, -2) -- (8.5, -1) -- (9, -1); \draw[red, ultra thick, dashed] (7.5, -5.5) -- (7.5, -5) -- (8.5, -5) -- (8.5, -4) -- (9, -4);
    \draw[blue, ultra thick, densely dotted] (3, 0) -- (3.5, 0) -- (3.5, .5); \draw[blue, ultra thick, densely dotted] (3, -3) -- (3.5, -3) -- (3.5, -2) -- (4.5, -2) -- (4.5, -1) -- (5.5, -1) -- (5.5, 0) -- (6.5, 0) -- (6.5, .5); \draw[blue, ultra thick, densely dotted] (3.5, -5.5) -- (3.5, -5) -- (4.5, -5) -- (4.5, -4) -- (5.5, -4) -- (5.5, -3) -- (6.5, -3) -- (6.5, -2) -- (7.5, -2) -- (7.5, -1) -- (8.5, -1) -- (8.5, 0) -- (9, 0);\draw[blue, ultra thick, densely dotted] (6.5, -5.5) -- (6.5, -5) -- (7.5, -5) -- (7.5, -4) -- (8.5, -4) -- (8.5, -3) -- (9, -3);
    \draw[orange, ultra thick] (3, -2) -- (3.5, -2) -- (3.5, -1) -- (4.5, -1) -- (4.5, 0) -- (5.5, 0) -- (5.5, .5); \draw[orange, ultra thick] (3, -5) -- (3.5, -5) -- (3.5, -4) -- (4.5, -4) -- (4.5, -3) -- (5.5, -3) -- (5.5, -2) -- (6.5, -2) -- (6.5, -1) -- (7.5, -1) -- (7.5, 0) -- (8.5, 0) -- (8.5, .5); \draw[orange, ultra thick] (5.5, -5.5) -- (5.5, -5) -- (6.5, -5) -- (6.5, -4) -- (7.5, -4) -- (7.5, -3) -- (8.5, -3) -- (8.5, -2) -- (9, -2); \draw[orange, ultra thick] (8.5, -5.5) -- (8.5, -5) -- (9, -5);
    \foreach \x in {3.5, 4.5, 5.5, 6.5, 7.5, 8.5}{\foreach \y in {0, -1, ..., -5}{\draw[draw = black, fill = white] (\x, \y) circle (3pt);}}
    
%B(1)
    \draw[red, ultra thick, dashed] (9.5, -1) -- (10, -1) -- (10, 0) -- (11, 0) -- (11, .5); \draw[red, ultra thick, dashed] (9.5, -4) -- (10, -4) -- (10, -3) -- (11, -3) -- (11, -2) -- (12, -2) -- (12, -1) -- (13, -1) -- (13, 0) -- (14, 0) -- (14, .5); \draw[red, ultra thick, dashed] (11, -5.5) -- (11, -5) -- (12, -5) -- (12, -4) -- (13, -4) -- (13, -3) -- (14, -3) -- (14, -2) -- (15, -2) -- (15, -1) -- (15.5, -1); \draw[red, ultra thick, dashed] (14, -5.5) -- (14, -5) -- (15, -5) -- (15, -4) -- (15.5, -4);
    \draw[blue, ultra thick, densely dotted] (9.5, 0) -- (10, 0) -- (10, .5); \draw[blue, ultra thick, densely dotted] (9.5, -3) -- (10, -3) -- (10, -2) -- (11, -2) -- (11, -1) -- (12, -1) -- (12, 0) -- (13, 0) -- (13, .5); \draw[blue, ultra thick, densely dotted] (10, -5.5) -- (10, -5) -- (11, -5) -- (11, -4) -- (12, -4) -- (12, -3) -- (13, -3) -- (13, -2) -- (14, -2) -- (14, -1) -- (15, -1) -- (15, 0) -- (15.5, 0);\draw[blue, ultra thick, densely dotted] (13, -5.5) -- (13, -5) -- (14, -5) -- (14, -4) -- (15, -4) -- (15, -3) -- (15.5, -3);
    \draw[orange, ultra thick] (9.5, -2) -- (10, -2) -- (10, -1) -- (11, -1) -- (11, 0) -- (12, 0) -- (12, .5); \draw[orange, ultra thick] (9.5, -5) -- (10, -5) -- (10, -4) -- (11, -4) -- (11, -3) -- (12, -3) -- (12, -2) -- (13, -2) -- (13, -1) -- (14, -1) -- (14, 0) -- (15, 0) -- (15, .5); \draw[orange, ultra thick] (12, -5.5) -- (12, -5) -- (13, -5) -- (13, -4) -- (14, -4) -- (14, -3) -- (15, -3) -- (15, -2) -- (15.5, -2); \draw[orange, ultra thick] (15, -5.5) -- (15, -5) -- (15.5, -5);
    \foreach \x in {10, 11, 12, 13, 14, 15}{\foreach \y in {0, -1, ..., -5}{\draw[draw = black, fill = white] (\x, \y) circle (3pt);}}
    
%B(N-1)
    \draw[red, ultra thick, dashed] (17.5, -1) -- (18, -1) -- (18, 0) -- (19, 0) -- (19, .5); \draw[red, ultra thick, dashed] (17.5, -4) -- (18, -4) -- (18, -3) -- (19, -3) -- (19, -2) -- (20, -2) -- (20, -1) -- (21, -1) -- (21, 0) -- (22, 0) -- (22, .5); \draw[red, ultra thick, dashed] (19, -5.5) -- (19, -5) -- (20, -5) -- (20, -4) -- (21, -4) -- (21, -3) -- (22, -3) -- (22, -2) -- (23, -2) -- (23, -1) -- (23.5, -1); \draw[red, ultra thick, dashed] (22, -5.5) -- (22, -5) -- (23, -5) -- (23, -4) -- (23.5, -4);
    \draw[blue, ultra thick, densely dotted] (17.5, 0) -- (18, 0) -- (18, .5); \draw[blue, ultra thick, densely dotted] (17.5, -3) -- (18, -3) -- (18, -2) -- (19, -2) -- (19, -1) -- (20, -1) -- (20, 0) -- (21, 0) -- (21, .5); \draw[blue, ultra thick, densely dotted] (18, -5.5) -- (18, -5) -- (19, -5) -- (19, -4) -- (20, -4) -- (20, -3) -- (21, -3) -- (21, -2) -- (22, -2) -- (22, -1) -- (23, -1) -- (23, 0) -- (23.5, 0);\draw[blue, ultra thick, densely dotted] (21, -5.5) -- (21, -5) -- (22, -5) -- (22, -4) -- (23, -4) -- (23, -3) -- (23.5, -3);
    \draw[orange, ultra thick] (17.5, -2) -- (18, -2) -- (18, -1) -- (19, -1) -- (19, 0) -- (20, 0) -- (20, .5); \draw[orange, ultra thick] (17.5, -5) -- (18, -5) -- (18, -4) -- (19, -4) -- (19, -3) -- (20, -3) -- (20, -2) -- (21, -2) -- (21, -1) -- (22, -1) -- (22, 0) -- (23, 0) -- (23, .5); \draw[orange, ultra thick] (20, -5.5) -- (20, -5) -- (21, -5) -- (21, -4) -- (22, -4) -- (22, -3) -- (23, -3) -- (23, -2) -- (23.5, -2); \draw[orange, ultra thick] (23, -5.5) -- (23, -5) -- (23.5, -5);
    \foreach \x in {18, 19, 20, 21, 22, 23}{\foreach \y in {0, -1, ..., -5}{\draw[draw = black, fill = white] (\x, \y) circle (3pt);}}
    
%C(0)
    \draw[blue, ultra thick, densely dotted] (0, -6) -- (0, -6.5) -- (1, -6.5) -- (1, -7.5) -- (0, -7.5) -- (0, -8); \draw[red, ultra thick, dashed] (1, -6) -- (1, -6.5) -- (2, -6.5) -- (2, -7.5) -- (1, -7.5) -- (1, -8); \draw[orange, ultra thick] (-.5, -6.5) -- (0, -6.5) -- (0, -7.5) -- (-.5, -7.5); \draw[orange, ultra thick] (2, -6) -- (2, -6.5) -- (2.5, -6.5); \draw[orange, ultra thick] (2.5, -7.5) -- (2, -7.5) -- (2, -8); \foreach \x in {0, 1, 2}{\foreach \y in {-6.5, -7.5}{\draw[draw = black, fill = white] (\x, \y) circle (3pt);}}

%C(1)
    \draw[blue, ultra thick, densely dotted] (0, -8.5) -- (0, -9) -- (1, -9) -- (1, -10) -- (0, -10) -- (0, -10.5); \draw[red, ultra thick, dashed] (1, -8.5) -- (1, -9) -- (2, -9) -- (2, -10) -- (1, -10) -- (1, -10.5); \draw[orange, ultra thick] (-.5, -9) -- (0, -9) -- (0, -10) -- (-.5, -10); \draw[orange, ultra thick] (2, -8.5) -- (2, -9) -- (2.5, -9); \draw[orange, ultra thick] (2.5, -10) -- (2, -10) -- (2, -10.5); \foreach \x in {0, 1, 2}{\foreach \y in {-9, -10}{\draw[draw = black, fill = white] (\x, \y) circle (3pt);}}

%C(M-1)
    \draw[blue, ultra thick, densely dotted] (0, -12.5) -- (0, -13) -- (1, -13) -- (1, -14) -- (0, -14) -- (0, -14.5); \draw[red, ultra thick, dashed] (1, -12.5) -- (1, -13) -- (2, -13) -- (2, -14) -- (1, -14) -- (1, -14.5); \draw[orange, ultra thick] (-.5, -13) -- (0, -13) -- (0, -14) -- (-.5, -14); \draw[orange, ultra thick] (2, -12.5) -- (2, -13) -- (2.5, -13); \draw[orange, ultra thick] (2.5, -14) -- (2, -14) -- (2, -14.5); \foreach \x in {0, 1, 2}{\foreach \y in {-13, -14}{\draw[draw = black, fill = white] (\x, \y) circle (3pt);}}
    
%D(0,0)
    \draw[blue, ultra thick, densely dotted] (3.5, -6) -- (3.5, -6.5) -- (4.5, -6.5) -- (4.5, -7.5) -- (3.5, -7.5) -- (3.5, -8); \draw[blue, ultra thick, densely dotted] (6.5, -6) -- (6.5, -6.5) -- (7.5, -6.5) -- (7.5, -7.5) -- (6.5, -7.5) -- (6.5, -8); \draw[red, ultra thick, dashed] (4.5, -6) -- (4.5, -6.5) -- (5.5, -6.5) -- (5.5, -7.5) -- (4.5, -7.5) -- (4.5, -8); \draw[red, ultra thick, dashed] (7.5, -6) -- (7.5, -6.5) -- (8.5, -6.5) -- (8.5, -7.5) -- (7.5, -7.5) -- (7.5, -8); \draw[orange, ultra thick] (5.5, -6) -- (5.5, -6.5) -- (6.5, -6.5) -- (6.5, -7.5) -- (5.5, -7.5) -- (5.5, -8); \draw[orange, ultra thick] (3, -6.5) -- (3.5, -6.5) -- (3.5, -7.5) -- (3, -7.5); \draw[orange, ultra thick] (8.5, -6) -- (8.5, -6.5) -- (9, -6.5); \draw[orange, ultra thick] (9, -7.5) -- (8.5, -7.5) -- (8.5, -8); \foreach \x in {3.5, 4.5, 5.5, 6.5, 7.5, 8.5}{\foreach \y in {-6.5, -7.5}{\draw[draw = black, fill = white] (\x, \y) circle (3pt);}}

%D(1,0)
    \draw[blue, ultra thick, densely dotted] (10, -6) -- (10, -6.5) -- (11, -6.5) -- (11, -7.5) -- (10, -7.5) -- (10, -8); \draw[blue, ultra thick, densely dotted] (13, -6) -- (13, -6.5) -- (14, -6.5) -- (14, -7.5) -- (13, -7.5) -- (13, -8); \draw[red, ultra thick, dashed] (11, -6) -- (11, -6.5) -- (12, -6.5) -- (12, -7.5) -- (11, -7.5) -- (11, -8); \draw[red, ultra thick, dashed] (14, -6) -- (14, -6.5) -- (15, -6.5) -- (15, -7.5) -- (14, -7.5) -- (14, -8); \draw[orange, ultra thick] (12, -6) -- (12, -6.5) -- (13, -6.5) -- (13, -7.5) -- (12, -7.5) -- (12, -8); \draw[orange, ultra thick] (9.5, -6.5) -- (10, -6.5) -- (10, -7.5) -- (9.5, -7.5); \draw[orange, ultra thick] (15, -6) -- (15, -6.5) -- (15.5, -6.5); \draw[orange, ultra thick] (15.5, -7.5) -- (15, -7.5) -- (15, -8); \foreach \x in {10, 11, 12, 13, 14, 15}{\foreach \y in {-6.5, -7.5}{\draw[draw = black, fill = white] (\x, \y) circle (3pt);}}

%D(N-1,0)
    \draw[blue, ultra thick, densely dotted] (18, -6) -- (18, -6.5) -- (19, -6.5) -- (19, -7.5) -- (18, -7.5) -- (18, -8); \draw[blue, ultra thick, densely dotted] (21, -6) -- (21, -6.5) -- (22, -6.5) -- (22, -7.5) -- (21, -7.5) -- (21, -8); \draw[red, ultra thick, dashed] (19, -6) -- (19, -6.5) -- (20, -6.5) -- (20, -7.5) -- (19, -7.5) -- (19, -8); \draw[red, ultra thick, dashed] (22, -6) -- (22, -6.5) -- (23, -6.5) -- (23, -7.5) -- (22, -7.5) -- (22, -8); \draw[orange, ultra thick] (20, -6) -- (20, -6.5) -- (21, -6.5) -- (21, -7.5) -- (20, -7.5) -- (20, -8); \draw[orange, ultra thick] (17.5, -6.5) -- (18, -6.5) -- (18, -7.5) -- (17.5, -7.5); \draw[orange, ultra thick] (23, -6) -- (23, -6.5) -- (23.5, -6.5); \draw[orange, ultra thick] (23.5, -7.5) -- (23, -7.5) -- (23, -8); \foreach \x in {18, 19, 20, 21, 22, 23}{\foreach \y in {-6.5, -7.5}{\draw[draw = black, fill = white] (\x, \y) circle (3pt);}}

%D(0,1)
    \draw[blue, ultra thick, densely dotted] (3.5, -8.5) -- (3.5, -9) -- (4.5, -9) -- (4.5, -10) -- (3.5, -10) -- (3.5, -10.5); \draw[blue, ultra thick, densely dotted] (6.5, -8.5) -- (6.5, -9) -- (7.5, -9) -- (7.5, -10) -- (6.5, -10) -- (6.5, -10.5); \draw[red, ultra thick, dashed] (4.5, -8.5) -- (4.5, -9) -- (5.5, -9) -- (5.5, -10) -- (4.5, -10) -- (4.5, -10.5); \draw[red, ultra thick, dashed] (7.5, -8.5) -- (7.5, -9) -- (8.5, -9) -- (8.5, -10) -- (7.5, -10) -- (7.5, -10.5); \draw[orange, ultra thick] (5.5, -8.5) -- (5.5, -9) -- (6.5, -9) -- (6.5, -10) -- (5.5, -10) -- (5.5, -10.5); \draw[orange, ultra thick] (3, -9) -- (3.5, -9) -- (3.5, -10) -- (3, -10); \draw[orange, ultra thick] (8.5, -8.5) -- (8.5, -9) -- (9, -9); \draw[orange, ultra thick] (9, -10) -- (8.5, -10) -- (8.5, -10.5); \foreach \x in {3.5, 4.5, 5.5, 6.5, 7.5, 8.5}{\foreach \y in {-9, -10}{\draw[draw = black, fill = white] (\x, \y) circle (3pt);}}

%D(1,1)
    \draw[blue, ultra thick, densely dotted] (10, -8.5) -- (10, -9) -- (11, -9) -- (11, -10) -- (10, -10) -- (10, -10.5); \draw[blue, ultra thick, densely dotted] (13, -8.5) -- (13, -9) -- (14, -9) -- (14, -10) -- (13, -10) -- (13, -10.5); \draw[red, ultra thick, dashed] (11, -8.5) -- (11, -9) -- (12, -9) -- (12, -10) -- (11, -10) -- (11, -10.5); \draw[red, ultra thick, dashed] (14, -8.5) -- (14, -9) -- (15, -9) -- (15, -10) -- (14, -10) -- (14, -10.5); \draw[orange, ultra thick] (12, -8.5) -- (12, -9) -- (13, -9) -- (13, -10) -- (12, -10) -- (12, -10.5); \draw[orange, ultra thick] (9.5, -9) -- (10, -9) -- (10, -10) -- (9.5, -10); \draw[orange, ultra thick] (15, -8.5) -- (15, -9) -- (15.5, -9); \draw[orange, ultra thick] (15.5, -10) -- (15, -10) -- (15, -10.5); \foreach \x in {10, 11, 12, 13, 14, 15}{\foreach \y in {-9, -10}{\draw[draw = black, fill = white] (\x, \y) circle (3pt);}}

%D(N-1,1)
    \draw[blue, ultra thick, densely dotted] (18, -8.5) -- (18, -9) -- (19, -9) -- (19, -10) -- (18, -10) -- (18, -10.5); \draw[blue, ultra thick, densely dotted] (21, -8.5) -- (21, -9) -- (22, -9) -- (22, -10) -- (21, -10) -- (21, -10.5); \draw[red, ultra thick, dashed] (19, -8.5) -- (19, -9) -- (20, -9) -- (20, -10) -- (19, -10) -- (19, -10.5); \draw[red, ultra thick, dashed] (22, -8.5) -- (22, -9) -- (23, -9) -- (23, -10) -- (22, -10) -- (22, -10.5); \draw[orange, ultra thick] (20, -8.5) -- (20, -9) -- (21, -9) -- (21, -10) -- (20, -10) -- (20, -10.5); \draw[orange, ultra thick] (17.5, -9) -- (18, -9) -- (18, -10) -- (17.5, -10); \draw[orange, ultra thick] (23, -8.5) -- (23, -9) -- (23.5, -9); \draw[orange, ultra thick] (23.5, -10) -- (23, -10) -- (23, -10.5); \foreach \x in {18, 19, 20, 21, 22, 23}{\foreach \y in {-9, -10}{\draw[draw = black, fill = white] (\x, \y) circle (3pt);}}

%D(0,M-1)
    \draw[blue, ultra thick, densely dotted] (3.5, -12.5) -- (3.5, -13) -- (4.5, -13) -- (4.5, -14) -- (3.5, -14) -- (3.5, -14.5); \draw[blue, ultra thick, densely dotted] (6.5, -12.5) -- (6.5, -13) -- (7.5, -13) -- (7.5, -14) -- (6.5, -14) -- (6.5, -14.5); \draw[red, ultra thick, dashed] (4.5, -12.5) -- (4.5, -13) -- (5.5, -13) -- (5.5, -14) -- (4.5, -14) -- (4.5, -14.5); \draw[red, ultra thick, dashed] (7.5, -12.5) -- (7.5, -13) -- (8.5, -13) -- (8.5, -14) -- (7.5, -14) -- (7.5, -14.5); \draw[orange, ultra thick] (5.5, -12.5) -- (5.5, -13) -- (6.5, -13) -- (6.5, -14) -- (5.5, -14) -- (5.5, -14.5); \draw[orange, ultra thick] (3, -13) -- (3.5, -13) -- (3.5, -14) -- (3, -14); \draw[orange, ultra thick] (8.5, -12.5) -- (8.5, -13) -- (9, -13); \draw[orange, ultra thick] (9, -14) -- (8.5, -14) -- (8.5, -14.5); \foreach \x in {3.5, 4.5, 5.5, 6.5, 7.5, 8.5}{\foreach \y in {-13, -14}{\draw[draw = black, fill = white] (\x, \y) circle (3pt);}}

%D(1,M-1)
    \draw[blue, ultra thick, densely dotted] (10, -12.5) -- (10, -13) -- (11, -13) -- (11, -14) -- (10, -14) -- (10, -14.5); \draw[blue, ultra thick, densely dotted] (13, -12.5) -- (13, -13) -- (14, -13) -- (14, -14) -- (13, -14) -- (13, -14.5); \draw[red, ultra thick, dashed] (11, -12.5) -- (11, -13) -- (12, -13) -- (12, -14) -- (11, -14) -- (11, -14.5); \draw[red, ultra thick, dashed] (14, -12.5) -- (14, -13) -- (15, -13) -- (15, -14) -- (14, -14) -- (14, -14.5); \draw[orange, ultra thick] (12, -12.5) -- (12, -13) -- (13, -13) -- (13, -14) -- (12, -14) -- (12, -14.5); \draw[orange, ultra thick] (9.5, -13) -- (10, -13) -- (10, -14) -- (9.5, -14); \draw[orange, ultra thick] (15, -12.5) -- (15, -13) -- (15.5, -13); \draw[orange, ultra thick] (15.5, -14) -- (15, -14) -- (15, -14.5); \foreach \x in {10, 11, 12, 13, 14, 15}{\foreach \y in {-13, -14}{\draw[draw = black, fill = white] (\x, \y) circle (3pt);}}

%D(N-1,M-1)
    \draw[blue, ultra thick, densely dotted] (18, -12.5) -- (18, -13) -- (19, -13) -- (19, -14) -- (18, -14) -- (18, -14.5); \draw[blue, ultra thick, densely dotted] (21, -12.5) -- (21, -13) -- (22, -13) -- (22, -14) -- (21, -14) -- (21, -14.5); \draw[red, ultra thick, dashed] (19, -12.5) -- (19, -13) -- (20, -13) -- (20, -14) -- (19, -14) -- (19, -14.5); \draw[red, ultra thick, dashed] (22, -12.5) -- (22, -13) -- (23, -13) -- (23, -14) -- (22, -14) -- (22, -14.5); \draw[orange, ultra thick] (20, -12.5) -- (20, -13) -- (21, -13) -- (21, -14) -- (20, -14) -- (20, -14.5); \draw[orange, ultra thick] (17.5, -13) -- (18, -13) -- (18, -14) -- (17.5, -14); \draw[orange, ultra thick] (23, -12.5) -- (23, -13) -- (23.5, -13); \draw[orange, ultra thick] (23.5, -14) -- (23, -14) -- (23, -14.5); \foreach \x in {18, 19, 20, 21, 22, 23}{\foreach \y in {-13, -14}{\draw[draw = black, fill = white] (\x, \y) circle (3pt);}}

\node at (16.5,-2.5){$\cdots$};
\node at (16.5,-7){$\cdots$};
\node at (16.5,-9.5){$\cdots$};
\node at (16.5,-13.5){$\cdots$};
\node at (1,-11.5){$\vdots$};
\node at (6,-11.5){$\vdots$};
\node at (12.5,-11.5){$\vdots$};
\node at (20.5,-11.5){$\vdots$};

%Vertices
\foreach \x\i in {0/0, 1/1, 2/2, 3.5/3, 4.5/4, 5.5/5, 6.5/6, 7.5/7, 8.5/8, 10/9, 11/10, 12/11, 13/12, 14/13, 15/14, 18/, 19/$n-5$, 20/, 21/$n-3$, 22/, 23/$n-1$}{\node at (\x, 1.25){\i};}
\foreach \y\j in {0/0, -1/1, -2/2, -3/3, -4/4, -5/5, -6.5/6, -7.5/7, -9/8, -10/9, -13/$m-2$, -14/$m-1$}{\node at (-1.5, \y){\j};}

\huge{
\node at (1, -17){$A$};
\node at (6,-17){$B(0)$};
\node at (12.5, -17){$B(1)$};
\node at (16.5, -17){$\cdots$};
\node at (20.5, -17){$B(N-1)$};
\node at (1, -19){$C(0)$};
\node at (6,-19){$D(0,0)$};
\node at (12.5, -19){$D(0,1)$};
\node at (16.5, -19){$\cdots$};
\node at (20.5, -19){$D(0,N-1)$};
\node at (1, -21){$C(1)$};
\node at (6,-21){$D(1,0)$};
\node at (12.5, -21){$D(1,1)$};
\node at (16.5, -21){$\cdots$};
\node at (20.5, -21){$D(1,N-1)$};
\node at (1,-23){$\vdots$};
\node at (6,-23){$\vdots$};
\node at (12.5,-23){$\vdots$};
%\node at (16.5,-23){$\ddots$};
\node at (20.5,-23){$\vdots$};
\node at (1, -25){$C(M-1)$};
\node at (6,-25){$D(M-1,0)$};
\node at (12.5, -25){$D(M-1,1)$};
\node at (16.5, -25){$\cdots$};
\node at (20.5, -25){$D(M-1,N-1)$};}
\end{tikzpicture}
\end{center}

\caption{Above: A decomposition of $C_m \sq C_n$ into three cycles in Case 7. Below: The arrangements of $A,B,C,$ and $D$ blocks in $C_m \sq C_n$ in Case 7.}\label{Case7}
\end{figure}
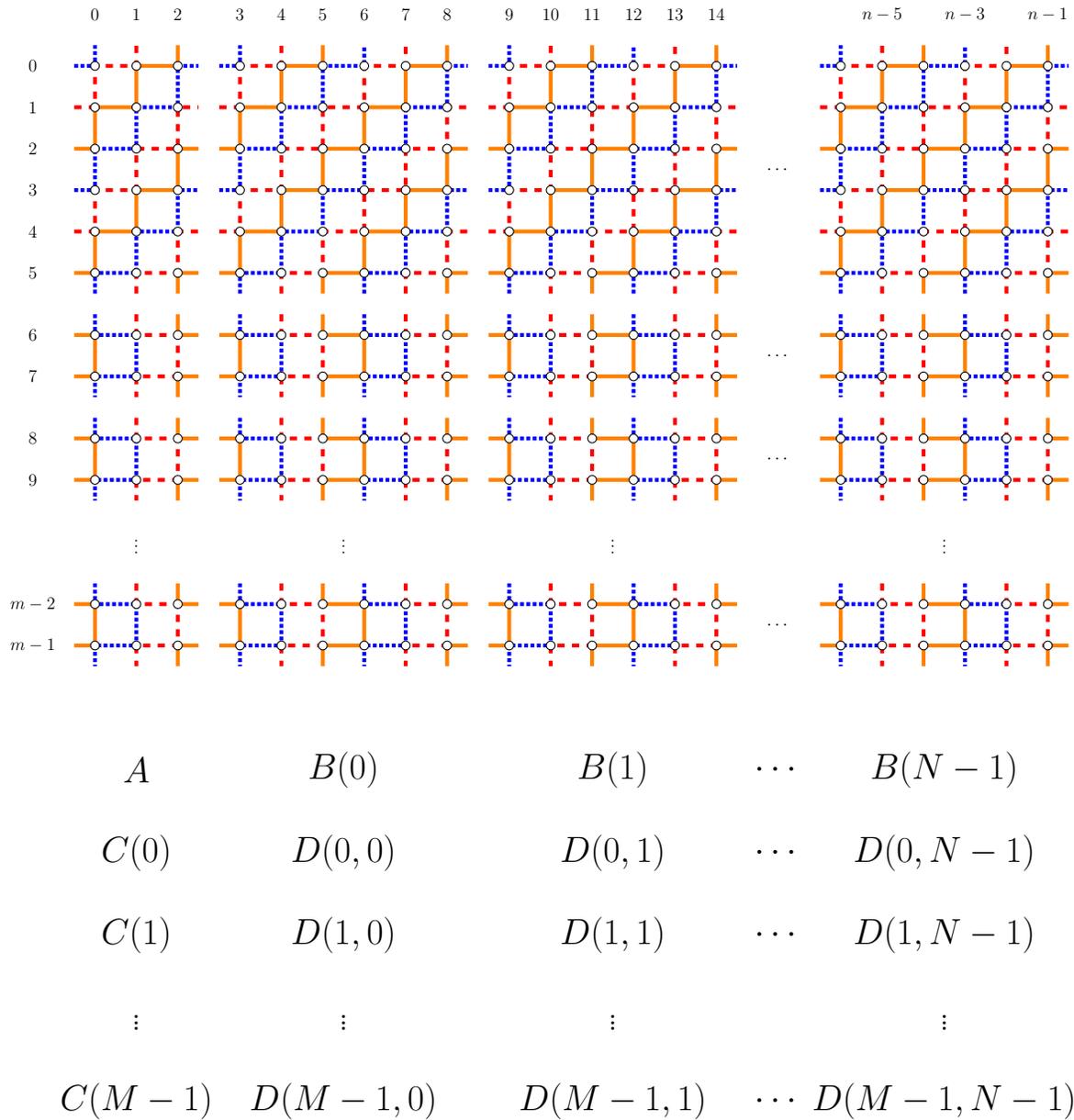

\end{proof}

\section{Decompositions of $C_{4m} \sq C_{4n}$}\label{mult4}

The goal of this section is to prove Theorem~\ref{4m4n} by generalizing a strategy used in \cite{GO22}.

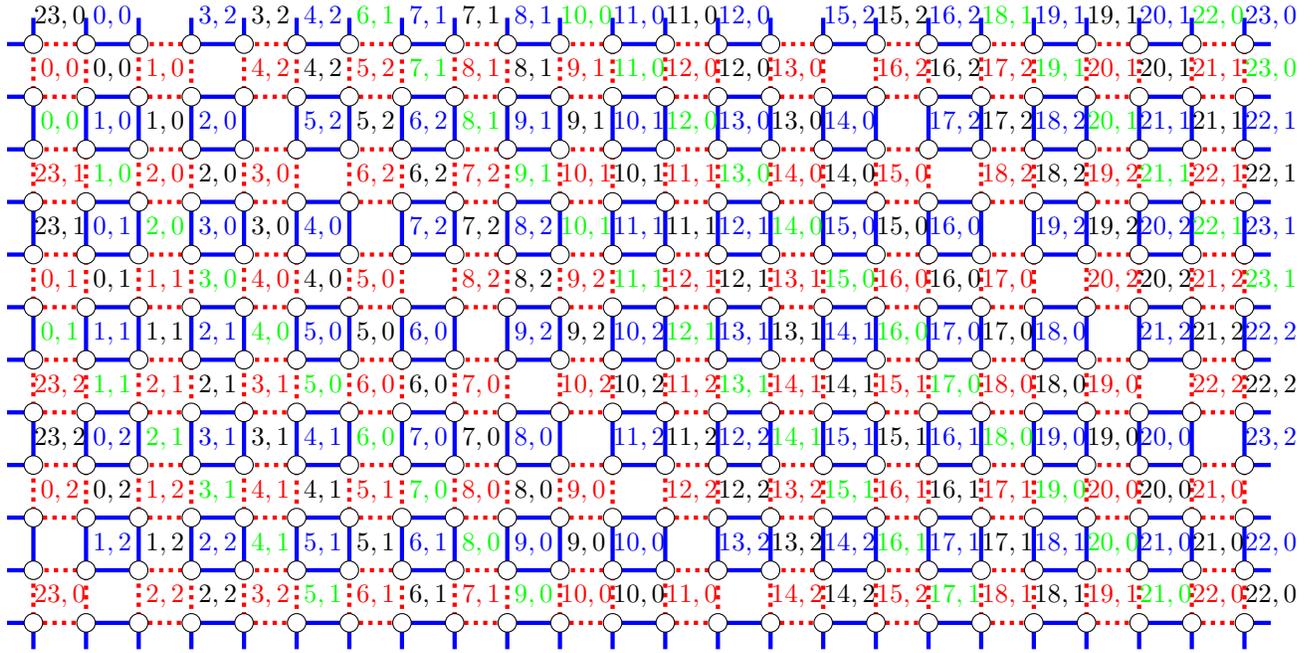
\begin{figure}

\begin{center}
    \begin{tikzpicture}[scale=.7]
        \foreach \i in {0,2,...,22}{\foreach \j in {0,2,...,10}{
            \draw [ultra thick, dotted, color=red] (\i, \j) -- (\i+1, \j) -- (\i+1, \j+1) -- (\i, \j+1) -- (\i, \j);
        }}
        \foreach \i in {1,3,...,21}{\foreach \j in {1,3,...,9}{
            \draw [ultra thick, color=blue] (\i, \j) -- (\i+1, \j) -- (\i+1, \j+1) -- (\i, \j+1) -- (\i, \j);
        }}
        \foreach \i in {1,3,...,21}{
            \draw [ultra thick, color=blue] (\i, 11.5) -- (\i, 11) -- (\i+1, 11) -- (\i +1, 11.5);
        }
        \foreach \i in {1,3,...,21}{
            \draw [ultra thick, color=blue] (\i, -0.5) -- (\i, 0) -- (\i+1, 0) -- (\i +1, -0.5);
        }
        \foreach \i in {1,3,...,9}{
            \draw [ultra thick, color=blue] (-0.5, \i) -- (0, \i) -- (0,\i+1) -- (-0.5, \i +1);
        }
        \foreach \i in {1,3,...,9}{
            \draw [ultra thick, color=blue] (23.5, \i) -- (23, \i) -- (23,\i+1) -- (23.5, \i +1);
        }
        \draw [ultra thick, color=blue] (-0.5, 0) -- (0,0) -- (0,-0.5);
        \draw [ultra thick, color=blue] (-0.5, 11) -- (0,11) -- (0,11.5);
        \draw [ultra thick, color=blue] (23.5, 0) -- (23,0) -- (23,-0.5);
        \draw [ultra thick, color=blue] (23.5, 11) -- (23,11) -- (23,11.5);

        \foreach [evaluate={
              \h=int(Mod(\i,24))+0.5;
              \v=10.5-int(Mod(\i,12));
              }]
                \i in {0,2,...,22}{
                \node [align=center] at (\h, \v) 
                {\color{red}
                \footnotesize $\i,0$};
            }   
            \foreach [evaluate={
              \h=int(Mod(\i+2,24))-0.5;
              \v=11.5-int(Mod(\i,12));
              }]
                \i in {1,3,...,23}{
                \node [align=center] at (\h, \v) 
                {\color{red}
                \footnotesize $\i,0$};
            }

            \foreach [evaluate={
              \h=int(Mod(\i,24))+0.5;
              \v=10.5-int(Mod(\i+4,12));
              }]
                \i in {0,2,...,22}{
                \node [align=center] at (\h, \v) 
                {\color{red}
                \footnotesize $\i,1$};
            }   
            \foreach [evaluate={
              \h=int(Mod(\i+2,24))-0.5;
              \v=11.5-int(Mod(\i+4,12));
              }]
                \i in {1,3,...,23}{
                \node [align=center] at (\h, \v) 
                {\color{red}
                \footnotesize $\i,1$};
            }   
            \foreach [evaluate={
              \h=int(Mod(\i,24))+0.5;
              \v=10.5-int(Mod(\i+8,12));
              }]
                \i in {0,2,...,22}{
                \node [align=center] at (\h, \v) 
                {\color{red}
                \footnotesize $\i,2$};
            }   
            \foreach [evaluate={
              \h=int(Mod(\i+2,24))-0.5;
              \v=11.5-int(Mod(\i+8,12));
              }]
                \i in {1,3,...,23}{
                \node [align=center] at (\h, \v) 
                {\color{red}
                \footnotesize $\i,2$};
            }   

            \foreach [evaluate={
              \h=int(Mod(\i+1,24))+0.5;
              \v=11.5-int(Mod(\i,12));
              }]
                \i in {0,2,...,22}{
                \node [align=center] at (\h, \v) 
                {\color{blue}
                \footnotesize $\i,0$};
            }   
            \foreach [evaluate={
              \h=int(Mod(\i-1,24))+1.5;
              \v=12.5-int(Mod(\i+2,12));
              }]
                \i in {1,3,...,23}{
                \node [align=center] at (\h, \v) 
                {\color{blue}
                \footnotesize $\i,0$};
            }   
            \foreach [evaluate={
              \h=int(Mod(\i+1,24))+0.5;
              \v=11.5-int(Mod(\i+4,12));
              }]
                \i in {0,2,...,22}{
                \node [align=center] at (\h, \v) 
                {\color{blue}
                \footnotesize $\i,1$};
            }   
            \foreach [evaluate={
              \h=int(Mod(\i-1,24))+1.5;
              \v=12.5-int(Mod(\i+6,12));
              }]
                \i in {1,3,...,23}{
                \node [align=center] at (\h, \v) 
                {\color{blue}
                \footnotesize $\i,1$};
            }   
            \foreach [evaluate={
              \h=int(Mod(\i+1,24))+0.5;
              \v=11.5-int(Mod(\i+8,12));
              }]
                \i in {0,2,...,22}{
                \node [align=center] at (\h, \v) 
                {\color{blue}
                \footnotesize $\i,2$};
            }   
            \foreach [evaluate={
              \h=int(Mod(\i-1,24))+1.5;
              \v=12.5-int(Mod(\i+10,12));
              }]
                \i in {1,3,...,23}{
                \node [align=center] at (\h, \v) 
                {\color{blue}
                \footnotesize $\i,2$};
            }

            \foreach [evaluate={
              \h=int(Mod(\i+1,24))+0.5;
              \v=11.5-int(Mod(\i+1,12));
              }]
                \i in {0,...,23}{
                \node [align=center] at (\h, \v) 
                {\color{black}
                \footnotesize $\i,0$};
            }   
            \foreach [evaluate={
              \h=int(Mod(\i+1,24))+0.5;
              \v=11.5-int(Mod(\i+5,12));
              }]
                \i in {0,...,23}{
                \node [align=center] at (\h, \v) 
                {\color{black}
                \footnotesize $\i,1$};
            }   
            \foreach [evaluate={
              \h=int(Mod(\i+1,24))+0.5;
              \v=11.5-int(Mod(\i+9,12));
              }]
                \i in {0,...,23}{
                \node [align=center] at (\h, \v) 
                {\color{black}
                \footnotesize $\i,2$};
            }  
            \foreach [evaluate={
              \h=int(Mod(\i,24))+0.5;
              \v=11.5-int(Mod(\i+2,12));
              }]
                \i in {0,...,23}{
                \node [align=center] at (\h, \v) 
                {\color{green}
                \footnotesize $\i,0$};
            }   
            \foreach [evaluate={
              \h=int(Mod(\i,24))+0.5;
              \v=11.5-int(Mod(\i+6,12));
              }]
                \i in {0,...,23}{
                \node [align=center] at (\h, \v) 
                {\color{green}
                \footnotesize $\i,1$};
            }

        \foreach \i in {0,1,...,23}{\foreach \j in {0,1,...,11}{\filldraw[fill=white, draw=black] (\i, \j) circle (5pt);}}
    
    \end{tikzpicture}
    \end{center}

\caption{The checkerboard decomposition of $C_{12} \sq C_{24}$ into cycles of length four, where the cycles $C_4(y,z)$ have dotted red or solid blue edges. Since $12=4\cdot 3$ and $24 = 4\cdot 6$, $m=3$ and $n=6$.  Thus $4\lcm(m,n) = 24$ and $\gcd(m,n) = 3$.  For $0 \le i < 24$ and $0 \le j < 3$, the cycles $R_4^{i,j}$ have dotted red edges and are labeled with a red $i,j$, the cycles $B_4^{i,j}$ have solid blue edges and are labeled with a blue $i,j$, the cycles $S_4^{i,j}$ are labeled with a black $i,j$. For $0 \le i < 24$ and $0 \le j < 2$, the cycles $T_4^{i,j}$ are labeled with a green $i,j$.
}\label{12by24checkerboard}
\end{figure}

We start with the observation that $C_m \sq C_n$ can be decomposed into cycles of length four if and only if $m$ and $n$ are both even.  The proof uses a specific ``checkerboard" $C_4$ decomposition that will be used as a starting point for decomposition into longer cycles. If $m$ and $n$ are not both four then a $C_4$-decomposition must follow the pattern shown in Figure~\ref{12by24checkerboard}.

\begin{proposition}\label{c4}
    The cartesian product $C_m \sq C_n$ can be decomposed into cycles of length four if and only if $m$ and $n$ are both even
\end{proposition}

\begin{proof}
    First note that four-cycles in $C_m \sq C_n$ can only take one of three forms. Given a vertex $(y,z) \in V(C_m \sq C_n)$, let the graph $C_4(y,z)$ be the four-cycle in $C_m \sq C_n$ with vertices $(y,z)$, $(y+1,z)$, $(y+1,z+1)$, and $(y,z+1)$.  Let the graph $V(z)$ be the four-cycle with vertices $(0,z)$, $(1,z)$, $(2,z)$ and $(3,z)$, which is only possible if $m=4$, and let $H(y)$ be the four-cycle with vertices $(y,0)$, $(y,1)$, $(y,2)$, and $(y,3)$, which is only possible only if $n=4$. Note these are the only four-cycles that are subgraphs of $C_m \sq C_n$.

    Suppose $m$ and $n$ are both even.   Then the graph $C_m \sq C_n$ is decomposed by the set of four-cycles $C_4(y,z)$ where $y+z$ is even.  
 We refer to this decomposition as the ``checkerboard'' decomposition.  The set of four-cycles $C_4(y,z)$ where $y+z$ is odd is another decomposition.  See Figure~\ref{12by24checkerboard}.

To show that no decomposition into 4-cycles is possible otherwise, suppose without loss of generality that $n$ is odd.  If a four-cycle of the form $C_4(y,z)$ is present in a decomposition, then each of the four-cycles sharing exactly one vertex with it ($C_4(y+1,z+1)$, $C_4(y+1,z-1)$, $C_4(y-1,z+1)$, and $C_4(y-1,z-1)$) must also be present, and by induction so must all four cycles $C_4(y',z')$ with $y'+z' = y+z+2a$ for some integer $a$.  But since $n$ is odd, this implies that $C_4(y,z+n-1) = C_4(y,z-1)$ must also be in the decomposition.  But since $C_4(y,z-1)$ and $C_4(y,z)$ share an edge this is a contradiction.

    Alternatively, if the decomposition contains no cycles of the form $C_4(y,z)$ then all cycles must be of the form $H(y)$ or $V(z)$, but this is impossible since $n$ is odd.
\end{proof}

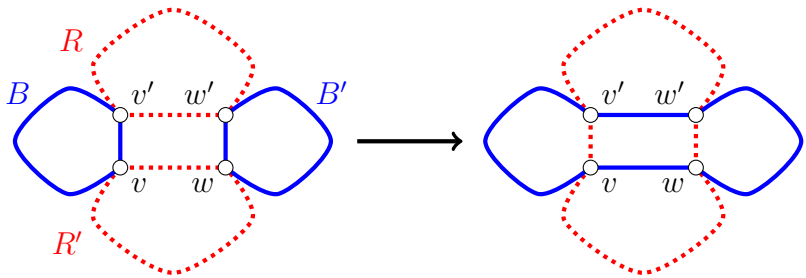
\begin{figure}

\begin{center}
    \begin{tikzpicture}[scale=.7]

        \draw[ultra thick, dotted, color=red]
        (0,0) -- (2,0);
        \draw[ultra thick, dotted, color=red]
        (0,1) -- (2,1);

        \draw[ultra thick, color=blue]
        (0,0) -- (0,1);
        \draw[ultra thick, color=blue]
        (2,1) -- (2,0);

        \draw[ultra thick, dotted, color=red]    
        plot [smooth] coordinates {(0,0) (-0.5,-1) (1,-2) (2.5,-1) (2,0)};

        \draw[ultra thick, dotted, color=red] 
        plot [smooth] coordinates {(0,1) (-0.5,2) (1,3) (2.5,2) (2,1)};

        \draw[ultra thick, color=blue]    
        plot [smooth] coordinates {(0,0) (-1,-0.5) (-2,0.5) (-1,1.5) (0,1)};

        \draw[ultra thick, color=blue]    
        plot [smooth] coordinates {(2,0) (3,-0.5) (4,0.5) (3,1.5) (2,1)};

        \node [below right] at (0,0){$v$};
        \node [below left] at (2,0){$w$};
        \node [above right] at (0,1){$v'$};
        \node [above left] at (2,1){$w'$};
        \node [above left, red] at (-.5,2){$R$};
        \node [below left, red] at (-.5,-1){$R'$};
        \node [above left, blue] at (-1.5,1){$B$};
        \node [above right, blue] at (3.5,1){$B'$};

        \foreach \i in {0,2}{\foreach \j in {0,1}{\filldraw[fill=white, draw=black] (\i, \j) circle (4pt);}}

        \draw[ultra thick, ->]
        (4.5,0.5) -- (6.5,0.5);

    \end{tikzpicture}
    \begin{tikzpicture}[scale=.7]
        \node at (-2,.5){};
        
        \draw[ultra thick, color=blue]
        (0,0) -- (2,0);
        \draw[ultra thick, color=blue]
        (0,1) -- (2,1);

        \draw[ultra thick, dotted, color=red]
        (0,0) -- (0,1);
        \draw[ultra thick, dotted, color=red]
        (2,1) -- (2,0);

        \draw[ultra thick, dotted, color=red] 
        plot [smooth] coordinates {(0,0) (-0.5,-1) (1,-2) (2.5,-1) (2,0)};

        \draw[ultra thick, dotted, color=red] 
        plot [smooth] coordinates {(0,1) (-0.5,2) (1,3) (2.5,2) (2,1)};

        \draw[ultra thick, color=blue]    
        plot [smooth] coordinates {(0,0) (-1,-0.5) (-2,0.5) (-1,1.5) (0,1)};

        \draw[ultra thick, color=blue]    
        plot [smooth] coordinates {(2,0) (3,-0.5) (4,0.5) (3,1.5) (2,1)};

        \node [below right] at (0,0){$v$};
        \node [below left] at (2,0){$w$};
        \node [above right] at (0,1){$v'$};
        \node [above left] at (2,1){$w'$};

        \foreach \i in {0,2}{\foreach \j in {0,1}{\filldraw[fill=white, draw=black] (\i, \j) circle (4pt);}}

    \end{tikzpicture}
    \end{center}

\caption{A schematic illustration of the cycle combination operation.  Recoloring the edges $vw$, $ww'$, $v'w'$, and $vv'$ "combines" $R$ and $R'$ into one long dotted red cycle and $B$ and $B'$ into one long solid blue cycle.}\label{cco}
\end{figure}

Our strategy for decomposing $C_{4m} \sq C_{4n}$ into cycles of length $4k$ is to begin with the checkerboard decomposition into $C_4$'s described in  Proposition~\ref{c4} and then apply a ``cycle combination operation" to obtain decompositions into longer cycles.  This operation is illustrated schematically in Figure~\ref{cco} and was described in \cite{GO22} as follows.

\begin{quote}
Suppose we are given vertex-disjoint red cycles $R$ and $R'$ where $v$ and $w$ are consecutive vertices in $R$ and
$v'$ and $w'$ are consecutive vertices in $R'$. We also have vertex-disjoint blue cycles $B$
and $B'$ where $v$ and $v'$ are consecutive vertices in $B$ and $w$ and $w'$ are consecutive
vertices in $B'$. Recoloring $vw$ and $v'w'$ to be blue and $vv'$ and $ww'$ to be red yields a
single red cycle whose length is the sum of the lengths of $R$ and $R'$ and a single blue cycle whose length is the sum of
the lengths of $B$ and $B'$. 
\end{quote}

Note that in the cycle combination operation, recoloring the 4-cycle $(v,w,w',v',v)$ combines the red cycles and combines the blue cycles in a way that each edge in the original cycles appears exactly once in the resulting ``combined'' cycles.  Thus applying this operation to some cycles in a decomposition yields another decomposition.

We are now ready to prove Theorem 2.

\begin{proof}[Prrof of Theorem~\ref{4m4n}]
Suppose $k$ divides $4mn$.

When describing vertices in $C_{4m} \sq C_{4n}$ using ordered pairs, as usual arithmetic is $\mod{4m}$ in the first coordinate, and $\mod{4n}$ in the second coordinate.

We begin by labeling the $C_4$'s in a checkerboard decomposition of $C_{4m} \sq C_{4n}$ in a way that will be useful when applying the cycle combination operation. We use the notation $R_4^{i,j}$, $B_4^{i,j}$, $S_4^{i,j}$, or $T_4^{i,j}$, and for these these cases, all arithmetic is $\mod{4\lcm(m,n)}$ in the first ($i$) index, and $\mod{\gcd(m,n)}$ in the second ($j$) index.  

We label each cycle in the checkerboard decomposition with $R_4^{i,j}$ or $B_4^{i,j}$ for $0 \leq i < 4\lcm(m,n)$ and $0 \leq j < \gcd(m,n)$, such that each cycle gets a unique label.  In all subsequent figures in this section, we represent edges in cycles labeled $R_4^{i,j}$ as dotted red and edges in cycles labeled $B_4^{i,j}$ as solid blue.  See Figure~\ref{12by24checkerboard}.

First, let $R_4^{0,0}$ be the cycle with vertices $\{(0,0), (0,1), (1,1), (1,0)\}$, $R_4^{1,0}$ be the cycle with vertices $\{(0,2), (0,3), (1,3), (1,2)\}$, $B_4^{0,0}$ be the cycle with vertices $\{(4m - 1, 1), (4m-1, 2), (0,2), (0,1)\}$, and $B_4^{1,0}$ be the cycle with vertices $\{(1,1), (1,2), (2,2), (2,1)\}$.

Then, for  $2 \leq i < 4\lcm(m,n)$, inductively define the cycle $C_4^{i,0}$ (where $C$ could be $R$ or $B$) to be the cycle with vertices 
\[\{(y + 2, z + 2): (y,z) \in V(C_4^{i-2, 0})\}.\]  Following our convention,  $R_4^{i+2,0}$ is two to the right and two below $R_4^{i,0}$, and $B_4^{i+2,0}$ is two to the right and two below $B_4^{i,0}$.

Finally, for $0 \leq i < 4\lcm(m,n)$ and for $1 \leq j < \gcd(m,n)$, inductively define the cycle $C_4^{i,j}$ (where $C$ could be $R$ or $B$) to be the cycle with vertices 
\[\{(y + 4j, z): (y,z) \in V(C_4^{i, 0})\}.\]  Following our convention, $R_4^{i,j+1}$ is four below $R_4^{i,j}$ and $B_4^{i,j+1}$ is four below $B_4^{i,j}$.  

We say cycles $R_4^{i,j}$ and $R_4^{i',j'}$ (or $B_4^{i,j}$ and $B_4^{i',j'}$) where $j=j'$ are in the same strand, and if $|j-j'|=1$ we say they are in consecutive strands.  Note $x = \lcm(4m,4n) = 4\lcm(m,n)$ is the smallest value where $C_4^{i,j} = C_4^{i+x,j}$, where $C$ can be $R$ or $B$, so there are $4\lcm(m,n)$ $C_4$'s of each color in each strand.  Also, $C_4(y,z)$ and $C_4(y',z)$ are in the same strand if and only if $y \equiv y' \mod{\gcd(4m,4n)}$. Specifically, if $z' = z + \ell\gcd(4m,4n)$ and $C_4(y,z) = R_4^{i,j}$, then $C_4(y',z) = R_4^{i+4nb,j}$ where $b$ is an integer such that $4nb \equiv \ell\gcd(4m,4n) \mod{4m}$. Since $\gcd(4m,4n) = 4\gcd(m,n)$ and each strand is four above the next, there are $\gcd(m,n)$ strands.  

Next, we label some $C_4$'s that are not in our decomposition with $S_4^{i,j}$ and $T_4^{i,j}$.  We will ``recolor" the edges in some of these cycles when we apply the cycle combination operation to the original red and blue cycles.

First, let $S_4^{0,0}$ be the cycle with vertices $\{(0,1),(0,2),(1,2),(1,1)\}$. Then, for $1 \leq i < 4\lcm(m,n)$, let $S_4^{i,0}$  be the cycle with vertices 
\[\{(y+1, z+1) : (y,z) \in V(S_4^{i-1, 0})\}.\]

For $0 \leq i < 4\lcm(m,n)$ and $1 \leq j < \gcd(m,n)$,
let $S_4^{i,j}$  be the cycle with vertices 
\[\{(y+4j, z) : (y,z) \in V(S_4^{i, 0})\}.\]

Following our convention, $S_4^{i+1,0}$ is one below and one to the right of $S_4^{i,0}$ and $S_4^{i,j+1}$ is four below $S_4^{i,j}$.

Finally, for $0 \leq i < 4\lcm(m,n)$ and $0 \leq j < \gcd(m,n)$, let $T_4^{i,j}$ be the cycle with vertices
\[\{(y+1 , z-1) : (y,z) \in V(S_4^{i,j})\}.\]
Following our convention, $T_4^{i,j}$ is one below and one to the left of $S_4^{i,j}$.

We now describe how to use the cycle combination operation in two phases to transform our checkerboard $C_4$-decomposition of $C_{4m} \sq C_{4n}$ into a $C_{4k}$ decomposition.

%Note that recoloring any $S_4^{i,j}$ performs a cycle combination operation on red and blue cycles in the same strand, and recoloring any $T_4^{i,j,}$ performs a cycle combination operation on red and blue cycles in consecutive strands. %TODO: Put this observation here or later?

Let $g$ and $h$ be positive integers, with $g \le h$, $g$ dividing $\gcd(m,n)$ and $h$ dividing $4\lcm(m,n)$, such that $k=gh$. We note that  this is always possible: If $k$ divides $4mn$ then there is some integer $s$ such that $ks=4mn=4\lcm(m,n)\gcd(m,n)$ for some integer $s$.  Thus the factors of $k$ can be split between $4\lcm(m,n)$ and $\gcd(m,n)$, i.e. there must exist $h$ dividing $4\lcm(m,n)$ and $g$ dividing $\gcd(m,n)$ such that $k=gh$. Further, every factor of $\gcd(m,n)$ is a factor of $\lcm(m,n)$ so $h$ can be chosen at least as large as $g$.  Specifically, $h$ will be 1 only if $k=1$.   In the case where $k=1$ the decomposition of $C_{4m} \sq C_{4n}$ into $C_{4k}$'s is the checkerboard decomposition.  So for the remainder of the proof we assume $h>1$.

%To get the desired decomposition, perform recolorings on the following four-cycles:  

Phase 1: In Phase 1 we perform the cycle combination operation on groups of $h$ cycles in the same strand to transform our checkerboard decomposition into a $C_{4h}$ decomposition. See Figures~\ref{12by24c12} and \ref{12by24c16} for examples where $h=3$ and $h=4$, respectively.

Note that recoloring the edges of $S_4^{i,j}$ performs a cycle combination operation on the red cycles $R^{i,j}$ and $R^{i+1,j}$ and the blue cycles $B^{i,j}$ and $B^{i+1,j}$.  

Recolor $S_4^{i,j}$ for all $i,j$ where $h$ does not divide $i+1$ and for all $j$.  This combines four-cycles $h$ at a time within strands and yields cycles of length $4h$.

Specifically, for $0 \le \alpha < 4\lcm(m,n)/h$, $0 \le j < \gcd(m,n)$, this creates cycles $C_{4h}^{\alpha,j}$ (where $C$ can be $R$ or $B$) where $C_{4h}^{\alpha,j}$ is a cycle of length $4h$, the result of combining cycles $C_4^{i,j}$ where $\alpha h \le i <\alpha(h+1)$.

Phase 2: In Phase 2 we perform the cycle combination operation on groups of $g$ cycles in consecutive strands to transform our $C_{4h}$ decomposition into a $C_{4gh}$ decomposition. See Figures~\ref{12by24c36} and \ref{12by24c48} for examples where $h=3,g=3$ and $h=4,g=3$, respectively.  Recall that for these examples, the results of Phase 1 are illustrated in Figures~\ref{12by24c12} and \ref{12by24c16}.

If $i$ is even and $0 \le j < \gcd(m,n)-1$, recoloring the edges of $T_4^{i,j}$ performs a cycle combination operation on the red cycles $R_4^{i,j}$ and $R_4^{i-1,j+1}$ and the blue cycles $B_4^{i+1,j}$ and $B_4^{i-2,j+1}$. 

If $i$ is odd and $0 \le j < \gcd(m,n)-1$, recoloring the edges of $T_4^{i,j}$ performs a cycle combination operation on the red cycles $R_4^{i+1,j}$ and $R_4^{i-2,j+1}$ and the blue cycles $B_4^{i,j}$ and $B_4^{i-1,j+1}$.

Thus if $h >1$, $0 \le \alpha < 4\lcm(m,n)/h$, and $0 \le j < \gcd(m,n)-1$, recoloring the edges of $T_4^{\alpha h,j}$ performs a cycle combination operation on the red cycles $R_{4h}^{\alpha ,j}$ and $R_{4h}^{\alpha  -1,j+1}$ and the blue cycles $B_{4h}^{\alpha ,j}$ and $B_{4h}^{\alpha-1,j+1}$.  (Note that $C_{4h}^{\alpha,j}$ is the result of combining cycles $C_4^{\alpha h,j}$ with $C_4^{\alpha h + 1,j}$ (and perhaps some other $C_4$'s), and $C_{4h}^{\alpha - 1,j+1}$ is the result of combining cycles $C_4^{\alpha h - 1,j+1}$ with $C_4^{\alpha h -2,j+1}$ (and perhaps some other $C_4$'s), where $C$ can be $R$ or $B$.)

So in Phase 2 we recolor $T_4^{i,j}$ for all $i,j$ where $h$ divides $i$ and $g$ does not divide $j+1$.  This combines the $C_{4h}$'s from Phase 1 $g$ at a time across strands, and yields cycles of length $4gh=4k$.

Thus, for $0 \le \alpha < 4\lcm(m,n)/h$, $0 \le \beta <\gcd(m,n)/g$, Phase 2 yields cycles $C_{4gh}^{\alpha,\beta}$ (where $C$ can be $R$ or $B$) where $C_{4gh}^{\alpha,\beta}$ is the result of combining the following set of cycles: \[\{C_{4h}^{\alpha-\ell,\beta g +\ell} ~:~  0 \le \ell < g\}.\]

Since all of the resulting cycles were obtained by applying the cycle combination to the cycles in the checkerboard decomposition, all edges are in exactly one cycle in the resulting decomposition.
\end{proof}

\begin{figure}

\begin{center}
    \begin{tikzpicture}[scale=.7]

    \foreach [evaluate={
              \h=int(Mod(\i+1,24));
              \v=12-int(Mod(\i+1,12));
              }]
                \i in {0,1,3,4, 6,7, 9, 10, 12, 13, 15, 16, 18, 19, 21, 22}{
                \fill [gray!15] (\h,\v) rectangle (\h+1,\v-1);;
            }
    \foreach [evaluate={
              \h=int(Mod(\i+1,24));
              \v=12-int(Mod(\i+5,12));
              }]
                \i in {0,1,3,4, 6,7, 9, 10, 12, 13, 15, 16, 18, 19, 21, 22}{
                \fill [gray!15] (\h,\v) rectangle (\h+1,\v-1);;
            }
    \foreach [evaluate={
              \h=int(Mod(\i+1,24));
              \v=12-int(Mod(\i+9,12));
              }]
                \i in {0,1,3,4, 6,7, 9, 10, 12, 13, 15, 16, 18, 19, 21, 22}{
                \fill [gray!15] (\h,\v) rectangle (\h+1,\v-1);;
            }   
            
        \foreach \i in {0,3}{
        \foreach \j in {0,1,2}{
        \draw [ultra thick, dotted, color=red] (0+4*\i,11-4*\j) -- (3+4*\i,11-4*\j) -- (3+4*\i,8-4*\j) --  (2+4*\i,8-4*\j) -- (2+4*\i,10-4*\j) -- (0+4*\i,10-4*\j) -- (0+4*\i,11-4*\j);
        }}

        \foreach \i in {0,3}{
        \foreach \j in {0,1}{
        \draw [ultra thick, dotted, color=red] (4+4*\i,9-4*\j) -- (5+4*\i,9-4*\j) -- (5+4*\i,7-4*\j) --  (7+4*\i,7-4*\j) -- (7+4*\i,6-4*\j) -- (4+4*\i,6-4*\j) -- (4+4*\i,9-4*\j);
        }}

        \foreach \i in {0,3}{
        \foreach \j in {0,1}{
        \draw [ultra thick, dotted, color=red] (6+4*\i,9-4*\j) -- (9+4*\i,9-4*\j) -- (9+4*\i,6-4*\j) --  (8+4*\i,6-4*\j) -- (8+4*\i,8-4*\j) -- (6+4*\i,8-4*\j) -- (6+4*\i,9-4*\j);
        }}

        \foreach \i in {0}{
        \foreach \j in {0,1,2}{
        \draw [ultra thick, dotted, color=red] (10+4*\i,11-4*\j) -- (11+4*\i,11-4*\j) -- (11+4*\i,9-4*\j) --  (13+4*\i,9-4*\j) -- (13+4*\i,8-4*\j) -- (10+4*\i,8-4*\j) -- (10+4*\i,11-4*\j);
        }}

        \foreach \i in {3}{
        \foreach \j in {0,1,2}{
        \draw [ultra thick, dotted, color=red] (11.5+4*\i,8-4*\j) --  (10+4*\i,8-4*\j) -- (10+4*\i,11-4*\j) -- (11+4*\i,11-4*\j) -- (11+4*\i,9-4*\j) -- (11.5+4*\i,9-4*\j) ;
        }}

        \foreach \i in {0}{
        \foreach \j in {0,1,2}{
        \draw [ultra thick, dotted, color=red] (-0.5+4*\i,8-4*\j) --  (1+4*\i,8-4*\j) -- (1+4*\i,9-4*\j) -- (-0.5+4*\i,9-4*\j);
        }}

        \foreach \i in {0,3}{
        \foreach \j in {0}{
        \draw [ultra thick, dotted, color=red] (4+4*\i,11.5-4*\j) --  (4+4*\i,10-4*\j) -- (7+4*\i,10-4*\j) -- (7+4*\i,11-4*\j)  -- (5+4*\i,11-4*\j) -- (5+4*\i,11.5-4*\j);
        }}
        \foreach \i in {0,3}{
        \foreach \j in {0}{
        \draw [ultra thick, dotted, color=red] (4+4*\i,-0.5-4*\j) --  (4+4*\i,1-4*\j) -- (5+4*\i,1-4*\j) -- (5+4*\i,-0.5-4*\j);
        }}

        \foreach \i in {0,3}{
        \foreach \j in {0}{
        \draw [ultra thick, dotted, color=red] (9+4*\i,-0.5-4*\j) --  (9+4*\i,1-4*\j) -- (6+4*\i,1-4*\j) -- (6+4*\i,0-4*\j)  -- (8+4*\i,0-4*\j) -- (8+4*\i,-0.5-4*\j);
        }}
        \foreach \i in {0,3}{
        \foreach \j in {0}{
        \draw [ultra thick, dotted, color=red] (8+4*\i,11.5-4*\j) --  (8+4*\i,10-4*\j) -- (9+4*\i,10-4*\j) -- (9+4*\i,11.5-4*\j);
        }}

        \foreach \i in {0,3}{
        \foreach \j in {0,1}{
        \draw [ultra thick, color=blue] (3+4*\i,8-4*\j) -- (6+4*\i,8-4*\j) -- (6+4*\i,5-4*\j) --  (5+4*\i,5-4*\j) -- (5+4*\i,7-4*\j) -- (3+4*\i,7-4*\j) -- (3+4*\i,8-4*\j);
        }}

        \foreach \i in {0}{
        \foreach \j in {0,1}{
        \draw [ultra thick, color=blue] (9+4*\i,10-4*\j) -- (12+4*\i,10-4*\j) -- (12+4*\i,7-4*\j) --  (11+4*\i,7-4*\j) -- (11+4*\i,9-4*\j) -- (9+4*\i,9-4*\j) -- (9+4*\i,10-4*\j);
        }}

        \foreach \i in {3}{
        \foreach \j in {0,1}{
        \draw [ultra thick, color=blue] (11.5+4*\i,10-4*\j)  -- (9+4*\i,10-4*\j) -- (9+4*\i,9-4*\j) --  (11+4*\i,9-4*\j) -- (11+4*\i,7-4*\j) -- (11.5+4*\i,7-4*\j);
        }}
        \foreach \i in {0}{
        \foreach \j in {0,1}{
        \draw [ultra thick, color=blue] (-0.5+4*\i,10-4*\j)  -- (0+4*\i,10-4*\j) -- (0+4*\i,7-4*\j) --  (-0.5+4*\i,7-4*\j);
        }}

        \foreach \i in {0,3}{
        \foreach \j in {0}{
        \draw [ultra thick, color=blue] (1+4*\i,11.5-4*\j) --  (1+4*\i,9-4*\j) -- (4+4*\i,9-4*\j) -- (4+4*\i,10-4*\j)  -- (2+4*\i,10-4*\j) -- (2+4*\i,11.5-4*\j);
        }}
        \foreach \i in {0,3}{
        \foreach \j in {0}{
        \draw [ultra thick, color=blue] (1+4*\i,-0.5-4*\j) --  (1+4*\i,0-4*\j) -- (2+4*\i,0-4*\j) -- (2+4*\i,-0.5-4*\j);
        }}

        \foreach \i in {0,3}{
        \foreach \j in {0,1}{
        \draw [ultra thick, color=blue] (7+4*\i,10-4*\j) --  (7+4*\i,7-4*\j) -- (10+4*\i,7-4*\j) -- (10+4*\i,8-4*\j)  -- (8+4*\i,8-4*\j) -- (8+4*\i,10-4*\j) -- (7+4*\i,10-4*\j) ;
        }}

        \foreach \i in {0,3}{
        \foreach \j in {0,1}{
        \draw [ultra thick, color=blue] (1+4*\i,8-4*\j) --  (1+4*\i,5-4*\j) -- (4+4*\i,5-4*\j) -- (4+4*\i,6-4*\j)  -- (2+4*\i,6-4*\j) -- (2+4*\i,8-4*\j) -- (1+4*\i,8-4*\j) ;
        }}

        \foreach \i in {0,3}{
        \draw [ultra thick, color=blue] 
        (3 + 4*\i, 11.5) --  (3+ 4*\i,11) -- (5+ 4*\i,11) -- (5+ 4*\i, 9) -- (6+ 4*\i,9) -- (6+ 4*\i,11.5) ;
        }
        \foreach \i in {0,3}{
        \draw [ultra thick, color=blue] 
        (3 + 4*\i, -0.5) --  (3+ 4*\i,0) -- (6+ 4*\i,0) -- (6+ 4*\i,-0.5) ;
        }

        \foreach \i in {0,3}{
        \draw [ultra thick, color=blue] 
        (7 + 4*\i, -0.5) --  (7+ 4*\i,2) -- (8+ 4*\i,2) -- (8+ 4*\i, 0) -- (10+ 4*\i,0) -- (10+ 4*\i,-0.5) ;
        }
        \foreach \i in {0,3}{
        \draw [ultra thick, color=blue] 
        (7 + 4*\i, 11.5) --  (7+ 4*\i,11) -- (10+ 4*\i,11) -- (10+ 4*\i,11.5) ;
        }

        \draw [ultra thick, color=blue]
        (11, -0.5) -- (11,1) -- (9,1) -- (9,2) -- (12,2) --(12,-0.5); 
        \draw [ultra thick, color=blue]
        (11, 11.5) -- (11,11) -- (12,11) --(12,11.5); 

        \draw [ultra thick, color=blue]
        (23, -0.5) -- (23, 1) -- (21,1) -- (21,2) -- (23.5,2);
        \draw [ultra thick, color=blue]
        (23, 11.5) -- (23, 11) -- (23.5,11);
        \draw [ultra thick, color=blue]
        (0, 11.5) -- (0, 11) -- (-0.5,11);
        \draw [ultra thick, color=blue]
        (0, -0.5) -- (0, 2) -- (-0.5,2);

        \foreach [evaluate={
              \h=int(Mod(\i,24))+0.5;
              \v=10.5-int(Mod(\i,12));
              }]
                \i in {0,2,...,22}{
                \node [align=center] at (\h, \v) 
                {\color{red}
                \footnotesize $\i,0$};
            }   
            \foreach [evaluate={
              \h=int(Mod(\i+2,24))-0.5;
              \v=11.5-int(Mod(\i,12));
              }]
                \i in {1,3,...,23}{
                \node [align=center] at (\h, \v) 
                {\color{red}
                \footnotesize $\i,0$};
            }

            \foreach [evaluate={
              \h=int(Mod(\i,24))+0.5;
              \v=10.5-int(Mod(\i+4,12));
              }]
                \i in {0,2,...,22}{
                \node [align=center] at (\h, \v) 
                {\color{red}
                \footnotesize $\i,1$};
            }   
            \foreach [evaluate={
              \h=int(Mod(\i+2,24))-0.5;
              \v=11.5-int(Mod(\i+4,12));
              }]
                \i in {1,3,...,23}{
                \node [align=center] at (\h, \v) 
                {\color{red}
                \footnotesize $\i,1$};
            }   
            \foreach [evaluate={
              \h=int(Mod(\i,24))+0.5;
              \v=10.5-int(Mod(\i+8,12));
              }]
                \i in {0,2,...,22}{
                \node [align=center] at (\h, \v) 
                {\color{red}
                \footnotesize $\i,2$};
            }   
            \foreach [evaluate={
              \h=int(Mod(\i+2,24))-0.5;
              \v=11.5-int(Mod(\i+8,12));
              }]
                \i in {1,3,...,23}{
                \node [align=center] at (\h, \v) 
                {\color{red}
                \footnotesize $\i,2$};
            }   

            \foreach [evaluate={
              \h=int(Mod(\i+1,24))+0.5;
              \v=11.5-int(Mod(\i,12));
              }]
                \i in {0,2,...,22}{
                \node [align=center] at (\h, \v) 
                {\color{blue}
                \footnotesize $\i,0$};
            }   
            \foreach [evaluate={
              \h=int(Mod(\i-1,24))+1.5;
              \v=12.5-int(Mod(\i+2,12));
              }]
                \i in {1,3,...,23}{
                \node [align=center] at (\h, \v) 
                {\color{blue}
                \footnotesize $\i,0$};
            }   
            \foreach [evaluate={
              \h=int(Mod(\i+1,24))+0.5;
              \v=11.5-int(Mod(\i+4,12));
              }]
                \i in {0,2,...,22}{
                \node [align=center] at (\h, \v) 
                {\color{blue}
                \footnotesize $\i,1$};
            }   
            \foreach [evaluate={
              \h=int(Mod(\i-1,24))+1.5;
              \v=12.5-int(Mod(\i+6,12));
              }]
                \i in {1,3,...,23}{
                \node [align=center] at (\h, \v) 
                {\color{blue}
                \footnotesize $\i,1$};
            }   
            \foreach [evaluate={
              \h=int(Mod(\i+1,24))+0.5;
              \v=11.5-int(Mod(\i+8,12));
              }]
                \i in {0,2,...,22}{
                \node [align=center] at (\h, \v) 
                {\color{blue}
                \footnotesize $\i,2$};
            }   
            \foreach [evaluate={
              \h=int(Mod(\i-1,24))+1.5;
              \v=12.5-int(Mod(\i+10,12));
              }]
                \i in {1,3,...,23}{
                \node [align=center] at (\h, \v) 
                {\color{blue}
                \footnotesize $\i,2$};
            }

            \foreach [evaluate={
              \h=int(Mod(\i+1,24))+0.5;
              \v=11.5-int(Mod(\i+1,12));
              }]
                \i in {0,1,3,4, 6,7, 9, 10, 12, 13, 15, 16, 18, 19, 21, 22}{
                \node [align=center] at (\h, \v) 
                {\color{black}
                \footnotesize $\i,0$};
            }   
            \foreach [evaluate={
              \h=int(Mod(\i+1,24))+0.5;
              \v=11.5-int(Mod(\i+5,12));
              }]
                \i in {0,1,3,4, 6,7, 9, 10, 12, 13, 15, 16, 18, 19, 21, 22}{
                \node [align=center] at (\h, \v) 
                {\color{black}
                \footnotesize $\i,1$};
            }   
            \foreach [evaluate={
              \h=int(Mod(\i+1,24))+0.5;
              \v=11.5-int(Mod(\i+9,12));
              }]
                \i in {0,1,3,4, 6,7, 9, 10, 12, 13, 15, 16, 18, 19, 21, 22}{
                \node [align=center] at (\h, \v) 
                {\color{black}
                \footnotesize $\i,2$};
            }  

    \junk{  
            \foreach [evaluate={
              \h=int(Mod(\i,24))+0.5;
              \v=11.5-int(Mod(\i+2,12));
              }]
                \i in {0,...,23}{
                \node [align=center] at (\h, \v) 
                {\color{green}
                \footnotesize $\i,0$};
            }   
            \foreach [evaluate={
              \h=int(Mod(\i,24))+0.5;
              \v=11.5-int(Mod(\i+6,12));
              }]
                \i in {0,...,23}{
                \node [align=center] at (\h, \v) 
                {\color{green}
                \footnotesize $\i,1$};
            } 
    }            

        \foreach \i in {0,1,...,23}{\foreach \j in {0,1,...,11}{\filldraw[fill=white, draw=black] (\i, \j) circle (5pt);}}
    
    \end{tikzpicture}
    \end{center}

\caption{Decomposition of $C_{12} \sq C_{24}$ into $C_{12}$'s. Since $12=4\cdot 3$, $k=3$, and we can choose $h=3$ and $g=1$. Cycles where the cycle combination operation has been applied are shaded in light gray. \newline
For example, applying the cycle combination operation on $S_4^{0,0}$ and $S_4^{1,0}$ combines the red cycles $R_4^{0,0}$, $R_4^{1,0}$, and $R_4^{2,0}$ into a single red cycle $R_{12}^{0,0}$ of length 12 and the blue cycles $B_4^{0,0}$, $B_4^{1,0}$, and $B_4^{2,0}$ into a single blue cycle $B_{12}^{0,0}$ of length 12. Similarly, applying the cycle combination operation on $S_4^{9,1}$ and $S_4^{10,1}$ combines the red cycles $R_4^{9,1}$, $R_4^{10,1}$, and $R_4^{11,1}$ into a single red cycle $R_{12}^{3,1}$ of length 12 and the blue cycles $B_4^{9,1}$, $B_4^{10,1}$, and $B_4^{11,1}$ into a single blue cycle $B_{12}^{3,1}$ of length 12.  }\label{12by24c12}
\end{figure}

\begin{figure}

\begin{center}
    \begin{tikzpicture}[scale=.7]

    \foreach [evaluate={
              \h=int(Mod(\i+1,24));
              \v=12-int(Mod(\i+1,12));
              }]
                \i in {0,1,3,4, 6,7, 9, 10, 12, 13, 15, 16, 18, 19, 21, 22}{
                \fill [gray!15] (\h,\v) rectangle (\h+1,\v-1);;
            }
    \foreach [evaluate={
              \h=int(Mod(\i+1,24));
              \v=12-int(Mod(\i+5,12));
              }]
                \i in {0,1,3,4, 6,7, 9, 10, 12, 13, 15, 16, 18, 19, 21, 22}{
                \fill [gray!15] (\h,\v) rectangle (\h+1,\v-1);;
            }
    \foreach [evaluate={
              \h=int(Mod(\i+1,24));
              \v=12-int(Mod(\i+9,12));
              }]
                \i in {0,1,3,4, 6,7, 9, 10, 12, 13, 15, 16, 18, 19, 21, 22}{
                \fill [gray!15] (\h,\v) rectangle (\h+1,\v-1);;
            }

    \foreach [evaluate={
              \h=int(Mod(\i,24));
              \v=12-int(Mod(\i+6,12));
              }]
                \i in {0,3,...,21}{
                \fill [gray!15] (\h,\v) rectangle (\h+1,\v-1);
            } 
    \foreach [evaluate={
              \h=int(Mod(\i,24));
              \v=12-int(Mod(\i+2,12));
              }]
                \i in {0,3,...,21}{
                \fill [gray!15] (\h,\v) rectangle (\h+1,\v-1);
            }

        \foreach \i in {0,2}{
        \draw [ultra thick, dotted, color=red]
        (0+6*\i, 3) -- (0+6*\i, 2) -- (2+6*\i,2) -- (2+6*\i,0) -- (3+6*\i,0) -- (3+6*\i,2) -- (7+6*\i,2) -- (7+6*\i,4) -- (8+6*\i,4) -- (8+6*\i,2) -- (9+6*\i,2) -- (9+6*\i,5) -- (6+6*\i,5) -- (6+6*\i,3) -- (5+6*\i,3) -- (5+6*\i,5) -- (4+6*\i,5) -- (4+6*\i,3) -- (0+6*\i,3);
        }

        \foreach \i in {1}{
        \draw [ultra thick, dotted, color=red]
        (0+6*\i, 3+6*\i) -- (0+6*\i, 2+6*\i) -- (2+6*\i,2+6*\i) -- (2+6*\i,0+6*\i) -- (3+6*\i,0+6*\i) -- (3+6*\i,2+6*\i) -- (7+6*\i,2+6*\i) -- (7+6*\i,4+6*\i) -- (8+6*\i,4+6*\i) -- (8+6*\i,2+6*\i) -- (9+6*\i,2+6*\i) -- (9+6*\i,5+6*\i) -- (6+6*\i,5+6*\i) -- (6+6*\i,3+6*\i) -- (5+6*\i,3+6*\i) -- (5+6*\i,5+6*\i) -- (4+6*\i,5+6*\i) -- (4+6*\i,3+6*\i) -- (0+6*\i,3+6*\i);
        }

        \draw [ultra thick, dotted, color=red]
        (10,7) -- (10,4) -- (13,4) -- (13,6) -- (14,6) -- (14,4) -- (15,4) -- (15,6) -- (19,6) -- (19,7) -- (17,7) -- (17,9) -- (16,9) -- (16,7) -- (12,7) -- (12,5) -- (11,5) -- (11,7) -- (10,7);

        \draw [ultra thick, dotted, color=red]
        (23.5,5) -- (23,5) -- (23,7) -- (22,7) -- (22,4) -- (23.5,4);
        \draw [ultra thick, dotted, color=red]
        (-0.5,5) -- (0,5) -- (0,7) -- (4,7) -- (4,9) -- (5,9) -- (5,7) -- (7,7) -- (7,6) -- (3,6) -- (3,4) -- (2,4) -- (2,6) -- (1,6) -- (1,4) -- (-0.5,4);

        \draw [ultra thick, dotted, color=red]
        (-0.5,9) -- (0,9) -- (0,11) -- (3,11) -- (3,8) --(2,8) -- (2,10) -- (1,10) -- (1,8) -- (-0.5,8);
        \draw [ultra thick, dotted, color=red]
        (23.5,9) -- (23,9) -- (23,11) -- (22,11) -- (22,9) -- (18,9) -- (18,8) -- (20,8) -- (20,6) -- (21,6) -- (21,8) -- (23.5,8);

        \foreach \i in {0,1}{
        \draw [ultra thick, dotted, color=red]
        (4+12*\i,-0.5) -- (4+12*\i,1) -- (5+12*\i,1) -- (5+12*\i, -0.5);
        \draw [ultra thick, dotted, color=red]
        (4+12*\i,11.5) -- (4+12*\i,10) -- (7+12*\i,10) -- (7+12*\i,11.5);
        \draw [ultra thick, dotted, color=red]
        (5+12*\i,11.5) -- (5+12*\i,11) -- (6+12*\i,11) -- (6+12*\i,11.5);
        \draw [ultra thick, dotted, color=red]
        (9+12*\i,11.5) -- (9+12*\i,10) -- (8+12*\i,10) -- (8+12*\i,11.5);
        \draw [ultra thick, dotted, color=red]
        (7+12*\i,-0.5) -- (7+12*\i,0) -- (8+12*\i,0) -- (8+12*\i,-0.5);
        }

        \draw [ultra thick, dotted, color=red]
        (6,-0.5) -- (6,1) -- (10,1) -- (10,3) -- (11,3) -- (11,1) -- (13,1) -- (13,0) -- (9,0) -- (9,-0.5);

        \draw [ultra thick, dotted, color=red]
        (18,-0.5) -- (18,1) -- (22,1) -- (22,3) -- (23,3) -- (23,1) -- (23.5,1);
        \draw [ultra thick, dotted, color=red]
        (-0.5, 1) -- (1,1) -- (1,0) -- (-0.5,0);
        \draw [ultra thick, dotted, color=red]
        (23.5, 0) -- (21,0) -- (21,-0.5);

        \foreach \i in {0,1}{
        \draw [ultra thick, color=blue]
        (1+12*\i, 4) -- (2+12*\i,4) -- (2+12*\i,2) -- (3+12*\i,2) -- (3+12*\i,4) -- (7+12*\i,4) -- (7+12*\i,6) -- (8+12*\i,6) -- (8+12*\i,4) -- (10+12*\i,4) -- (10+12*\i,3) -- (6+12*\i,3) -- (6+12*\i,1) -- (5+12*\i,1) -- (5+12*\i,3) -- (4+12*\i,3) -- (4+12*\i,1) -- (1+12*\i,1) -- (1+12*\i,4);
        }

        \foreach \i in {0.5}{
        \draw [ultra thick, color=blue]
         (8+12*\i,6+11*\i) -- (8+12*\i,4+12*\i) -- (10+12*\i,4+12*\i) -- (10+12*\i,3+12*\i) -- (6+12*\i,3+12*\i) -- (6+12*\i,1+12*\i) -- (5+12*\i,1+12*\i) -- (5+12*\i,3+12*\i) -- (4+12*\i,3+12*\i) -- (4+12*\i,1+12*\i) -- (1+12*\i,1+12*\i) -- (1+12*\i,4+12*\i) -- (1+12*\i, 4+12*\i) -- (2+12*\i,4+12*\i) -- (2+12*\i,2+12*\i) -- (3+12*\i,2+12*\i) -- (3+12*\i,4+12*\i) -- (7+12*\i,4+12*\i) -- (7+12*\i,6+11*\i);
         }
         
        \foreach \i in {0,1}{
        \draw [ultra thick, color=blue]
        (13-12*\i, -0.5)-- (13-12*\i,0) -- (14-12*\i,0) -- (14-12*\i, -0.5);
        }

        \draw [ultra thick, color=blue]
        (2,11.5) -- (2,10) -- (4,10) -- (4,9) -- (0,9) -- (0,7) -- (-0.5,7);
        \draw [ultra thick, color=blue]
        (23.5,7) -- (23,7) -- (23,9)-- (22,9) -- (22,7) -- (19,7) --(19,10) -- (20,10) -- (20,8) -- (21,8) -- (21,10) -- (23.5,10); 
        \draw [ultra thick, color=blue]
        (-0.5,10) -- (1,10) -- (1,11.5);

        \draw [ultra thick, color=blue]
        (9,6) -- (13,6) -- (13,8) -- (14,8) -- (14,6) -- (15,6) -- (15,8) -- (18,8) -- (18,5) -- (17,5) -- (17,7) -- (16,7) -- (16,5) -- (12,5) -- (12,3) -- (11,3) --(11,5) -- (9,5) -- (9,6);

        \draw [ultra thick, color=blue]
        (23.5,6) -- (21,6) -- (21,5) -- (23,5) -- (23,3) -- (23.5,3);
        \draw [ultra thick, color=blue]
        (-0.5,6) -- (1,6) -- (1,8) -- (2,8) -- (2,6) -- (3,6) -- (3,8) -- (6,8) -- (6,5) -- (5,5) -- (5,7) --(4,7) -- (4,5) -- (0,5) -- (0,3) -- (-0.5,3);

        \draw [ultra thick, color=blue]
        (3,-0.5) -- (3,0) -- (7,0) -- (7,2) -- (8,2) -- (8,0) -- (9,0)-- (9,2) -- (12,2) -- (12, -0.5);

        \foreach \i in {0,1}{
        \draw [ultra thick, color=blue]
        (11+12*\i,-0.5) -- (11+12*\i,1) -- (10+12*\i,1) -- (10+12*\i, -0.5);
        \draw [ultra thick, color=blue]
        (3+12*\i,11.5) -- (3+12*\i,11) -- (5+12*\i,11) -- (5+12*\i, 9) --(6+12*\i,9) -- (6+12*\i,11) -- (10+12*\i,11) -- (10+12*\i,11.5);
        }

        \draw [ultra thick, color=blue]
        (11,11.5) -- (11,11) -- (12,11) -- (12, 11.5);

        \draw [ultra thick, color=blue]
        (15,-0.5) -- (15,0) -- (19,0) -- (19,2) -- (20,2) -- (20,0) -- (21,0)-- (21,2) -- (23.5,2);
        \draw [ultra thick, color=blue]
        (-0.5,2) -- (0,2) -- (0,-0.5);
        \draw [ultra thick, color=blue]
        (-0.5,11) -- (0,11) -- (0,11.5);
        \draw [ultra thick, color=blue]
        (23.5,11) -- (23,11) -- (23,11.5);

        \foreach [evaluate={
              \h=int(Mod(\i,24))+0.5;
              \v=10.5-int(Mod(\i,12));
              }]
                \i in {0,2,...,22}{
                \node [align=center] at (\h, \v) 
                {\color{red}
                \footnotesize $\i,0$};
            }   
            \foreach [evaluate={
              \h=int(Mod(\i+2,24))-0.5;
              \v=11.5-int(Mod(\i,12));
              }]
                \i in {1,3,...,23}{
                \node [align=center] at (\h, \v) 
                {\color{red}
                \footnotesize $\i,0$};
            }

            \foreach [evaluate={
              \h=int(Mod(\i,24))+0.5;
              \v=10.5-int(Mod(\i+4,12));
              }]
                \i in {0,2,...,22}{
                \node [align=center] at (\h, \v) 
                {\color{red}
                \footnotesize $\i,1$};
            }   
            \foreach [evaluate={
              \h=int(Mod(\i+2,24))-0.5;
              \v=11.5-int(Mod(\i+4,12));
              }]
                \i in {1,3,...,23}{
                \node [align=center] at (\h, \v) 
                {\color{red}
                \footnotesize $\i,1$};
            }   
            \foreach [evaluate={
              \h=int(Mod(\i,24))+0.5;
              \v=10.5-int(Mod(\i+8,12));
              }]
                \i in {0,2,...,22}{
                \node [align=center] at (\h, \v) 
                {\color{red}
                \footnotesize $\i,2$};
            }   
            \foreach [evaluate={
              \h=int(Mod(\i+2,24))-0.5;
              \v=11.5-int(Mod(\i+8,12));
              }]
                \i in {1,3,...,23}{
                \node [align=center] at (\h, \v) 
                {\color{red}
                \footnotesize $\i,2$};
            }   

            \foreach [evaluate={
              \h=int(Mod(\i+1,24))+0.5;
              \v=11.5-int(Mod(\i,12));
              }]
                \i in {0,2,...,22}{
                \node [align=center] at (\h, \v) 
                {\color{blue}
                \footnotesize $\i,0$};
            }   
            \foreach [evaluate={
              \h=int(Mod(\i-1,24))+1.5;
              \v=12.5-int(Mod(\i+2,12));
              }]
                \i in {1,3,...,23}{
                \node [align=center] at (\h, \v) 
                {\color{blue}
                \footnotesize $\i,0$};
            }   
            \foreach [evaluate={
              \h=int(Mod(\i+1,24))+0.5;
              \v=11.5-int(Mod(\i+4,12));
              }]
                \i in {0,2,...,22}{
                \node [align=center] at (\h, \v) 
                {\color{blue}
                \footnotesize $\i,1$};
            }   
            \foreach [evaluate={
              \h=int(Mod(\i-1,24))+1.5;
              \v=12.5-int(Mod(\i+6,12));
              }]
                \i in {1,3,...,23}{
                \node [align=center] at (\h, \v) 
                {\color{blue}
                \footnotesize $\i,1$};
            }   
            \foreach [evaluate={
              \h=int(Mod(\i+1,24))+0.5;
              \v=11.5-int(Mod(\i+8,12));
              }]
                \i in {0,2,...,22}{
                \node [align=center] at (\h, \v) 
                {\color{blue}
                \footnotesize $\i,2$};
            }   
            \foreach [evaluate={
              \h=int(Mod(\i-1,24))+1.5;
              \v=12.5-int(Mod(\i+10,12));
              }]
                \i in {1,3,...,23}{
                \node [align=center] at (\h, \v) 
                {\color{blue}
                \footnotesize $\i,2$};
            }

            \foreach [evaluate={
              \h=int(Mod(\i+1,24))+0.5;
              \v=11.5-int(Mod(\i+1,12));
              }]
                \i in {0,1,3,4, 6,7, 9, 10, 12, 13, 15, 16, 18, 19, 21, 22}{
                \node [align=center] at (\h, \v) 
                {\color{black}
                \footnotesize $\i,0$};
            }   
            \foreach [evaluate={
              \h=int(Mod(\i+1,24))+0.5;
              \v=11.5-int(Mod(\i+5,12));
              }]
                \i in {0,1,3,4, 6,7, 9, 10, 12, 13, 15, 16, 18, 19, 21, 22}{
                \node [align=center] at (\h, \v) 
                {\color{black}
                \footnotesize $\i,1$};
            }   
            \foreach [evaluate={
              \h=int(Mod(\i+1,24))+0.5;
              \v=11.5-int(Mod(\i+9,12));
              }]
                \i in {0,1,3,4, 6,7, 9, 10, 12, 13, 15, 16, 18, 19, 21, 22}{
                \node [align=center] at (\h, \v) 
                {\color{black}
                \footnotesize $\i,2$};
            }  
            \foreach [evaluate={
              \h=int(Mod(\i,24))+0.5;
              \v=11.5-int(Mod(\i+2,12));
              }]
                \i in {0,3,...,21}{
                \node [align=center] at (\h, \v) 
                {\color{green}
                \footnotesize $\i,0$};
            }   
            \foreach [evaluate={
              \h=int(Mod(\i,24))+0.5;
              \v=11.5-int(Mod(\i+6,12));
              }]
                \i in {0,3,...,21}{
                \node [align=center] at (\h, \v) 
                {\color{green}
                \footnotesize $\i,1$};
            }

        \foreach \i in {0,1,...,23}{\foreach \j in {0,1,...,11}{\filldraw[fill=white, draw=black] (\i, \j) circle (5pt);}}
    
    \end{tikzpicture}
    \end{center}

\caption{Decomposition of $C_{12} \sq C_{24}$ into $C_{36}$'s.
Since $36= 4\cdot 9$, $k=9$.  So we choose $h=3$, $g=3$. Cycles where the cycle combination operation has been applied are shaded in light gray. The cycles $C_{12}^{\alpha,\beta}$ obtained at the end of Phase 1 ($0 \le \alpha < 8$, $0 \le \beta < 3$) are shown in Figure~\ref{12by24c12}. \newline
For example, applying the cycle combination operation on $T_4^{9,1}$ and $T_4^{12,0}$ combines the red cycles $R_{12}^{2,2}$ ($R_4^{6,2}$, $R_4^{7,2}$, $R_4^{8,2}$), $R_{12}^{3,1}$ ($R_4^{9,1}$, $R_4^{10,1}$, $R_4^{11,1}$), and $R_{12}^{4,0}$ ($R_4^{12,0}$, $R_4^{13,0}$, $R_4^{14,0}$) into a single red cycle $R_{36}^{4,0}$ of length 36 and combines the blue cycles $B_{12}^{2,2}$ ($B_4^{6,2}$, $B_4^{7,2}$, $B_4^{8,2}$), $B_{12}^{3,1}$ ($B_4^{9,1}$, $B_4^{10,1}$, $B_4^{11,1}$), and $B_{12}^{4,0}$ ($B_4^{12,0}$, $B_4^{13,0}$, $B_4^{14,0}$) into a single blue cycle $B_{36}^{4,0}$ of length 36.}\label{12by24c36}
\end{figure}

\begin{figure}

\begin{center}
    \begin{tikzpicture}[scale=.7]

    \foreach [evaluate={
              \h=int(Mod(\i+1,24));
              \v=12-int(Mod(\i+1,12));
              }]
                \i in {0,1,2,4, 5, 6, 8, 9, 10, 12, 13, 14, 16, 17, 18, 20, 21, 22}{
                \fill [gray!15] (\h,\v) rectangle (\h+1,\v-1);;
            }
    \foreach [evaluate={
              \h=int(Mod(\i+1,24));
              \v=12-int(Mod(\i+5,12));
              }]
                \i in {0,1,2,4, 5, 6, 8, 9, 10, 12, 13, 14, 16, 17, 18, 20, 21, 22}{
                \fill [gray!15] (\h,\v) rectangle (\h+1,\v-1);;
            }
    \foreach [evaluate={
              \h=int(Mod(\i+1,24));
              \v=12-int(Mod(\i+9,12));
              }]
                \i in {0,1,2,4, 5, 6, 8, 9, 10, 12, 13, 14, 16, 17, 18, 20, 21, 22}{
                \fill [gray!15] (\h,\v) rectangle (\h+1,\v-1);;
            }   
            
        \foreach \i in {0,1,...,4}{
        \foreach \j in {0,1,2}{
        \draw [ultra thick, dotted, color=red] (0+4*\i,11-4*\j) -- (3+4*\i,11-4*\j) -- (3+4*\i,9-4*\j) -- (5+4*\i,9-4*\j) -- (5+4*\i,8-4*\j) -- (2+4*\i,8-4*\j) -- (2+4*\i,10-4*\j) -- (0+4*\i,10-4*\j) -- (0+4*\i,11-4*\j);
        }}
        \foreach \j in {0,1,2}{
        \draw [ultra thick, dotted, color=red] (23.5,8-4*\j) -- (22,8-4*\j) -- (22,10-4*\j) -- (20,10-4*\j) -- (20,11-4*\j) -- (23,11-4*\j) -- (23,9-4*\j) -- (23.5,9-4*\j);
        \draw [ultra thick, dotted, color=red] (-0.5,9-4*\j) -- (1,9-4*\j) -- (1,8-4*\j) -- (-0.5,8-4*\j);
        }

        \foreach \i in {0,1,...,4}{
        \draw [ultra thick, color=blue] (1+4*\i,11.5) -- (1+4*\i,9) -- (3+4*\i,9) -- (3+4*\i,7) -- (4+4*\i,7) -- (4+4*\i,10) -- (2+4*\i,10) -- (2+4*\i,11.5);
        }

        \foreach \i in {0,1,...,5}{
        \draw [ultra thick, color=blue] (1+4*\i,-0.5) -- (1+4*\i,0) -- (2+4*\i,0) -- (2+4*\i,-0.5);
        }

        \draw [ultra thick, color=blue] (21,11.5) -- (21,9) -- (23,9) -- (23,7) -- (23.5,7);
        \draw [ultra thick, color=blue] (22,11.5) -- (22,10) -- (23.5,10);
        \draw [ultra thick, color=blue] (-0.5,10) -- (0,10) -- (0,7) -- (-0.5,7);

        \foreach \i in {0,1,...,4}{
        \draw [ultra thick, color=blue] (1+4*\i,8) -- (1+4*\i,5) -- (3+4*\i,5) -- (3+4*\i,3) -- (4+4*\i,3) -- (4+4*\i,6) -- (2+4*\i,6) -- (2+4*\i,8) -- (1+4*\i,8);
        }

        \draw [ultra thick, color=blue] (21,8) -- (21,5) -- (23,5) -- (23,3) -- (23.5,3);
        \draw [ultra thick, color=blue] (21,8) -- (22,8) -- (22,6) -- (23.5,6);
        \draw [ultra thick, color=blue] (-0.5,6) -- (0,6) -- (0,3) -- (-0.5,3);

        \foreach \i in {0,1,...,4}{
        \draw [ultra thick, color=blue] (4+4*\i,-0.5) -- (4+4*\i,2) -- (2+4*\i,2) -- (2+4*\i,4) -- (1+4*\i,4) -- (1+4*\i,1) -- (3+4*\i,1) -- (3+4*\i,-0.5);
        }

        \foreach \i in {0,1,...,4}{
        \draw [ultra thick, color=blue] (3+4*\i,11.5) -- (3+4*\i,11) -- (4+4*\i,11) -- (4+4*\i,11.5);
        }

        \draw [ultra thick, color=blue] (23.5,2) -- (22,2) -- (22,4) -- (21,4) -- (21,1) -- (23,1) -- (23, -0.5);

        \draw [ultra thick, color=blue] (-0.5, 2) -- (0,2) -- (0,-0.5);
        \draw [ultra thick, color=blue] (-0.5, 11) -- (0,11) -- (0,11.5);

        \foreach [evaluate={
              \h=int(Mod(\i,24))+0.5;
              \v=10.5-int(Mod(\i,12));
              }]
                \i in {0,2,...,22}{
                \node [align=center] at (\h, \v) 
                {\color{red}
                \footnotesize $\i,0$};
            }   
            \foreach [evaluate={
              \h=int(Mod(\i+2,24))-0.5;
              \v=11.5-int(Mod(\i,12));
              }]
                \i in {1,3,...,23}{
                \node [align=center] at (\h, \v) 
                {\color{red}
                \footnotesize $\i,0$};
            }

            \foreach [evaluate={
              \h=int(Mod(\i,24))+0.5;
              \v=10.5-int(Mod(\i+4,12));
              }]
                \i in {0,2,...,22}{
                \node [align=center] at (\h, \v) 
                {\color{red}
                \footnotesize $\i,1$};
            }   
            \foreach [evaluate={
              \h=int(Mod(\i+2,24))-0.5;
              \v=11.5-int(Mod(\i+4,12));
              }]
                \i in {1,3,...,23}{
                \node [align=center] at (\h, \v) 
                {\color{red}
                \footnotesize $\i,1$};
            }   
            \foreach [evaluate={
              \h=int(Mod(\i,24))+0.5;
              \v=10.5-int(Mod(\i+8,12));
              }]
                \i in {0,2,...,22}{
                \node [align=center] at (\h, \v) 
                {\color{red}
                \footnotesize $\i,2$};
            }   
            \foreach [evaluate={
              \h=int(Mod(\i+2,24))-0.5;
              \v=11.5-int(Mod(\i+8,12));
              }]
                \i in {1,3,...,23}{
                \node [align=center] at (\h, \v) 
                {\color{red}
                \footnotesize $\i,2$};
            }   

            \foreach [evaluate={
              \h=int(Mod(\i+1,24))+0.5;
              \v=11.5-int(Mod(\i,12));
              }]
                \i in {0,2,...,22}{
                \node [align=center] at (\h, \v) 
                {\color{blue}
                \footnotesize $\i,0$};
            }   
            \foreach [evaluate={
              \h=int(Mod(\i-1,24))+1.5;
              \v=12.5-int(Mod(\i+2,12));
              }]
                \i in {1,3,...,23}{
                \node [align=center] at (\h, \v) 
                {\color{blue}
                \footnotesize $\i,0$};
            }   
            \foreach [evaluate={
              \h=int(Mod(\i+1,24))+0.5;
              \v=11.5-int(Mod(\i+4,12));
              }]
                \i in {0,2,...,22}{
                \node [align=center] at (\h, \v) 
                {\color{blue}
                \footnotesize $\i,1$};
            }   
            \foreach [evaluate={
              \h=int(Mod(\i-1,24))+1.5;
              \v=12.5-int(Mod(\i+6,12));
              }]
                \i in {1,3,...,23}{
                \node [align=center] at (\h, \v) 
                {\color{blue}
                \footnotesize $\i,1$};
            }   
            \foreach [evaluate={
              \h=int(Mod(\i+1,24))+0.5;
              \v=11.5-int(Mod(\i+8,12));
              }]
                \i in {0,2,...,22}{
                \node [align=center] at (\h, \v) 
                {\color{blue}
                \footnotesize $\i,2$};
            }   
            \foreach [evaluate={
              \h=int(Mod(\i-1,24))+1.5;
              \v=12.5-int(Mod(\i+10,12));
              }]
                \i in {1,3,...,23}{
                \node [align=center] at (\h, \v) 
                {\color{blue}
                \footnotesize $\i,2$};
            }

            \foreach [evaluate={
              \h=int(Mod(\i+1,24))+0.5;
              \v=11.5-int(Mod(\i+1,12));
              }]
                \i in {0,1,2,4, 5, 6, 8, 9, 10, 12, 13, 14, 16, 17, 18, 20, 21, 22}{
                \node [align=center] at (\h, \v) 
                {\color{black}
                \footnotesize $\i,0$};
            }   
            \foreach [evaluate={
              \h=int(Mod(\i+1,24))+0.5;
              \v=11.5-int(Mod(\i+5,12));
              }]
                \i in {0,1,2,4, 5, 6, 8, 9, 10, 12, 13, 14, 16, 17, 18, 20, 21, 22}{
                \node [align=center] at (\h, \v) 
                {\color{black}
                \footnotesize $\i,1$};
            }   
            \foreach [evaluate={
              \h=int(Mod(\i+1,24))+0.5;
              \v=11.5-int(Mod(\i+9,12));
              }]
                \i in {0,1,2,4, 5, 6, 8, 9, 10, 12, 13, 14, 16, 17, 18, 20, 21, 22}{
                \node [align=center] at (\h, \v) 
                {\color{black}
                \footnotesize $\i,2$};
            }

\junk{            
            \foreach [evaluate={
              \h=int(Mod(\i,24))+0.5;
              \v=11.5-int(Mod(\i+2,12));
              }]
                \i in {0,...,23}{
                \node [align=center] at (\h, \v) 
                {\color{green}
                \footnotesize $\i,0$};
            }   
            \foreach [evaluate={
              \h=int(Mod(\i,24))+0.5;
              \v=11.5-int(Mod(\i+6,12));
              }]
                \i in {0,...,23}{
                \node [align=center] at (\h, \v) 
                {\color{green}
                \footnotesize $\i,1$};
            } 
    }

        \foreach \i in {0,1,...,23}{\foreach \j in {0,1,...,11}{\filldraw[fill=white, draw=black] (\i, \j) circle (5pt);}}
    
    \end{tikzpicture}
    \end{center}

\caption{Decomposition of $C_{12} \sq C_{24}$ into $C_{16}$'s. Since $16=4\cdot 4$, $k=4$, and we can choose $h=4$ and $g=1$. Cycles where the cycle combination operation has been applied are shaded in light gray. \newline
For example, applying the cycle combination operation on $S_4^{0,0}$, $S_4^{1,0}$, and $S_4^{2,0}$ combines the red cycles $R_4^{0,0}$, $R_4^{1,0}$, $R_4^{2,0}$, and $R_4^{3,0}$ into a single red cycle $R_{16}^{0,0}$ of length 16 and the blue cycles $B_4^{0,0}$, $B_4^{1,0}$, $B_4^{2,0}$, and $B_4^{3,0}$ into a single blue cycle $B_{16}^{0,0}$ of length 16. Similarly, applying the cycle combination operation on $S_4^{12,1}$, $S_4^{13,1}$, and $S_4^{14,1}$ combines the red cycles $R_4^{12,1}$, $R_4^{13,1}$, $R_4^{14,1}$, and $R_4^{15,1}$ into a single red cycle $R_{16}^{4,1}$ of length 16 and the blue cycles $B_4^{12,1}$, $B_4^{13,1}$, $B_4^{14,1}$, and $B_4^{15,1}$ into a single blue cycle $B_{16}^{4,1}$ of length 16}\label{12by24c16}
\end{figure}

\begin{figure}

\begin{center}
    \begin{tikzpicture}[scale=.7]

    \foreach [evaluate={
              \h=int(Mod(\i+1,24));
              \v=12-int(Mod(\i+1,12));
              }]
                \i in {0,1,2,4, 5, 6, 8, 9, 10, 12, 13, 14, 16, 17, 18, 20, 21, 22}{
                \fill [gray!15] (\h,\v) rectangle (\h+1,\v-1);
            }
    \foreach [evaluate={
              \h=int(Mod(\i+1,24));
              \v=12-int(Mod(\i+5,12));
              }]
                \i in {0,1,2,4, 5, 6, 8, 9, 10, 12, 13, 14, 16, 17, 18, 20, 21, 22}{
                \fill [gray!15] (\h,\v) rectangle (\h+1,\v-1);
            }
    \foreach [evaluate={
              \h=int(Mod(\i+1,24));
              \v=12-int(Mod(\i+9,12));
              }]
                \i in {0,1,2,4, 5, 6, 8, 9, 10, 12, 13, 14, 16, 17, 18, 20, 21, 22}{
                \fill [gray!15] (\h,\v) rectangle (\h+1,\v-1);
            }   

    \foreach [evaluate={
              \h=int(Mod(\i,24));
              \v=12-int(Mod(\i+6,12));
              }]
                \i in {0,4,...,20}{
                \fill [gray!15] (\h,\v) rectangle (\h+1,\v-1);
            } 
    \foreach [evaluate={
              \h=int(Mod(\i,24));
              \v=12-int(Mod(\i+2,12));
              }]
                \i in {0,4,...,20}{
                \fill [gray!15] (\h,\v) rectangle (\h+1,\v-1);
            }

        \draw [ultra thick, dotted, color=red] (23.5,8) -- (22,8) -- (22,10) -- (21,10) -- (21,8) -- (18,8) -- (18,10) -- (16,10) -- (16,11) -- (19,11) -- (19, 9) -- (20, 9) -- (20,11) -- (23,11) -- (23,9) -- (23.5,9);
        \draw [ultra thick, dotted, color=red] (-0.5,8) -- (1,8) -- (1,10) -- (2,10) -- (2,8) -- (5,8) -- (5,9) -- (3,9) -- (3,11) -- (0,11) -- (0, 9) -- (-0.5, 9);

        \draw [ultra thick, dotted, color=red] (23.5,4) -- (22,4) -- (22,6) -- (20,6) --(20,7) --  (23,7) -- (23,5) --(23.5,5);
        \draw [ultra thick, dotted, color=red] (-0.5,4) -- (1,4) -- (1,6) -- (2,6) --(2,4) --  (5,4) -- (5,6) --(6,6) -- (6,4) -- (9,4) -- (9,5) -- (7,5) -- (7,7) -- (4,7) -- (4,5) -- (3,5) -- (3,7) -- (0,7) -- (0,5) -- (-0.5,5);

        \draw [ultra thick, dotted, color=red] (23.5,0) -- (22,0) -- (22,2) -- (21,2) -- (21,0) -- (18,0) -- (18,2) -- (17,2) -- (17,0) -- (14,0) -- (14,2) -- (12,2) -- (12,3) -- (15,3) -- (15,1) -- (16,1) -- (16,3) -- (19,3) -- (19,1) -- (20,1) -- (20, 3) -- (23,3) -- (23,1) --(23.5,1);
        \draw [ultra thick, dotted, color=red] (-0.5,0) -- (1,0) -- (1,1) -- (-0.5,1);

        \draw [ultra thick, dotted, color=red]
        (21,4) -- (18,4) -- (18,6) -- (17,6) -- (17,4) -- (14,4) -- (14,6) -- (13,6) -- (13,4) -- (10,4) -- (10,6) -- (8,6) -- (8,7) -- (11,7) -- (11,5) -- (12,5) -- (12,7) -- (15,7) -- (15,5) -- (16,5) -- (16,7) -- (19,7) -- (19,5) -- (21,5) -- (21, 4);
        
        \draw [ultra thick, dotted, color=red]
        (17,8) -- (14,8) -- (14,10) -- (13,10) -- (13,8) -- (10,8) -- (10,10) -- (9,10) -- (9,8) -- (6,8) -- (6,10) -- (4,10) -- (4,11) -- (7,11) -- (7,9) -- (8,9) -- (8,11) -- (11,11) -- (11,9) -- (12,9) -- (12,11) -- (15,11) -- (15,9) -- (17,9) -- (17, 8);

        \draw [ultra thick, dotted, color=red]
        (13,0) -- (10,0) -- (10,2) -- (9,2) -- (9,0) -- (6,0) -- (6,2) -- (5,2) -- (5,0) -- (2,0) -- (2,2) -- (0,2) -- (0,3) -- (3,3) -- (3,1) -- (4,1) -- (4,3) -- (7,3) -- (7,1) -- (8,1) -- (8,3) -- (11,3) -- (11,1) -- (13,1) -- (13, 0);

        \draw [ultra thick, color=blue]
        (2,11.5) -- (2,10) -- (4,10) --(4,7) -- (3,7) -- (3,9) -- (0,9) -- (0,7) -- (-0.5,7);
        \draw [ultra thick, color=blue]
        (1,11.5) -- (1,10) -- (-0.5,10);
        \draw [ultra thick, color=blue]
        (23.5,7) -- (23,7) -- (23,9) -- (20,9) -- (20,7) -- (19,7) -- (19,9) -- (17,9) --(17, 11.5);
        \draw [ultra thick, color=blue]
        (18, 11.5) -- (18,10) -- (21,10) -- (21,11.5);
        \draw [ultra thick, color=blue]
        (22, 11.5) -- (22,10) -- (23.5,10);

        \draw [ultra thick, color=blue]
        (5,11.5) -- (5,9) -- (7,9) -- (7,7) -- (8,7) -- (8,9) -- (11,9) -- (11,7) --(12,7) -- (12,9) -- (15,9) -- (15,7) -- (16,7) -- (16,10) -- (14,10) -- (14,11.5);
        \draw [ultra thick, color=blue]
        (6,11.5) -- (6,10) -- (9,10) -- (9,11.5);
        \draw [ultra thick, color=blue]
        (10,11.5) -- (10,10) -- (13,10) -- (13,11.5);

        \draw [ultra thick, color=blue]
        (9,8) -- (9,5) -- (11,5) -- (11,3) -- (12,3) -- (12,5) -- (15,5) -- (15,3) --(16,3) -- (16,5) -- (19,5) -- (19,3) -- (20,3) -- (20,6) -- (18,6) -- (18,8) -- (17,8) -- (17,6) -- (14,6) -- (14,8)  -- (13,8) -- (13,6) -- (10,6) -- (10,8) -- (9,8);

        \draw [ultra thick, color=blue]
        (-0.5, 3) -- (0,3) -- (0,5) -- (3,5) -- (3,3) --(4,3) -- (4,5) -- (7,5) -- (7,3) -- (8,3) -- (8,6) -- (6,6) -- (6,8) -- (5,8) -- (5,6) -- (2,6) -- (2,8)  -- (1,8) -- (1,6) -- (-0.5,6);
        \draw [ultra thick, color=blue]
        (23.5, 3) -- (23,3) -- (23,5) -- (21,5) -- (21,8) -- (22,8) -- (22,6) -- (23.5,6);

        \draw [ultra thick, color=blue]
        (3, -0.5) -- (3,1) -- (1,1) -- (1,4) -- (2,4) -- (2,2) -- (5,2) -- (5,4) -- (6,4) -- (6,2) -- (9,2) -- (9,4) -- (10,4) -- (10,2) -- (12,2) -- (12,-0.5);
        \draw [ultra thick, color=blue]
        (4, -0.5) -- (4,1) -- (7,1) -- (7,-0.5);
        \draw [ultra thick, color=blue]
        (8, -0.5) -- (8,1) -- (11,1) -- (11,-0.5);

        \draw [ultra thick, color=blue]
        (15, -0.5) -- (15,1) -- (13,1) -- (13,4) -- (14,4) -- (14,2) -- (17,2) -- (17,4) -- (18,4) -- (18,2) -- (21,2) -- (21,4) -- (22,4) -- (22,2) -- (23.5,2);
        \draw [ultra thick, color=blue]
        (-0.5,2) -- (0,2) -- (0,-0.5);
        \draw [ultra thick, color=blue]
        (-0.5,11) -- (0,11) -- (0,11.5);
        \draw [ultra thick, color=blue]
        (23.5,11) -- (23,11) -- (23,11.5);
        \draw [ultra thick, color=blue]
        (23,-0.5) -- (23,1) -- (20,1) -- (20, -0.5);
        \draw [ultra thick, color=blue]
        (19,-0.5) -- (19,1) -- (16,1) -- (16, -0.5);

        \foreach \i in {0,1,...,4}{
        \draw [ultra thick, color=blue]
        (3 + 4*\i, 11.5) -- (3+ 4*\i,11) -- (4+ 4*\i,11) -- (4+ 4*\i,11.5);
        }

        \foreach \i in {0,1,...,5}{
        \draw [ultra thick, color=blue]
        (1 + 4*\i, -0.5) -- (1+ 4*\i,0) -- (2+ 4*\i,0) -- (2+ 4*\i,-0.5);
        }

        \foreach [evaluate={
              \h=int(Mod(\i,24))+0.5;
              \v=10.5-int(Mod(\i,12));
              }]
                \i in {0,2,...,22}{
                \node [align=center] at (\h, \v) 
                {\color{red}
                \footnotesize $\i,0$};
            }   
            \foreach [evaluate={
              \h=int(Mod(\i+2,24))-0.5;
              \v=11.5-int(Mod(\i,12));
              }]
                \i in {1,3,...,23}{
                \node [align=center] at (\h, \v) 
                {\color{red}
                \footnotesize $\i,0$};
            }

            \foreach [evaluate={
              \h=int(Mod(\i,24))+0.5;
              \v=10.5-int(Mod(\i+4,12));
              }]
                \i in {0,2,...,22}{
                \node [align=center] at (\h, \v) 
                {\color{red}
                \footnotesize $\i,1$};
            }   
            \foreach [evaluate={
              \h=int(Mod(\i+2,24))-0.5;
              \v=11.5-int(Mod(\i+4,12));
              }]
                \i in {1,3,...,23}{
                \node [align=center] at (\h, \v) 
                {\color{red}
                \footnotesize $\i,1$};
            }   
            \foreach [evaluate={
              \h=int(Mod(\i,24))+0.5;
              \v=10.5-int(Mod(\i+8,12));
              }]
                \i in {0,2,...,22}{
                \node [align=center] at (\h, \v) 
                {\color{red}
                \footnotesize $\i,2$};
            }   
            \foreach [evaluate={
              \h=int(Mod(\i+2,24))-0.5;
              \v=11.5-int(Mod(\i+8,12));
              }]
                \i in {1,3,...,23}{
                \node [align=center] at (\h, \v) 
                {\color{red}
                \footnotesize $\i,2$};
            }   

            \foreach [evaluate={
              \h=int(Mod(\i+1,24))+0.5;
              \v=11.5-int(Mod(\i,12));
              }]
                \i in {0,2,...,22}{
                \node [align=center] at (\h, \v) 
                {\color{blue}
                \footnotesize $\i,0$};
            }   
            \foreach [evaluate={
              \h=int(Mod(\i-1,24))+1.5;
              \v=12.5-int(Mod(\i+2,12));
              }]
                \i in {1,3,...,23}{
                \node [align=center] at (\h, \v) 
                {\color{blue}
                \footnotesize $\i,0$};
            }   
            \foreach [evaluate={
              \h=int(Mod(\i+1,24))+0.5;
              \v=11.5-int(Mod(\i+4,12));
              }]
                \i in {0,2,...,22}{
                \node [align=center] at (\h, \v) 
                {\color{blue}
                \footnotesize $\i,1$};
            }   
            \foreach [evaluate={
              \h=int(Mod(\i-1,24))+1.5;
              \v=12.5-int(Mod(\i+6,12));
              }]
                \i in {1,3,...,23}{
                \node [align=center] at (\h, \v) 
                {\color{blue}
                \footnotesize $\i,1$};
            }   
            \foreach [evaluate={
              \h=int(Mod(\i+1,24))+0.5;
              \v=11.5-int(Mod(\i+8,12));
              }]
                \i in {0,2,...,22}{
                \node [align=center] at (\h, \v) 
                {\color{blue}
                \footnotesize $\i,2$};
            }   
            \foreach [evaluate={
              \h=int(Mod(\i-1,24))+1.5;
              \v=12.5-int(Mod(\i+10,12));
              }]
                \i in {1,3,...,23}{
                \node [align=center] at (\h, \v) 
                {\color{blue}
                \footnotesize $\i,2$};
            }

            \foreach [evaluate={
              \h=int(Mod(\i+1,24))+0.5;
              \v=11.5-int(Mod(\i+1,12));
              }]
                \i in {0,1,2,4, 5, 6, 8, 9, 10, 12, 13, 14, 16, 17, 18, 20, 21, 22}{
                \node [align=center] at (\h, \v) 
                {\color{black}
                \footnotesize $\i,0$};
            }   
            \foreach [evaluate={
              \h=int(Mod(\i+1,24))+0.5;
              \v=11.5-int(Mod(\i+5,12));
              }]
                \i in {0,1,2,4, 5, 6, 8, 9, 10, 12, 13, 14, 16, 17, 18, 20, 21, 22}{
                \node [align=center] at (\h, \v) 
                {\color{black}
                \footnotesize $\i,1$};
            }   
            \foreach [evaluate={
              \h=int(Mod(\i+1,24))+0.5;
              \v=11.5-int(Mod(\i+9,12));
              }]
                \i in {0,1,2,4, 5, 6, 8, 9, 10, 12, 13, 14, 16, 17, 18, 20, 21, 22}{
                \node [align=center] at (\h, \v) 
                {\color{black}
                \footnotesize $\i,2$};
            }  
            \foreach [evaluate={
              \h=int(Mod(\i,24))+0.5;
              \v=11.5-int(Mod(\i+2,12));
              }]
                \i in {0,4,...,20}{
                \node [align=center] at (\h, \v) 
                {\color{green}
                \footnotesize $\i,0$};
            }   
            \foreach [evaluate={
              \h=int(Mod(\i,24))+0.5;
              \v=11.5-int(Mod(\i+6,12));
              }]
                \i in {0,4,...,20}{
                \node [align=center] at (\h, \v) 
                {\color{green}
                \footnotesize $\i,1$};
            }

        \foreach \i in {0,1,...,23}{\foreach \j in {0,1,...,11}{\filldraw[fill=white, draw=black] (\i, \j) circle (5pt);}}
    
    \end{tikzpicture}
    \end{center}

\caption{Decomposition of $C_{12} \sq C_{24}$ into $C_{48}$'s. Since $48= 4\cdot 12$, $k=12$.  So we choose $h=4$, $g=3$. Cycles where the cycle combination operation has been applied are shaded in light gray. The cycles $C_{16}^{\alpha,\beta}$ obtained at the end of Phase 1 ($0 \le \alpha < 6$, $0 \le \beta < 3$) are shown in Figure~\ref{12by24c16}. \newline
For example, applying the cycle combination operation on $T_4^{12,1}$ and $T_4^{16,0}$ combines the red cycles $R_{16}^{2,2}$ ($R_4^{8,2}$, $R_4^{9,2}$, $R_4^{10,2}$, $R_4^{11,2}$), $R_{16}^{3,1}$ ($R_4^{12,1}$, $R_4^{13,1}$, $R_4^{14,1}$, $R_4^{15,1}$), and $R_{16}^{4,0}$ ($R_4^{16,0}$, $R_4^{17,0}$, $R_4^{18,0}$, $R_4^{19,0}$) into a single red cycle $R_{36}^{4,0}$ of length 36 and combines the blue cycles $B_{16}^{2,2}$ ($B_4^{8,2}$, $B_4^{9,2}$, $B_4^{10,2}$, $B_4^{11,2}$), $B_{16}^{3,1}$ ($B_4^{12,1}$, $B_4^{13,1}$, $B_4^{14,1}$, $B_4^{15,1}$), and $B_{16}^{4,0}$ ($B_4^{16,0}$, $B_4^{17,0}$, $B_4^{18,0}$, $B_4^{19,0}$) into a single blue cycle $B_{36}^{4,0}$ of length 36.}\label{12by24c48}
\end{figure}

\section{The wrapping equation}\label{odd}

\begin{figure}
    \begin{center}
        \begin{tikzpicture}[scale=.8, transform shape, square/.style={regular polygon,regular polygon sides=4}]
            \draw [ultra thick, dashed, color=red] (-0.5, 2) -- (0, 2) -- (0, 1) -- (1, 1) -- (1, 2) -- (2,2) -- (2, -0.5);
            \draw [ultra thick, dashed, color=red] (2, 3.5) -- (2,3) -- (3,3) -- (3,1) -- (4,1) -- (4,-0.5);
            \draw [ultra thick, dashed, color=red] (4, 3.5) -- (4,3) -- (5,3) -- (5,2) -- (5.5, 2);
        
            \draw [ultra thick, dotted, color=blue] (-0.5, 3) -- (2, 3) -- (2, 2) -- (5,2) -- (5,1) -- (5.5,1);
            \draw [ultra thick, dotted, color=blue] (-0.5, 1) -- (0, 1) -- (0,0) -- (5,0) -- (5, -0.5);
            \draw [ultra thick, dotted, color=blue] (5, 3.5) -- (5,3) -- (5.5, 3);
        
            \draw [ultra thick, -, color=orange] (0,3.5) -- (0, 2) -- (1,2) -- (1, 3.5);
            \draw [ultra thick, -, color=orange] (1, -0.5) -- (1,1) -- (3,1) -- (3, -0.5);
            \draw [ultra thick, -, color=orange] (3,3.5) -- (3,3) -- (4,3) -- (4,1) -- (5,1) -- (5,0) -- (5.5, 0);
            \draw [ultra thick, -, color=orange] (-0.5, 0) -- (0,0) -- (0,-0.5);
        
            \foreach \i in {0,1,...,5}{\foreach \j in {0,1,2,3}{\filldraw[fill=white, draw=black] (\i, \j) circle (4pt);}}
            \node[square, draw, fill=black] at (0,3){};

            \foreach \y/\i in {3/0,2/1,1/2, 0/3}{\node at (-1,\y){$\i$};}
            \foreach \x in {0,1,...,5}{\node at (\x,4){$\x$};}
            \end{tikzpicture}
\hfill
            \begin{tikzpicture}[scale = .8]
            \node at (10,3) {Cycle};
            \node at (15,3) {$\ell$};
            \node at (13,3) {$v$};
            \node at (14,3) {$h$};
            \node at (16,3) {$\ell_v$};
            \node at (17,3) {$\ell_h$};
            \node at (10,2) {Dotted blue};
            \node at (10,1) {Dashed red};
            \node at (10,0) {Solid orange};
            \node at (15,2) {$0$};
            \node at (13,2) {$1$};
            \node at (14,2) {$2$};
            \node at (16,2) {$0$};
            \node at (17,2) {$0$};
            \node at (15,1) {$1$};
            \node at (13,1) {$2$};
            \node at (14,1) {$1$};
            \node at (16,1) {$1$};
            \node at (17,1) {$0$};
            \node at (13,0) {$1$};
            \node at (14,0) {$1$};
            \node at (15,0) {$3$};
            \node at (16,0) {$3$};
            \node at (17,0) {$0$};
            
            \draw [thick, -] (8.5,2.5) -- (17.5,2.5);
            \draw [thick, -] (12, 3.5) -- (12,-0.5);
        \end{tikzpicture}

\vspace{.2in}
        \begin{tikzpicture}[scale=.8, transform shape, square/.style={regular polygon,regular polygon sides=4}]
            \draw [ultra thick, dashed, color=red] (0, 2) -- (0, 2) -- (0, 1) -- (1, 1) -- (1, 2) -- (2,2) -- (2, -1) -- (3,-1) -- (3,-3) -- (4,-3) -- (4,-5) -- (5,-5) -- (5,-6) -- (6, -6);
        
            \draw [ultra thick, dotted, color=blue] (0, 3) -- (2, 3) -- (2, 2) -- (5,2) -- (5,1) -- (6,1) -- (6, 0) -- (11,0) -- (11, -1)-- (12, -1);
        
            \draw [ultra thick, -, color=orange] (0,3) -- (0, 2) -- (1,2) -- (1, 5) -- (3,5) -- (3,3) -- (4,3) -- (4,1) -- (5,1) -- (5,0) -- (6, 0) -- (6,-1);
        
            \foreach \i in {0,1,...,12}{\foreach \j in {-6,-5,...,5}{\filldraw[fill=white, draw=black] (\i, \j) circle (4pt);}}
            \foreach \x in {0,6,12}{
            \foreach \y in {-1,3,-5}{\node[style=square, draw, fill=black] at (\x,\y){};}}

            \foreach \y/\i in {5/-2,4/-1,3/0, 2/1,1/2,0/3,-1/4,-2/5,-3/6,-4/7, -5/8, -6/9}{\node at (-1,\y){$\i$};}
            \foreach \x in {0,1,...,12}{\node at (\x,6){$\x$};}
            \end{tikzpicture}
    \end{center}
    \caption{Top left: A $C_{16}$-decomposition of $C_4 \sq C_6$, where $(0,0)$ is marked as a black square.  Top right: A table showing the the vertical and horizontal wrapping numbers $v$ and $h$ for each $C_{16}$.  Note that for each cycle, $4v + 6h + 2\ell= 16$ and $\ell = \ell_v + \ell_h$. Bottom: A preimage of each $C_{16}$ in $\mathbb{Z} \sq \mathbb{Z}$ under the mapping $\phi$. Edges not in the preimages are omitted, and the vertices mapped to $(0,0)$ by $\phi$ are marked as black squares.  The preimage of the dotted blue cycle starts at $(0,0)$ and ends at $(4,12)$, so $v=1$ and $h=2$. The preimage of the dashed red cycle starts at $(1,0)$ and ends at $(9,6)$, so $v=2$ and $h=1$. The preimage of the solid orange cycle starts at $(0,0)$ and ends at $(4,6)$, so $v=1$ and $h=1$.
    }\label{wrapping}
\end{figure}
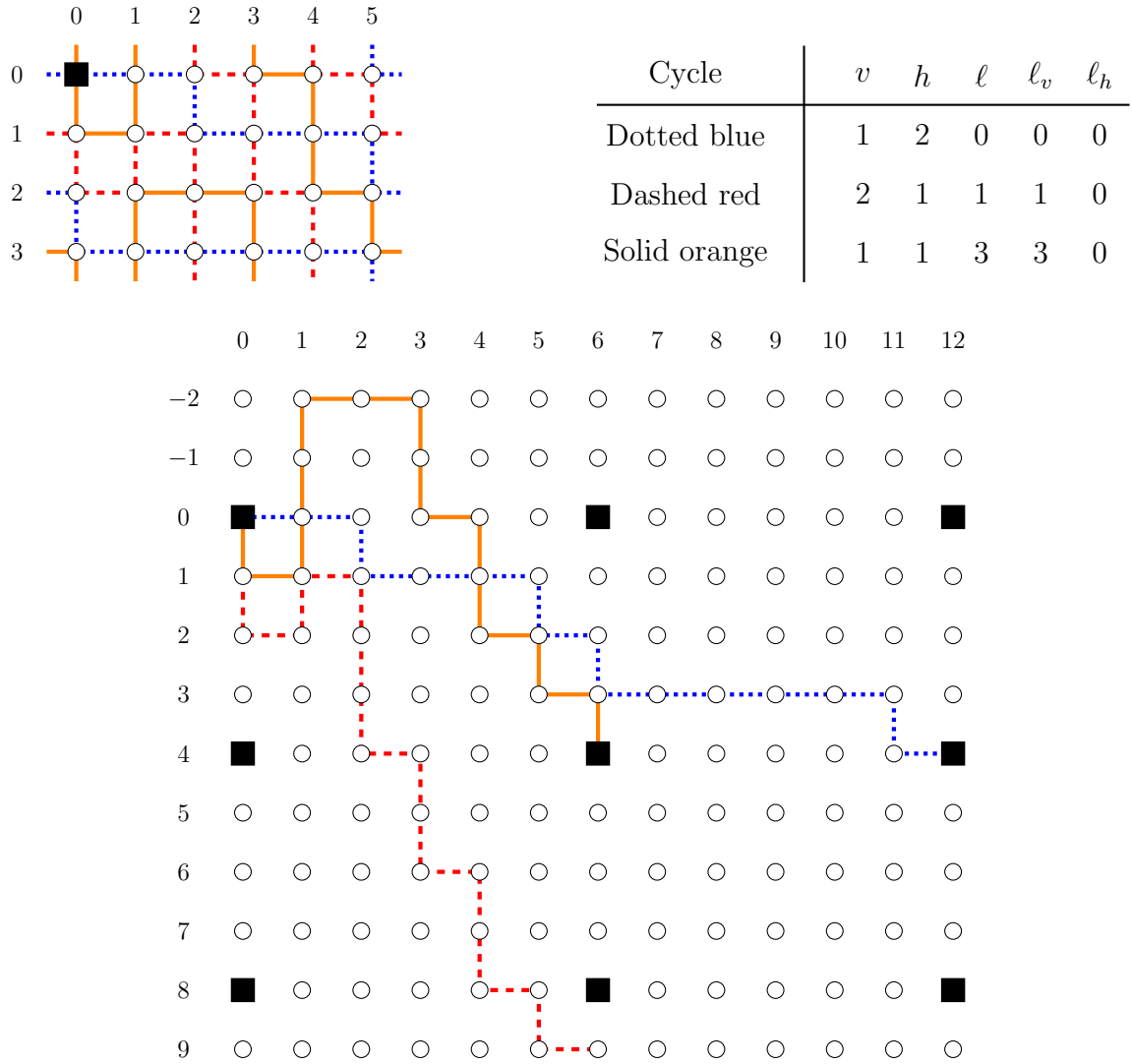

In this section we introduce a condition necessary for $C_k$ to be a subgraph of $C_m \sq C_n$, called the wrapping equation. For each cycle in $C_m \sq C_n$ we define two parameters $h$ and $v$ that we think of as horizontal and vertical wrapping numbers, i.e. $h$ measures how many times the cycle wraps around the torus horizontally, and $v$ measures how many times the cycle wraps around the torus vertically. See Figure~\ref{wrapping}.

To formalize this definition, let $\mathbb{Z}$ denote the graph whose vertices are the integers and whose edges are pairs of consecutive integers. We define a graph homomorphism $\phi$ from the infinite grid $\mathbb{Z} \sq \mathbb{Z}$ to the torus $C_m \sq C_n$ as follows. Define $\phi:\mathbb{Z} \sq \mathbb{Z} \rightarrow C_m \sq C_n$ such that $\phi(i,j) = (i',j')$ if  $i \equiv i' \pmod{m}$ and $j \equiv j' \pmod{n}$, which induces a corresponding mapping of the edges from $E(\mathbb{Z} \sq \mathbb{Z})$ to $E(C_m \sq C_n)$.

To define the wrapping numbers for a given cycle $C_k$ in $C_m \sq C_n$, we look at a preimage $\phi^{-1}(C_k)$ in $\mathbb{Z} \sq \mathbb{Z}$.  Following the cycle $C_k$ in the torus from a specified vertex $w \in V(C_m \sq C_n)$, we eventually return to $w$.  In $\mathbb{Z} \sq \mathbb{Z}$, we start at a vertex $(i_0,j_0)$ where $\phi(i_0,j_0) = w$, and we follow a preimage of the cycle until arriving at another vertex $(i_1,j_1)$ where $\phi(i_1,j_1) = w$.  Let $v = |(i_1-i_0)/m|$ and $h = |(j_1-j_0)/n|$. Further, let $2\ell_v$ be the total number of vertical edges (1st coordinate differs) in the cycle minus $mv$, let $2\ell_h$ be the total number of horizontal edges (second coordinate differs) in the cycle minus $nh$, and $\ell = \ell_h + \ell_v$

%TODO: Describe "following the cycle" more formally?

Since it requires $n$ edges to wrap around the torus horizontally, and $m$ edges to wrap around the torus vertically, a cycle with wrapping numbers $h$ and $v$ must contain at least $nh + mv$ edges.  Further, aside from the edges involved in wrapping, each step in a horizontal or vertical direction must be balanced by a step in the opposite direction, so each of $\ell_v$ and $\ell_h$ must be integers, and $k-(nh + mv)$ must be even.  We now state the wrapping condition for a cycle in a torus:

\textbf{Wrapping equation:} If a cycle $C_k$ is a subgraph of $C_m \sq C_n$, wrapping $v$ times vertically and $h$ times horizontally, then for some nonnegative integer $\ell = \ell_v + \ell_h$, 
\begin{equation}k = nh + mv + 2\ell = nh + mv + 2\ell_v + 2\ell_h.\end{equation}\label{wrap}

The wrapping equation allows us to rule out certain types of cycles from being part of decompositions, and in the following subsections we will use it to restrict the types of cycles we need to consider in a decomposition.

\subsection{$C_6$-decompositions}

\begin{proposition}\label{c6thm}
The cycle $C_6$ decomposes $C_m \sq C_n$ if and only if
\begin{itemize}
\item $m=n=3$, or
\item $m=n=6$, or
\item One of $m$ and $n$ is a multiple of four, and the other is 6.
\end{itemize}
\end{proposition}

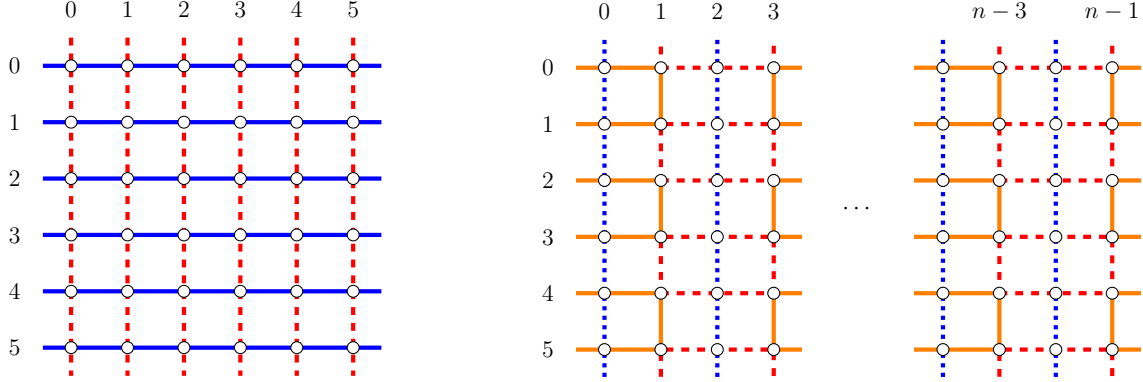
\begin{figure}\begin{center}
\begin{tikzpicture}[scale=.75,transform shape]
    \foreach \i in {0,1,...,5}{\draw[red, dashed, ultra thick] (\i,.5) -- (\i,-5.5);\draw[blue, ultra thick] (-.5,-\i) -- (5.5,-\i);\node at (-1,-\i){\i};\node at (\i,1){\i};}\foreach \x in {0,1,...,5}{\foreach\y in {0,1,...,5}{\filldraw[fill=white, draw=black] (\x,-\y) circle (3pt);}}
\end{tikzpicture}
\hfill
\begin{tikzpicture}[scale=.75, transform shape]
    \begin{scope} 
        \clip(-.5,.5) rectangle (3.5,-5.5);
        \foreach\r in {-1,1,3,5}{\draw[red, dashed, ultra thick] (1,-\r) -- (1,-\r-1) -- (3,-\r-1) -- (3,-\r) -- cycle;}
        \foreach\y in {0,2,4}{\foreach\x in {-1,3}{\draw[orange, ultra thick] (\x,-\y) -- (\x,-\y-1) -- (\x+2,-\y-1) -- (\x+2,-\y) -- cycle;}}
        \foreach \b in {0,2}{\draw[blue, dotted, ultra thick] (\b, 1) -- (\b, -6);}
    \end{scope}
    \node at (4.5,-2.5) {$\cdots$};
    \begin{scope}
        \clip(5.5,.5) rectangle (9.5,-5.5);
        \foreach\r in {-1,1,3,5}{\draw[red, dashed, ultra thick] (7,-\r) -- (7,-\r-1) -- (9,-\r-1) -- (9,-\r) -- cycle;}
        \foreach\y in {0,2,4}{\foreach\x in {5,9}{\draw[orange, ultra thick] (\x,-\y) -- (\x,-\y-1) -- (\x+2,-\y-1) -- (\x+2,-\y) -- cycle;}}
        \foreach \b in {6,8}{\draw[blue, dotted, ultra thick] (\b, 1) -- (\b, -6);}
    \end{scope}
    \foreach \x\n in {0,1,2,3,6/,7/$n-3$,8/, 9/$n-1$}{\node at (\x, 1){\n};}
    \foreach \y in {0,1,...,5}{\node at (-1,-\y){$\y$};\foreach \x in {0,1,2,3,6,7,8,9}{\filldraw[fill=white, draw=black] (\x,-\y) circle (3pt);}}
    \end{tikzpicture}
\end{center}
\caption{Left: A $C_6$-decomposition of  $C_6 \sq C_6$. Right: A $C_6$-decomposition of $C_6 \sq C_n$ for $n$ a multiple of four.}\label{c6_1}
\end{figure}

\begin{figure}
\begin{center}
    \begin{tikzpicture}[scale=.5]
            \draw [ultra thick, -, color=blue] (1,0) -- (1,3);
            \draw [ultra thick, dashed, color=red] (2,1) -- (2,2) -- (0,2) -- (0,1) -- (2,1);
            
            \foreach \i in {0,1,2}{\foreach \j in {0,1,2,3}{\filldraw[fill=white, draw=black] (\i, \j) circle (5pt);}}
        \end{tikzpicture}
    \hspace{.5in}
    \begin{tikzpicture}[scale=.5]
            \draw [ultra thick, -, color=blue] (0,1) -- (3,1);
            \draw [ultra thick, dashed, color=red] (1,2) -- (2,2) -- (2,0) -- (1,0) -- (1,2);
            
            \foreach \i in {0,1,2,3}{\foreach \j in {0,1,2}{\filldraw[fill=white, draw=black] (\i, \j) circle (5pt);}}
        \end{tikzpicture}   
        \hfill
        \begin{tikzpicture}[scale=.5]
            \draw [ultra thick, -, color=blue] (-0.5,2) -- (0,2) --(0,1) -- (3,1) -- (3,2) -- (3.5,2);
            \draw [ultra thick, dashed, color=red] (1,2) -- (2,2) -- (2,0) -- (1,0) -- (1,2);
            
            \foreach \i in {0,1,2,3}{\foreach \j in {0,1,2}{\filldraw[fill=white, draw=black] (\i, \j) circle (5pt);}}
        \end{tikzpicture}
        \hspace{.5in}
        \begin{tikzpicture}[scale=.5]
            \draw [ultra thick, -, color=blue] (-0.5,1) -- (5.5,1);
            \draw [ultra thick, dashed, color=red] (1,2) -- (2,2) -- (2,0) -- (1,0) -- (1,2);
            
            \foreach \i in {0,1,2,3,4,5}{\foreach \j in {0,1,2}{\filldraw[fill=white, draw=black] (\i, \j) circle (5pt);}}
        \end{tikzpicture}
    \end{center}
        \caption{Left: The two possible $C_6$ subgraphs of $C_m \sq C_n$ with $h=v=0$ are shown in dashed red, with one of the form $W(i,j)$ on the left and one of the form $T(i,j)$ on the right.  In either example, if the dashed red $C_6$ is part of a $C_6$-decomposition, then the three solid blue edges must be part of another $C_6$, which is impossible unless $m \in \{4,6\}$ (for $W(i,j)$) or $n\in \{4,6\}$ (for $T(i,j)$). Right: The possible $C_6$ subgraphs that the solid blue edges can be a part of when $n=4$ or $n=6$ (up to symmetry).
    }\label{c6rb}
\end{figure}
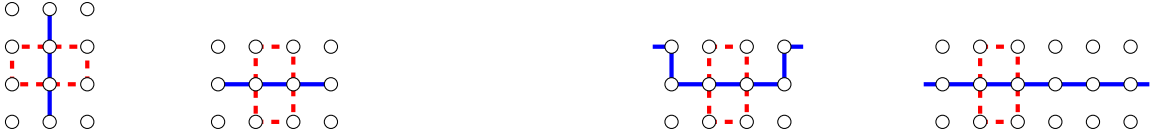

\begin{proof} The wrapping equation implies that a $C_6$ in $C_m \sq C_n$ with $v=h=1$ is only possible when $m=n=3$.
The $C_6$-decomposition of $C_3\sq C_3$ that consists of such cycles is described in Case 1 of Theorem~\ref{3cycles}.  See Figure~\ref{ABcase1}.

To decompose $C_6\sq C_6$ into copies of $C_6,$ we use horizontal and vertical cycles. For $0\leq i \leq 5$, let $H(i)$ be the cycle with vertices $\{(i,0), (i,1), (i,2), (i,3), (i,4), (i,5)\}$ and let $V(j)$ be the cycle with vertices $\{(0,j), (1,j), (2,j), (3,j), (4,j), (5,j)\}$.  The the union $\left(\bigcup_{i=0}^5H(i)\right) \cup \left(\bigcup_{j=0}^5 V(j)\right)$ is a $C_6$-decomposition of $C_6\sq C_6$. See Figure~\ref{c6_1} (left).

Finally, suppose $m=6$ and $n=4N$ for some $N \ge 1$.  Let $W(i,j)$ be the cycle with edges 
\[\{(i,j)(i,j+1), (i,j+1)(i,j+2),(i,j+2)(i+1,j+2), (i+1,j+2)(i+1,j+1),(i+1,j+1)(i+1,j), (i+1,j)(i,j)\}\]
See Figure~\ref{c6rb}.

Then the union
\[\left(\bigcup_{i = 0}^2\bigcup_{j=0}^{N-1} W(2i-1,4j+1) \right)
\cup 
\left(\bigcup_{i = 0}^2\bigcup_{j=0}^{N-1} W(2i,4j-1) \right)
\cup 
\left(\bigcup_{j=0}^{2N-1} V(2j)\right)\]
is a decomposition of $C_m \sq C_n$.  See Figure~\ref{c6_1} (right).

We now show there are no other toruses with $C_6$-decompositions.

Let $T(i,j)$ be the cycle with edges 
\[\{(i,j)(i,j+1), (i,j+1)(i+1,j+1),(i+1,j+1)(i+2,j+1), (i+2,j+1)(i+2,j),(i+2,j)(i+1,j), (i+1,j)(i,j)\}.\]
See Figure~\ref{c6rb}.

Suppose  $m=3$ and $n > 3$. Then no $C_6$ is possible with $v=1$ and $h =0$, since then $mv + nh + 2 \ell$ is odd.  Also, no cycle is possible with $v=1$ and $h > 0$, since then $mv + nh + 2 \ell > 6$.  Thus since no $C_6$ is possible with $v \ge 2$, all $C_6$'s must have $v=0$. The cycles $W(i,j)$ and $T(i,j)$ are the two types of $C_6$'s with $v=h=0$, and as illustrated in Figure~\ref{c6rb} (left), if there is a cycle of type $W(i,j)$, then there is another cycle with three consecutive vertical edges, which is not possible if $m=3$. If there is a cycle of type $T(i,j)$, then $n \in \{4,6\}$ and there is another cycle with $h=1$ containing the three consecutive horizontal edges $\{(i+1,j-1)(i+1,j), (i+1,j)(i+1,j+1), (i+1,j+1)(i+1,j+2)\}$.  Thus the only cycles with $v=h=0$ can be of type $T(i,j)$, and these are only possible if $m=3$ and $n \in \{4,6\}$.  See Figure~\ref{c6rb} (right).  For $C_3 \sq C_4$ and $C_3 \sq C_6$, a quick search shows that the cycles in Figure~\ref{c6rb} (right) can not be part of a full decomposition of these toruses (In $C_3 \sq C_4$, there are no other edge-disjoint $C_6$'s, while in $C_3 \sq C_4$ it is possible to find more edge-disjoint $W(i,j)$'s, but they do not form a decomposition).  Since a decomposition is not possible with only cycles with $h \ge 1$, there is no $C_6$-decomposition of $C_3 \sq C_n$ if $n > 3$.

If $m=6$, $n\ge 5$, and $n$ odd then there are no possible cycles with $h\ge 1$.  The cycles $W(i,j)$ have zero or two horizontal edges in each horizontal cycle, so since there are five edges in each horizontal cycle, it is impossible to cover them using only cycles of the form $W(i,j)$.  Since any cycle of the form $T(i,j)$ would force a cycle with $h\ge 1$, no decomposition is possible.

If $m=6$, $n\ge 10$, and $n \equiv 2 \pmod{4}$ then since a cycle of the form $T(i,j)$ in the decomposition would force a cycle with $h \ge 1$, which is impossible, the only possible cycles are of the form $W(i,j)$ and $V(j)$.  But similar to the checkerboard argument in Proposition~\ref{c4}, if $W(i,j)$ is in a decomposition, then $W(i+2a,j+1)$ and $W(i+4a,j)$ must also be in the decomposition, which inductively implies that $W(i+n-2,j) = W(i-2,j)$ is in the decomposition.  This is a contradiction since $W(i,j)$ and $W(i-2,j)$ share an edge.

Finally, if $m,n \ge 7$, then all $C_6$ subgraphs of $C_m \sq C_n$ are of the form $T(i,j)$ or $W(i,j)$, and no decomposition can be made from just these two types of cycle.
\end{proof}

\subsection{Decompositions of toruses into cycles of odd length}

\begin{proposition}
If $n$ and $m$ are odd and $m \le n < 2m$, then $C_n$ decomposes $C_m \sq C_n$.
\end{proposition}

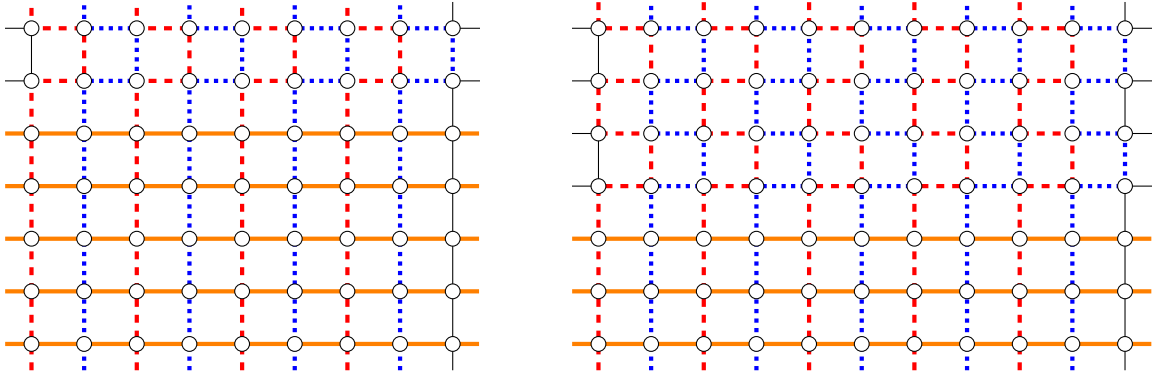
\begin{figure}
\begin{center}
    \begin{tikzpicture}[scale=.7]
        %\small{\node at (3.25,-1) {$C_9$ decomposing $C_7 \sq C_9$.};}
        \begin{scope}
            \clip (-0.5, -0.5) rectangle (8.5,6.5);
            
            \foreach \i in {1, 3, 5, 7}{
                \draw [ultra thick, dotted, color=blue] (\i, -0.5) -- (\i, 5) -- (\i+1, 5) -- (\i+1, 6) -- (\i, 6) -- (\i, 6.5);
            }
            \foreach \i in {-1,8}{
                \draw [thin, color=black] (\i, -0.5) -- (\i, 5) -- (\i+1, 5) -- (\i+1, 6) -- (\i, 6) -- (\i, 6.5);
            }
            \foreach \i in {0, 2, 4, 6}{
                \draw [ultra thick, dashed, color=red] (\i, -0.5) -- (\i, 5) -- (\i+1, 5) -- (\i+1, 6) -- (\i, 6) -- (\i, 6.5);
            }
            \foreach \j in {0,1,...,4}{
                \draw [ultra thick, -, color=orange] (-0.5, \j) -- (8.5, \j);
            }
          
            \foreach \i in {0,1,...,8}{\foreach \j in {0,1,...,6}{
                \filldraw[fill=white, draw=black] (\i, \j) circle (4pt);
            }}
        \end{scope}
        \end{tikzpicture}
        \hfill
        \begin{tikzpicture}[scale=.7]
        %\small{\node at (13.5,-1) {$C_{11}$ decomposing $C_7 \sq C_{11}$.};}
        \begin{scope}
            \clip (9.5,-0.5) rectangle (20.5,6.5);
            
            \foreach \i in {11,13,15,17,19}{
                \draw [ultra thick, dotted, color=blue] (\i, -0.5) -- (\i, 3) -- (\i+1, 3) -- (\i+1,4) -- (\i,4) -- (\i,5) -- (\i+1, 5) -- (\i+1, 6) -- (\i,6) -- (\i,6.5);
            }
            \foreach \i in {10,12, 14, 16, 18}{
                \draw [ultra thick, dashed, color=red] (\i, -0.5) -- (\i, 3) -- (\i+1, 3) -- (\i+1,4) -- (\i,4) -- (\i,5) -- (\i+1, 5) -- (\i+1, 6) -- (\i,6) -- (\i,6.5);
            }
            \foreach \i in {9, 20}{
                \draw [thin, color=black] (\i, -0.5) -- (\i, 3) -- (\i+1, 3) -- (\i+1,4) -- (\i,4) -- (\i,5) -- (\i+1, 5) -- (\i+1, 6) -- (\i,6) -- (\i,6.5);
            }
            
            \foreach \j in {0,1,2}{
                \draw [ultra thick, -, color=orange] (9.5,\j) -- (20.5,\j);
            }
            
            \foreach \i in {10,11,...,20}{\foreach \j in {0,1,...,6}{
                \filldraw[fill=white, draw=black] (\i, \j) circle (4pt);
            }}
        \end{scope}
        \end{tikzpicture}
    \end{center}
    \caption{Left: A $C_9$-decomposition of  $C_7 \sq C_9$. Right: A $C_{11}$-decomposition of $C_7 \sq C_{11}$.
    }\label{oddfig}
\end{figure}

\begin{proof}
See Figure~\ref{oddfig} for two examples.

Let $m= 2M+1$ and $n = 2N+1$.  

For $0 \le \ell < N-M$, define the edge set 
\[E(\ell,j) = \{(2\ell,j)(2\ell,j+1), (2\ell,j+1)(2\ell+1,j+1), (2\ell+1,j+1)(2\ell+1,j), (2\ell+1,j)(2\ell+2,j)\}.
\]

For $1 \le j \le n$, define the cycle $C_n(j)$ to be 
\[C_n(j) = \left(\bigcup_{\ell = 0}^{N-M-1}E(\ell,j)\right) \cup \left(\bigcup_{i=2(N-M)}^{n} (i,j)(i+1,j)\right).\]
These are the cycles  with $v=1$ and $h=0$ in Figure~\ref{oddfig}.

For $1\le i \le m$, let $H_n(i) = \bigcup_{j=1}^n(i,j)(i,j+1)$.  These are the horizontal cycles ($v=0$, $h=1$) in Figure~\ref{oddfig}.

Then the union 
\[\left(\bigcup_{j=1}^n C_n(j) \right)\cup \left(\bigcup_{i=m-n}^m H_n(i)\right)\] 
is a $C_n$-decomposition of $C_m \sq C_n$.
\end{proof}

\begin{proposition}
If $k$ is odd and $m, n > k$, then $C_k$ does not decompose $C_m \sq C_n$.
\end{proposition}
    
\begin{proof}
Since $k<m$ and $k < n$, any cycle must have $h=v=0$.  But since $k$ is odd, $k$ cannot equal $2\ell$. Thus by the wrapping equation, $C_m \sq C_n$ does not contain $C_k$ as a subgraph.
\end{proof}

\section*{Acknowledgements}
This research was conducted in Summer 2022 as a part of the Summer Undergraduate Applied
Mathematics Institute at Carnegie Mellon University,

%%%%%%%%%%%%%%%%%%%%%%%%%%%%%%%%%%%%%%%%%%%%%%
%%%%%%%%%%%%%%%%%%%%%%%%%%%%%%%%%%%%%%%%%%%%%%
%%%%%%%%%%%%%%%%%%%%%%%%%%%%%%%%%%%%%%%%%%%%%%
%%%%%%%%%%%%%%%%%%%%%%%%%%%%%%%%%%%%%%%%%%%%%%

\end{document}